\documentstyle [twoside]{article}
\newtheorem{thm}{Th\'eor\`eme}[section]
\newtheorem{prop}[thm]{Proposition}
\newtheorem{cor}[thm]{Corollaire}
\newtheorem{lem}[thm]{Lemme}
\newtheorem{remark}{Remarque} 
\newtheorem{dfn}{D\'efinition} 

\def\b0{{\bf 0}}
\def\ba{{\bf a}}
\def\bi{{\bf i}}
\def\bI{{\bf J}}
\def\bK{{\bf K}}

\def\bz{{\bf z}}
\def\b1{{\bf 1}}
\def\bA{{\bf A}}
\def\bT{{\bf T}}
\def\bV{{\bf V}}
\def\bP{{\bf P}}
\def\bQ{{\bf Q}}
\def\bR{{\bf R}}
\def\bM{{\bf M}}
\def\bD{{\bf D}}
\def\bN{{\bf N}}
\def\mW{{\mathcal W}}
\def\bZ{{\bf Z}}
\def\bC{{\bf C}}
\def\bzeta{\zeta}

\def\Blat{\mbox{\it \raise2pt\hbox{"}\kern-2pt H}}
\def\lvup{\rlap{\ ${}^{\ell\atop{\hbox{${}^{\vee}$}}}$}\cdots}
\begin{document}
\begin{center}
{\center{\Large{\bf
Transform\'ee  de Mellin des int\'egrales- fibres \\
associ\'ees \`a l'intersection
compl\`ete non-d\'eg\'en\'er\'ee} }}

 \vspace{1pc}
{ \center{\large{ Susumu TANAB\'E }}}
\end{center}

\noindent
\begin{center}
 \begin{minipage}[t]{10.2cm}
{\sc R\'esum\'e.} {\em La transform\'ee de Mellin de  l'int\'egrale
-fibre est calcul\'ee pour certaines classes d'intersection
compl\`ete non-d\'eg\'en\'er\'ee affine, surtout
les cas nomm\'es simpliciables.
On met \`a jour
la structure de r\'eseaux des p\^oles de
la transform\'ee de Mellin \`a l'aide des donn\'ees
topologiques qui d\'ecrivent la structure de Hodge de la vari\'et\'e
affine.
On \'etablit la relation de l'int\'egrale-fibre avec la fonction
hyperg\'eom\'etrique de Horn et de Gel'fand-Kapranov-Zelevinski.
Comme application, nous d\'emontrons une hypoth\`ese propos\'ee par
Berglund et autres sur la sym\'etrie de miroir.}
 \end{minipage} \hfill
\end{center}
 \vspace{1pc}
{
\center{\section{Introduction}}
}

Ici on calcule concr\`etement
l'int\'egrale-fibre associ\'ee \`a la  vari\'et\'e affine
 d'intersection compl\`ete (IC) non-d\'eg\'en\'er\'ee au moyen de sa
 transform\'ee de Mellin.

D'abord nous fixons la situation. Pour les deux vari\'et\'es
complexes $ X =$ soit ${\bf C}^{N}$ soit $({\bf C^\times})^{N}$ et
$S = {\bf C}^k, $ on regarde une application
$$ f: X \rightarrow S  \leqno(0.1)$$
telle que $X_s :=\{(x_1,\cdots, x_{N})\in X ; f_1(x) +s_1=0,
\ldots, f_k(x)+ s_k=0 \}.$ Soient $f_1(x), \cdots,f_k(x)$ des
polyn\^omes qui d\'efinissent  l'intersection compl\`ete,
non-d\'eg\'en\'er\'ee au sens de Danilov-Khovanski \cite{DX1}.
Soit $n$ la dimension de la vari\'et\'e $X_0 =n \geq 0.$
Regardons une hypersurface de dimension
$N+k-1,$ $W_s :=\{(x_1,\cdots, x_{N}, y_1,\cdots, y_k) \in X
\times ({\bf C^\times})^k; y_1(f_1(x) +s_1) + \ldots + y_k( f_k(x) +s_k)
=0 \}.$ Pour la compactification ${\mathcal X}_s$ (resp. ${\mathcal W}_s$)
des vari\'et\'es $X_s$ (resp. $W_s$) dans les vari\'et\'es toriques
proprement construites, les r\'esultats de Terasoma (cas homog\`ene)
\cite{Ter1}, Dimca (cas quasihomog\`ene) \cite{Dim}, Mavlyutov (cas
g\'en\'eral) \cite{Mav} impliquent l'existence d'un isomorphisme entre
les parties primitives de la cohomologie $PH^{N-k}({\mathcal X}_s)$ et
$PH^{N+k-2}({\mathcal W}_s).$
R\'ecemment Terasoma a \'etabli la version affine de cet isomorphisme
pour $X_s$ et  $W_s$ d\'efinis dans les tores \cite{Ter2}, \cite{Ter3}.
L'utilisation effective de cet
isomorphisme (et de ses analogues) s'est baptis\'ee sous le nom de
``technique de Cayley'' (Cayley trick).  Dans ce travail, on regarde
cet isomorphisme entre les structures de Hodge d'une fa\c{c}on un peu
plus d\'etaill\'ee.  C'est \`a  dire, on calcule la transform\'ee de
Mellin des int\'egrales fibres associ\'ees \`a $X_s$ au moyen de
l'isomorphisme en se servant de la technique de Cayley pour la
vari\'et\'e affine.

%
 \vspace{1.5pc}
\footnoterule

\footnotesize{AMS Subject Classification: 14M10, 32S25, 32S40.

Key words and phrases: Gauss-Manin connexion, complete intersection,
Hodge structure.}

 \footnotesize{${}^1$  Travail   r\'ealis\'e   par   le   soutien
 du
Max Planck
Institut f\"ur Mathematik.
}
\normalsize

\newpage

Surtout, dans cette deuxi\`eme partie du travail, \`a part des cas
quasihomog\`enes, nous nous interessons au cas de
la singularit\'e du type d'A'Campo:
$$f_i(x) = x_1^{\alpha_{1i}}+x_2^{\alpha_{2i}}+ \cdots +
x_{N}^{\alpha_{N,i}}\;, 1 \leq i \leq k. \leqno(0.2)$$
Pour le cas $N=2, ({\alpha_{11}},{\alpha_{21}})=(3,2), ({\alpha_{12}},
{\alpha_{22}})=(2,3),$  la monodromie locale
de la singularit\'e $f_1(x)f_2(x)=0$
a \'et\'e \'etudi\'ee par A'Campo dans \cite{AC}.
L'IC de ce type l\`a nous fournit des exemples importants des cas
simpliciables (voir \S 1.4).

Ici nous nous servons de la transform\'ee de Mellin
d'int\'egrales-fibres parce qu'elle permet de mieux visualiser
les propri\'et\'es importantes des vari\'et\'es
d'IC non-d\'eg\'en\'en\'er\'ees par rapport aux autres repr\'esentations
soit d'int\'egrales elles m\^emes soit des \'equations
diff\'erentielles (connexion de Gauss-Manin).  Cette situation a
incit\'e  certains auteurs comme C.Sabbah \cite{Sab2}, \cite{Sab3},
D.Barlet\cite{Bar},
F.Loeser \cite{LS} \`a poursuivre
des recherches sur la  transform\'ee de Mellin d'int\'egrales fibres
qui etaient, entre autres, motiv\'ees par une id\'ee de P.Deligne
reproduite dans \cite{Sab1}.

\par
Dans le \S 1 nous nous souvenons  des faits utiles sur la
technique de Cayley pour l'IC afin de  d\'efinir l'int\'egrale
fibre associ\'ee \`a $f$.  Dans le \S 2, nous interpr\'etons
l'int\'egrale- fibre comme fonction hyperg\'eom\'etrique
g\'en\'eralis\'ee au sens de Mellin-Barnes-Pincherle. Du surcroit,
on montre que l'int\'egrale fibre satisfait le syst\`eme des
\'equations diff\'erentielles du type de Horn.  Pourtant le
degr\'e de chaque \'equation diff\'erentielle est exprim\'e au
moyen du caract\'eristique d'Euler d'une vari\'et\'e d'IC.
 Dans le \S 3, les p\^oles de la transform\'ee de Mellin de
l'int\'egrale- fibre sont exprim\'es  au moyen de la structure de
Hodge mixte d\'ecrite dans \cite{Baty1}.

Dans le \S 4.1, on \'etablit une relation entre les notions des deux
fonctoins hypeg\'eom\'etriques g\'en\'eralis\'ees: celle de Horn et
celle de Gel'fand-Kapranov-Zelevinski. Dans le \S 4.2, 4.3, nous
regardons les applications de notre calcul \`a la sym\'etrie de
miroir: les cas d'IC projective de Givental' et les cas
d'IC \`a des poids multiples quasihomog\`enes. En se servant des r\'esultats du
\S 2, nous d\'emontrons une hypoth\`ese propos\'ee par Berglund, Candelas 
et les autres.

Je tiens \`a remercier T.Terasoma, V.Golyshev de leurs discussions
utiles. J'aimerais aussi mentionner les int\^erets sinc\`eres de
H.Kimura et de Y.Haraoka comme motivation importante pour
l'am\'elioration des \S2, \S 3 et \S 4.1.

 \vspace{2pc}
{
\center{\section{ Technique de Cayley}}
}

{\bf 1.1. Quasihomog\'enisation}
On reprend la situation et les notations de \S 0.
On introduit encore les notations suivantes.
Soit ${\bf T}^m=( {\bf C}\setminus \{0\})^m$ $=(\bC^{\times})^m$
un tore complexe de dimension $m.$ On note par
$x^{\bi}$ le mon\^ome $x^{\bi}:= x_1^{i_1} \cdots x_N^{i_N}$ pour le
multi-indice ${\bi}=({i_1}, \cdots, {i_N}) \in {\bf Z}^{N},$ et par
$dx$ la $N-$forme volume $dx :=dx_1 \wedge\cdots \wedge dx_N.$ Avec la
lettre $\bzeta,$ on note le multi-indice ${\bzeta}= (\zeta_1, \cdots,
\zeta_k) \in ({\bf Z}_{\geq 0})^k.$ Nous nous servons aussi des
notations $x^\b1 :=x_1\cdots x_N, $ $ y^{\bzeta}= y_1^{\zeta_1}\cdots
y_k^{\zeta_k},$ $s^\bz = s_1^{z_1} \cdots s_k^{z_k}$ et $ds =
ds_1\wedge \cdots \wedge ds_k.$ Dans cette section on consid\`ere
l'extension  de l'application $f$ de la vari\'et\'e torique
$\bP_\Sigma$ dans $\bC^k.$ Nous poursuivons la  construction de
\cite{BC} et \cite{Mav}.  On d\'efinit le nombre $M$ comme la dimension
 de l'espace ambiant minimum pour que nous puissons rendre les
polyn\^omes $(f_1(x),\cdots,f_{k}(x))$ simultan\'ement quasihomog\`enes
en multipliant certains termes par des nouvelles variables:  $$ x^\bi
\rightarrow x'_j x^\bi\; j=1,2, \cdots.  $$ Notons par
$(f_1(x,x'),\cdots,f_{k}(x,x'))$ les polyn\^omes ainsi obtenus.  Ils so
nt quasihomog\`enes par rapport \`a certain poids c.\`a.d.
il existe des nombres entiers positifs
$(w_1, \cdots, w_N, w_1', \cdots, w'_{M-N})$ tels que
leur
pgdc soit \'egal \`a 1 et
qu'on ait la relation suivante:
$$ E(x,x')= \sum_{i=1}^N w_ix_i \frac{\partial}{\partial x_i}
+ \sum_{j=1}^{M-N}w'_j x'_j \frac{\partial}{\partial x'_j}
,\leqno(1.1.1)$$
de telle sorte que
$$ E(x,x')(f_\ell(x,x')) = p_\ell f_\ell(x,x')
\;\;\;\mbox{pour}\;\ell=1,\cdots, k,$$
o\`u on appelle $E$ le
champ d'Euler.
\par {\bf Exemple}  On ajoute  au polyn\^ome
$ f(x)=x_1^a + x_1x_2 + x_2^b,$
avec $a,b >2,$ la variable nouvelle $x'_1$ en telle sorte que
le polyn\^ome
$ f(x,x')=x_1^a + x'_1x_1x_2 + x_2^b, $
devienne quasihomog\`ene par rapport au poids $(b,a, ab-a-b)$.

Il y a en g\'en\'eral des choix des termes que l'on modifie pour
r\'ealiser une quasihomog\'en\'eit\'e (voir \S 1.5, 1.6, \S 3.2).
\par D\'esormais nous nous servirons de la notation des variables
$X:= (X_1 , \cdots X_M):= (x_1, \cdots, x_N,$$
x_1', \cdots,$
$ x'_{M-N})$ et celle du polyn\^ome par $f_\ell(X):= f_\ell(x,x').$
Si on introduit un champ d'Euler,
$$ E(X')= \sum_{i=1}^Nw_ix_i \frac{\partial}{\partial x_i}
+ \sum_{j=1}^{M-N}w'_j x'_j \frac{\partial}{\partial x'_j}+
X_{M+1} \frac{\partial}{\partial X_{M+1}},$$
alors on a la relation suivante:
$$ E(X')(f_\ell(X)+X_{M+1}^{p_\ell}s_\ell) =
p_\ell (f_\ell(X)+X_{M+1}^{p_\ell}s_\ell)
\;\;\;\mbox{pour}\;\ell=1,\cdots, k.$$ D'ici bas on note $X':=(X,X_{M+1}).$
 Soit $\bM_\bZ$ un r\'eseau int\'egral de rang $N$ et $\bN_\bZ$
son r\'eseau dual, $\bN_\bZ= Hom(\bM_\bZ, \bZ).$
Nous notons par $\bM_\bR$ (resp. $\bN_\bR$) l'extension
naturelle \`a l'espace r\'eel de $\bM_\bZ$ (resp. $\bN_\bZ$).
Nous consid\'erons la situation pour
$\vec{e}_1, \cdots, \vec{e}_{M+1}$ g\'en\'erateurs des c\^ones unidimensionels
tels que $ \sum_{\ell=1}^{M+1}\bR \vec{e}_\ell = \bN_\bR.$
On peut d\'efinir un \'eventail rationel simplicial $\Sigma$
dans $\bN_\bR$ comme ensemble des c\^ones
simpliciaux engendr\'es par $\vec{e}_1, \cdots, \vec{e}_{M+1}$
ci-dessus.
Notre construction ci-dessus du champ $E(X')$ 
correspond \`a l'exhaussement d'un espace $\bN_\bR \times \bN_\bR'$ 
avec 
une base des g\'en\'erateurs 
$\vec{\tilde e}_{N+1}, \cdots \vec{\tilde e}_{M+1}$ tels que
$$ \sum_{i=1}^Nw_i\vec{\tilde e}_i +  \sum_{j=1}^{M-N}w'_{j} \vec{\tilde 
e}_j +\vec{\tilde e}_{M+1}=0.$$ Ici on a 
$ p_\bN(\vec{\tilde e}_j) =\vec{e}_j $
pour la projection
$p_\bN: \bN_\bR \times \bN_\bR' \rightarrow \bN_\bR.$
Pourtant la dimension de l'espace 
vectoriel 
$\bN_\bR \times \bN_\bR'$ doit \^etre minimale c.\`a .d. 
$dim(\bN_\bR \times \bN_\bR')=$
$M.$

Dans la suite, regardons l'homomorphisme injectif
$$ \varphi: \bM_\bZ \rightarrow \bZ^{M+1},$$
d\'efini par
$$\varphi(\vec{m}) =(<\vec{m}, {\vec{e}_1}>, \cdots,<\vec{m}, {\vec{e}_{M+1}}>).$$
Le conoyau de cette application est un groupe ab\'elien fini
$$ {  C l}(\Sigma) = \bZ^{M+1}/\varphi(\bM_\bZ)$$
dont on peut d\'efinir le groupe (dit le tore de N\'eron-Severi),
$$\bD (\Sigma) := Spec \bC [{  C l}(\Sigma) ].$$
 En plus on introduit l'espace affine,
$$ \bA^{M+1}= Spec \bC [X'].$$
Soit $\hat{X_\sigma} := \prod_{\vec{e}_i \not \in \sigma} X_i,$
et on d\'efinit l'id\'eal
$$B(\Sigma) = <\hat{X_\sigma}: \sigma \in \Sigma> \subset
\bC [X_1, \cdots, X_{M+1}].$$
Soit $ Z(\Sigma):=\bV(B(\Sigma)) \subset \bA^{M+1} $
la vari\'et\'e d\'efinie par $B(\Sigma).$

A l'instar de la m\'ethode initi\'ee par  M.Audin, on construit la vari\'et\'e
torique $\bP_\Sigma$ comme le quotient de $U(\Sigma):=\bA^{M+1} \setminus
Z(\Sigma)$ par l'action du groupe $\bD(\Sigma)$:
$$ \bP_\Sigma = U(\Sigma )/ \bD(\Sigma),$$
avec $ dim \;\bD(\Sigma)= M+1-N, dim\;U(\Sigma )= M+1.$
\par

Gr\^ace \`a la compatibilit\'e de l'action de $\bD(\Sigma)$ sur
$U(\Sigma)$ et celle sur $V_s:=\bV(f_1(X) +X_{M+1}^{p_1}s_1, \cdots, f_k(X)+
X_{M+1}^{p_k}s_k) \subset \bA^{M+1},$
on est ramen\'e \`a d\'efinir le sous-ensemble ferm\'e ${\mathcal X}_s$
de $\bP_\Sigma$
vu que celui-ci est un quotient g\'eom\'etrique,
$$ {\mathcal X}_s := V_s
\cap U(\Sigma)/ \bD (\Sigma).$$
Nous notons aussi $V_{s_i}:=\bV(f_i(X) +X_{M+1}^{p_i}s_i),$
${\mathcal X}_{s_i} := \bV_{s_i}\cap U(\Sigma)/ \bD (\Sigma).$
\begin{dfn} (\cite{Mav})
On dit que $ {\mathcal X}_s$ est une IC
non-d\'eg\'en\'er\'ee si pour chaque $\{i_1, \cdots, i_\ell\} \subset
\{1, \cdots, k\}$ l'ensemble $ {\mathcal X}_{s_{i_1}} \cap \cdots
\cap {\mathcal X}_{s_{i_\ell}} \cap \bT_\tau $ d\'efinit une
sous-vari\'et\'e lisse de codimension $\ell$ pour n'importe quel
$\tau \subset \Sigma,$ o\`u $\bT_\tau = \{X \in \bA^{M+1}: \hat
X_\tau \not =0, X_{j} =0 \;si \; \vec{e}_j \subset \tau\}.$
\label{dfn1}
\end{dfn}
Dor\'enavant on prend pour $s$ telle valeur de $\bC^k$ que sa
fibre ${\mathcal X}_s$ soit non singuli\`ere. On dit que $
{\mathcal X}_s$ est une intersection quasi-lisse si $V_s \cap
U(\Sigma)$ est soit vide soit une sous-vari\'et\'e lisse de
codimension $k$ dans $U(\Sigma).$ Il est connu qu'il existe un
ensemble ferm\'e $D(f)$ de $\bC^k$ tel que $ {\mathcal X}_s$ pour
$s \in \bC^k \setminus D(f)$ est une intersection non-
d\'eg\'en\'er\'ee si $V_0$ d\'efinit la singularit\'e isol\'ee
d'intersection compl\`ete dans $\bA^M = Spec \bC [X_1, \cdots ,
X_M].$
\par
Ici on se souvient de l'espace projectif quasihomog\`ene
$\bP^{(k-1)}(w_1'', \cdots, w_k'') $ \cite{Dolg}, \cite{BC} qui
est d\'efini comme la vari\'et\'e torique associ\'ee aux c\^ones
simpliciaux de l'espace vectoriel r\'eel de dimension $k$
engendr\'e par les vecteurs non-z\'eros $\{\vec\epsilon_1,\cdots,
\vec \epsilon_k \}$ qui satisfont la relation $w_1''\vec\epsilon_1
+\cdots + w_k''\vec\epsilon_k=0.$
 Consid\'erons maintenant un polyn\^ome
$$G(X',s,y):=y_1(f_1(X)+X_{M+1}^{p_1}s_1)+\cdots + y_k (f_k(X)+
X_{M+1}^{p_k}s_k), \leqno(1.1.4)$$ d\'efini sur une vari\'et\'e
torique $\bP (E)$ qui est un espace fibr\'e par
$\bP^{(k-1)}(w_1'', \cdots, w_k'')  $ de la vari\'et\'e
$\bP_{\Sigma}:$ $\bP(E) =\bP_{\Sigma}\times \bP^{(k-1)}(w_1'',
\cdots, w_k'')$ muni d'un bidegr\'e
$$(w(y_1), \cdots, w(y_k))=(1,1,\cdots,1),
(w(X_1), \cdots, w(X_{M+1}))=(0,0,\cdots,0).$$
L'anneau des coordon\'ees homog\`enes au sens de
\cite {Mav}
de $\bP(E)$ est $\bC[X_1, \cdots, X_{M+1}, y_1, \cdots, y_k].$
Cela veut dire qu'il existe le champ d'Euler
$$ E(X',y)=
\sum_{i=1}^Nw_ix_i \frac{\partial}{\partial x_i}
+ \sum_{j=1}^{M-N}w'_j x'_j
\frac{\partial}{\partial x'_j}+
X_{M+1} \frac{\partial}{\partial X_{M+1}}
+\sum_{\ell=1}^k w_\ell'' y_\ell \frac{\partial}{\partial y_\ell},
\leqno(1.1.5)$$
en sorte que
$$ E(X',y) y_b f_b(X) = (w_b'' + p_b)y_b f_b(X),\; 1 \leq b \leq k,$$
$$E(X',y)G(X',s,y) = w(F) G(X',s,y) , $$
 pour $w(G) =w_1'' + p_1 = \cdots = w_k'' + p_k. $
Dans cette situation on peut d\'efinir la vari\'et\'e:
$$ \mW_s := \{(X',y) \in \bP(E): G(X',s,y)=0   \},$$
avec $ dim \; \mW_s = N+k-2.$
D\'esormais on va se servir de la notation $w(\bullet)$ pour le
poids quasihomog\`ene d'un polyn\^ome, d'une forme etc.
determin\'e par le champ d'Euler $E(X', y)$ un objet
quasihomog\`ene $\bullet.$
\par
On peut consid\'erer le technique de Cayley sans r\'ef\'erence au \cite{Mav}.
Notamment on introduit un polyn\^ome, 
$$H(x,y):=y_1f_1(x)+\cdots + y_k f_k(x) \in \bZ[x_1, \cdots, x_n, 
y_1, \cdots, y_k].$$
Soient $ \vec{n}_1, \cdots, \vec{n}_{M+k}$
les \'el\'ements de l'ensemble $ supp( H(x,y)) \subset \bZ^{N+k}.$
On d\'efinit un \'eventail rationel simplicial $\tilde \Sigma$
dans $\bR^{N+k}$ comme ensemble des c\^ones
simpliciaux engendr\'es par $\vec{n}_1, \cdots, \vec{n}_{M+k}.$
Pour ${\tilde \bM}_\bZ = \bM_\bZ \times \bZ^k,$ 
on regarde l'homomorphisme injectif
$$ \varphi: {\tilde \bM}_\bZ \rightarrow \bZ^{M+k},$$
d\'efini par
$$\varphi(\vec{\tilde m}) =(<\vec{\tilde m}, {\vec{n}_1}>, \cdots,
<\vec{\tilde m}, {\vec{n}_{M+k}}>).$$
Le conoyau de cette application est un groupe ab\'elien fini
$$ {  C l}(\tilde \Sigma) = \bZ^{M+k}/\varphi(\tilde \bM_\bZ)$$
dont on peut d\'efinir le groupe 
$$\bD (\tilde \Sigma) := Spec \bC [{  C l}(\tilde \Sigma) ]. \leqno(1.1.6)$$
 On peut construire la vari\'et\'e torique $\bP_{\tilde \Sigma}$
associ\'ee \`a  l'espace affine,
$$ \bA^{M+k}= Spec \bC [X_1, \cdots, X_M, 
y_1, \cdots, y_k ],$$
d'une fa\c{c}on parall\`elle \`a la construction de la vari\'et\'e 
$\bP_{\Sigma}.$

\par
{\bf 1.2. Situation torique}
Dor\'enavant nous notons par $F$ le polyn\^ome suivant,
$$F(X,s,y):=y_1(f_1(X)+s_1)+\cdots + y_k (f_k(X)+s_k), \leqno(1.2.1)$$

Pour la  hypersurface $\mW_s$ on se souvient du fait fondamental
(\cite{DX1} 6.2, \cite{Mav} lemma 4.3) lors d'application du technique de
Cayley.
 \begin{lem}
Si ${\mathcal X}_s$ est une intersection
non-d\'eg\'en\'er\'ee au sens de la D\'efinition 1,
alors l'hypersurface $\mW_s$
est aussi non-d\'eg\'en\'er\'ee.
\end{lem}
Pour $\mW_s$ non-d\'eg\'en\'er\'ee on d\'efinit l'anneau de Jacobi
comme suit:
$$J(G):= < \frac{\partial G(X',s,y)}{\partial X_1}, \cdots,
\frac{\partial G(X',s,y)}{\partial X_M},
\frac{\partial G(X',s,y)}{\partial X_{M+1}},
\frac{\partial G(X',s,y)}{\partial y_1}, \cdots,
\frac{\partial G(X',s,y)}{\partial y_k} >, \leqno(1.2.2)$$
pour $s \in {\bf C}^k \setminus D(f)$ d\'efini ind\'ependamment
du choix concr\`et de  $s.$
 De plus, il est connu
qu'il existe une correspondance entre variation inifinitesimale de
structure de Hodge $ Gr_F^i \;PH^n({\mathcal X}_s)$ et $
Gr_F^{i+k-1} \;PH^{n+2k-2}(\mW_s),$ pour $0 \leq i \leq n.$ Voir
\cite{Ter1} (pour le cas $f$ homog\`ene), \cite {Dim} (cas
quasihomog\`ene) et \cite {Mav} (cas g\'en\'eral).

On note par $PH^{N}(\bP_{\Sigma} \setminus \cup_{i=1}^k {\mathcal
X}_{s_i})$ la partie primitive d\'efinie comme quotient
$$H^{N}(\bP_{\Sigma } \setminus \cup_{i=1}^k {\mathcal X}_{s_i})/
\sqcup_{\ell =1}^k (j_{\ell}^ {\ast}H^{N}(\bP_{\Sigma} \setminus
\cup_{i=1}^k {\mathcal X}_{s_i}) ) $$
  o\`u $j_{\ell}:
\bP_{\Sigma} \setminus \cup_{i=1}^k {\mathcal X}_{s_i} \rightarrow
\bP_{\Sigma}
\setminus \cup_{i\not= \ell} {\mathcal X}_{s_i}. $
Selon \cite{Dim} on trouve l'isomorphisme provoqu\'e par le
r\'esidu de Poincar\'e:
$$ R_{{\mathcal X}_s}: H^{N}(\bP_{\Sigma}
\setminus \cup_{i=1}^k{\mathcal X}_{s_i}) \rightarrow
PH^{N-k}({\mathcal X}_s),$$ dont le noyau est \'egal \`a
$\sqcup_{\ell =1}^k j_{\ell}^{\ast}(H^{N}(\bP_{\Sigma} \setminus
\cup_{i=1}^k {\mathcal X}_{s_i})).$ D'autre part on a la suite
exacte (\cite{BC}, 10.11):
$$ 0 \rightarrow H^{N+k-3}(\bP(E)) \stackrel{[\mW_s]} {\rightarrow}
H^{N+k-1}(\bP(E))\stackrel{i_\Sigma}{\rightarrow}
H^{N+k-1}(\bP(E) \setminus \mW_s) \stackrel{Res}{\rightarrow}
PH^{N+k-2}(\mW_s)  \rightarrow 0.$$
Nous introduisons encore une autre partie primitive
$PH^{N+k-1}(\bP (E) \setminus \mW_s)$
que l'on d\'efinit comme conoyau de l'injection
$$ i_\Sigma:H^{N+k-1}(\bP (E))  \rightarrow H^{N+k-1}(\bP (E) \setminus \mW_s).$$
D'autre part d'apr\`es \cite{Mav} on a l'isomorphisme entre
$PH^{N+k-2}(\mW_s)$ et $PH^{N-k}({\mathcal X}_s).$ Avant de
formuler l'\'enonc\'e central de la section, on introduit quelques
notations associ\'ees \`a un ensemble des indices $I \subset \{1,
\cdots, M+1\}$ avec $|I|= N.$ Notamment,
$$  \hat X_I = \prod _{i \not \in I} X_i,\;\;
dX_I = dX_{i_1} \wedge \cdots \wedge dX_{i_N},\;\; det(e_I) =
det (<\vec{m}_j, \vec{e}_{i_k}>)_{1 \leq j \leq N, 1 \leq i_1 <i_2< \cdots <i_N \leq M+1
} ,$$
o\`u $\vec{m}_1, \cdots , \vec{m}_N$ une base du r\'eseau $\bf M_\bZ.$
On d\'efinit $$\Omega_0(X') = \sum_{|I|=N} det( e_I) {\hat X'_I} dX'_I,$$
qui s'av\`ere le g\'en\'erateur de $H^0\left (\bP_\Sigma,
\Omega^N_{\bP_\Sigma}
(w_1, \cdots, w_N, w_1', \cdots, w'_{M-N})\right ).$
D'ici bas on utilise la notation $X'^{\bI' }y^{\bzeta}=
x_1^{i_1}\cdots x_N^{i_N}  {x'}_{1}^{i'_1}
\cdots {x'}_{M-N+1}^{i'_{M-N+1}}y_1^{\zeta_1}\cdots y_k^{\zeta_k}.$
Alors on obtient la Proposition suivante.
\begin{prop}
(Transformation cohomologique de Radon cf. \cite {Mav},
\cite{Ter2})
\par
Il existe l'isomorphisme entre les parties primitives des
groupes de cohomologie
comme suit:
$$  \leqno(1.2.3) PH^{N}(\bP_{\Sigma} \setminus \cup_{i=1}^k
{\mathcal X}_{s_i}) \rightarrow
PH^{N+k-1}( \bP (E) \setminus \mW_s )$$
$$
\frac{X'^{\bI'} \Omega_0(X')}{ (f_1(X)+X_{M+1}^{p_1}s_1)^{\zeta_1+1}
\cdots (f_k(X) +X_{M+1}^{p_k}s_k)^{\zeta_k+1}}
\;\;\; |\rightarrow
\frac{X'^{\bI' }\omega^{\bzeta}\Omega_0(X') \wedge \Omega_1(\omega)}
{F(X',\omega, s)^{\zeta_1 +\cdots +\zeta_k+k}} ,$$
o\`u
$\Omega_1(\omega)= \sum^k_{\ell=1}(-1)^{\ell+1}
w_\ell''\omega_\ell d\omega_1\wedge \lvup \wedge d\omega_k. $
Et l'\'el\'ement
$ X'^{\bI'} y^\zeta \in \bC [X',y]/J(F),$
avec  $w( X'^{\bI' + \b1} y^{ \zeta + \b1})$ $ = $
$(\zeta_1 +\cdots +\zeta_k+k)w(F). $
\label{prop11}
\end{prop}

{\bf 1.3. Situation affine}
\par Si on restreint la vari\'et\'e torique $\bP(E)$
ci-dessus \`a son intersection avec le tore maximal
$\bT^{N+k-1},$ la partie primitive de la cohomologie
$H^{N+k-2}(\mW_s \cap \bT^{N+k-1})$ peut h\'eriter certaines
propri\'et\'es de $PH^{N+k-2}(\mW_s).$ Notamment il est connu
(\cite{BC}, 11.6) que $W_{N+k-2}H^{N+k-2}(\mW_s \cap \bT^{N+k-1})
\cong PH^{N+k-2}(\mW_s).$ Il est naturel d'esp\'erer que l'on
trouve dans des propri\'et\'es de l'int\'egrale $ I_{X'^{\bI'},
\partial\gamma}^{\bzeta} (s)$ certaine r\'eflexion de celles de la structure
de Hodge $PH^{N+k-2}(\mW_s).$
Dans la suite, nous \'etudions l'int\'egrale fibre d\'efinie sur un 
$(M+1)-$cycle dans $\bT^{M+1}.$
\par
On prend un cycle $\gamma(s) $ $\in H_{M+1-k}(V_s \cap
\bT^{M+1})$ et d\'efinit l'int\'egrale fibre y associ\'ee,
$$ I_{X'^{\bI'}, \partial \gamma}^{\bzeta} (s): =
(\frac{1}{2\pi \sqrt -1})^k  \int_{\partial \gamma(s)}\frac
{X'^{\bI' }dX' }{(f_1(X) +X_{M+1}^{p_1}s_1)^{\zeta_1+1} \cdots
(f_k(X) +X_{M+1}^{p_k}s_k)^{\zeta_k+1}}, \leqno(1.3.1)$$ o\`u
${\partial \gamma(s)}$ $\in H_{M+1}(\bT^{M+1} \setminus {\cup_{i=1}^k
V_{s_i}}\cap \bT^{M+1})$ est un cycle obtenu \`a l'aide de $\partial,$
l'op\'erateur de cobord de Leray. Quant \`a l'op\'eration de
Leray, on renvoie au livre de F.Pham \cite{Ph}, ou bien \`a celui
de V.A.Vassiliev \cite{Vas}. Nous notons par $\psi_{\bI'}$ une
$(M+1-k)$-forme m\'eromorphe telle que
$$ X'^{\bI'} dX' =   df_1\wedge \cdots \wedge df_k
\wedge \psi_{\bI'}. $$ Alors le th\'eor\`eme de r\'esidus de Leray implique
l'\'egalit\'e suivante:
       $$ I_{X'^{\bI'}, \partial\gamma}^{\bzeta} (s) =
\frac{1}{\Gamma(\zeta_1+1) \ldots \Gamma(\zeta_k+1)}
(\frac{\partial}{\partial s})^{\bzeta}  \int_{\gamma(s)}
\psi_{\bI'} \leqno(1.3.2)$$
$$=  \frac{1}{\Gamma(\zeta_1+1) \ldots \Gamma(\zeta_k+1)}
(\frac{\partial}{\partial s})^{\bzeta}
I_{X'^{\bI'}, \partial \gamma}^{0} (s)
$$
o\`u $ (\frac{\partial}{\partial s})^{\bzeta} =
(\frac{\partial}{\partial s_1})^{\zeta_1} \ldots
(\frac{\partial}{\partial s_k})^{\zeta_k}.$

Dans cette situation, on formule la version int\'egrale
de la Proposition ~\ref{prop11}.
\begin{prop}
 Pour un cycle de Leray $\Gamma_s \subset H_{M+1}({\bf T}^{M+1} \setminus
\cup_{i=1}^k V_{s_i} \cap \bT^{M+1})$ tel que
$\Re (f_i(X) +X_{M+1}^{p_i}s_i)|_{\Gamma_s} <0$
hors sa partie compacte,
l'\'egalit\'e suivante a lieu:
$$ \Gamma(\zeta_1+1) \cdots \Gamma(\zeta_k+1) \int_{\Gamma_s}
(f_1(X) + X_{M+1}^{p_1}s_1)^{-\zeta_1-1}
\cdots (f_k(X) + X_{M+1}^{p_k}s_k)^{-\zeta_k-1}\frac{X'^{\bI'+\b1} dX'}
{ {X'}^\b1 }$$
$$
=\Gamma(\zeta_1 + \cdots + \zeta_k+k) \int_{S^{k-1}_+ (w'')\times
\Gamma_s} \frac{X'^{\bI'}\omega^{\bzeta} dX' \wedge
\Omega_0(\omega)} {(\omega_1(f_1(X) + X_{M+1}^{p_1}s_1)+\cdots+
\omega_k(f_k(X) +X_{M+1}^{p_k}s_k))^{\zeta_1 +\cdots
+\zeta_k+k}}.$$ o\`u $S^{k-1}_+(w") =\{(\omega_1, \cdots,
\omega_k): \omega_1^{\frac{\bf w''}{w_1''}}+\cdots
+\omega_k^{\frac{\bf w''}{w_k''}}=1, \omega_\ell >0\;\; pour\;
tout\; \ell, {\bf w''} = \prod_{1 \leq i \leq k}w''_i \}$ et
$\Omega_0(\omega)$ la $(k-1)-$forme de volume sur $S^{k-1}_+
(w''),$
$$\Omega_0(\omega) =\frac{1}{w(F)} \sum_{\ell=1}^k (-1)^\ell w"_\ell \omega_\ell d\omega_1\wedge \lvup \wedge
d\omega_k.$$
\label{prop12}
\end{prop}
{\bf D\'emonstration}
D'abord on se souvient de la relation classique:
$$\int_{\Gamma_s} \int_{\bR_+}
e^{ y_j(f_j(X) + X_{M+1}^{p_j}s_j)} y_j^{\zeta_j} dy_j \wedge dX'
= \Gamma(\zeta_j+1) \int_{\Gamma_s}
 (f_j(X) + X_{M+1}^{p_j}s_j)^{-\zeta_j-1} dX', \leqno(1.3.3)$$
qui se d\'eduit de la d\'efinition de la fonction d\'elta de Dirac comme
r\'esidu. Alors on a  l'\'egalit\'e:
$$
\int_{(\bR_+)^k} \int_{ {\Gamma_s}}
 e^{(y_1(f_1(X)+X_{M+1}^{p_1}s_1)+\cdots
+y_k(f_k(X)+X_{M+1}^{p_k}s_k))}
X'^{\bI'} y^{\zeta} dy \wedge dX' \leqno(1.3.4)$$
$$ =\Gamma(\zeta_1+1)\cdots
\Gamma(\zeta_k+1)\int_{\Gamma_s} X'^{\bI'}
(f_1(X)+X_{M+1}^{p_1}s_1)^{-\zeta_1-1}
\cdots (f_k(X)+X_{M+1}^{p_k}s_k)^{-\zeta_k-1} dX'$$
$$=
\Gamma(\zeta_1+1)\cdots \Gamma(\zeta_k+1)(\frac{1}{2\pi \sqrt
-1})^k I_{X'^{\bI'}, \Gamma_s}^{\bzeta} (s).$$ Si on transforme le
terme gauche de $(1.3.4)$ en effectuant le
 changement des variables $(y_1, \cdots, y_k)
:= (\sigma^{\frac{w_1''}{w(F)}}\omega_1, \cdots,
\sigma^{\frac{w_k''}{w(F)}}\omega_k),$ avec $\sigma \in \bR,$
$(\omega_1, \cdots, \omega_k) \in S^{k-1}_+(w")$ et l'homoth\'etie
naturelle induite $X' \rightarrow (\sigma^{\frac{w_1}{w(F)}}X_1,
\cdots, \sigma^\frac{w_N}{w(F)}X_N,
\sigma^\frac{w_1'}{w(F)}X_{N+1}, \cdots,
\sigma^\frac{w_M'}{w(F)}X_M, \sigma^\frac{1}{w(F)} X_{M+1} ),$ on
aura
 $$ \int_{\Gamma_s} \int_{S^{k-1}_+(w")} \int_{\bR_+}
 e^{\sigma(\omega_1(f_1(X) +X_{M+1}^{p_1}s_1)+\cdots
+\omega_k(f_k(X)+X_{M+1}^{p_k}s_k))}\sigma^{\zeta_1 + \cdots + \zeta_k+k-1}
X'^{\bI'}\omega^{\zeta} d\sigma \wedge  \Omega_0 (\omega)\wedge \ dX'$$
$$= \Gamma(\zeta_1+1)\cdots
\Gamma(\zeta_k+1)\int_{ \Gamma_s}\int_{S^{k-1}_+(w")}
\frac{X'^{\bI'} \omega^{\zeta} \Omega_0(\omega) \wedge dX'}
{\left(\sum_{\ell=1}^k\omega_\ell (f_\ell(X)
+X_{M+1}^{p_\ell}s_\ell) \right)^{\zeta_1 +\cdots \zeta_k +k}}.$$
Ici on se sert des formules $(1.3.3)$, $ dy_1 \wedge \cdots \wedge
dy_k = \sigma^{\frac{w_1''+\cdots +w_k''}{w(F)}-1} d\sigma \wedge
\Omega_0(\omega)$ et de la relation $w(X'^{\bI' + \b1})+
\sum_{\ell=1}^k(\zeta_\ell +1)w''_\ell - (\zeta_1+ \cdots +\zeta_k
+k)w(F)=0.$ {\bf C.Q.F.D.}
\par
On remarque ici que si on applique le changement de variables
$$(X',y) \rightarrow (X_{M+1}^{\frac{w_1}{w(F)}}X_1, \cdots,
X_{M+1}^\frac{w_N}{w(F)}X_N,
X_{M+1}^\frac{w_1'}{w(F)
}X_{N+1}, \cdots,X_{M+1}^\frac{w_M'}{w(F)}X_M,
X_{M+1}, X_{M+1}^{\frac{w_1''}{w(F)}}y_1, \cdots,
X_{M+1}^{\frac{w_k''}{w(F)}} y_k),$$
au terme gauche de (1.3.4), on obtient l'int\'egrale
$$\int_{(\bR_+)^k} \int_{ \Gamma_s}
 e^{X_{M+1}(y_1(f_1(X) +s_1)+\cdots
+y_k(f_k(X) +s_k))}
X^{\bI} y^{\zeta} X_{M+1}^{\zeta_1 +\cdots \zeta_k +k-1}
dy \wedge dX \wedge dX_{M+1}$$
qui est \'egale \`a
$$ I_{X^{\bI}, \gamma}^{\zeta} (s):= \int_{ \gamma_s} X^{\bI +\b1}
(f_1(X) +s_1)^{-\zeta_1-1}
\cdots (f_k(X) +s_k)^{-\zeta_k-1} \frac{dX}{X} \leqno(1.3.5)$$
\`a constante pr\`es.
Ici l'indice d\'enote $\bI =(i_1, \cdots, i_M)$ et le cycle $\gamma_s
\in H^M(\bT^M \setminus {\mathcal X}_s|_{X_{M+1}=1})$
est une projection de $\Gamma_s$ sur $\bT^M).$
On va donc pr\'eferer le travaille avec
l'int\'egrale $(1.3.5)$ \`a celui avec $(1.3.1).$
L'indice $i_{M+1}$ perdue dans l'expression $(1.3.5)$
peut \^etre recuper\'ee par la relation
$i_{M+1}= -1-w(X^{\bI + \b1})
-\sum_{\ell=1}^k(\zeta_\ell +1)w''_\ell +
(\zeta_1+ \cdots +\zeta_k +k)w(F).$

Dans cette situation nous pouvons formuler la version affine de
la Proposition ~\ref{prop11}. On se sert de la notation de la
vari\'et\'e affine
 $$ Z_{F(X,0,y)} =\{(X,y) \in \bT^{M+k}; F(X,0,y)=0 \}.$$
\begin{prop} (\cite{Ter2})
Soit $(f_1(X),  \cdots, f_k(X))$ comme dans \S 1.1. On d\'enote par
$X_q =\{ X \in \bT^M ; f_q(X)= 0 \}, 1 \leq q \leq k$ l'hypersurface
dans $\bT^M.$
Soit $F(X,s,y)$ un polyn\^ome comme $(1.2.1),$ alors on a
l'isomorphisme ci-dessous comme structure de Hodge,
$$\begin {array}{ccc}
 PH^{M+k-1}(Z_{F(X,0,y)} )&\rightarrow &  PH^M( \bT^M\setminus X_1
\cup \cdots \cup X_k) \\
X^\bI y^\zeta \frac{dX\wedge dy}{dF}& \mapsto & \frac{X^\bI
dX}{ f_1(X)^{\zeta_1+1}  \cdots f_k(X)^{\zeta_k+1} } \\
\end {array}  \leqno(1.3.6) $$
o\`u
$$ X^\bI y^\zeta \in \frac{\bC [X,y]}{ \langle
X_1 \frac{\partial F(X,0,y)}{\partial X_1}, \cdots,
X_M \frac{\partial F(X,0,y)}{\partial X_M},
y_1 f_1(X), \cdots, y_k f_k(X)  \rangle}.$$
\label{prop124} \end{prop}

La structure de Hodge mixte de $PH^{M+k-1}(Z_{F(X,0,y)} )$
sera d\'ecrite dans \S 3.1 d'une fa\c{c}on plus d\'etaill\'ee.

En se servant de la Proposition ~\ref{prop12}
 on reduit la transform\'ee de Mellin de
l'int\'egrale $(1.3.1)$ \`a une int\'egrale du type d'Euler:  $$
M_{{\bI},\gamma_s}^\zeta ({\bz} ):=\int_\Pi s^{\bz } I_{X^{\bI},
\gamma_s}^{\zeta} (s ) \frac{ds}{s^{\b1}} \leqno(1.3.7)$$ $$
=\int_{S^{k-1}_+(w'') \times \gamma^\Pi} \frac{X^{\bI}
\omega^\zeta s^{\bz -\b1} dX \wedge \Omega_0(\omega) \wedge ds}
{(\omega_1(f_1(X) +s_1)+\cdots +\omega_k(f_k(X) +s_k))^ {\zeta_1+
\cdots +\zeta_k +k}} $$ $$= \int_{\bR_+} \sigma^{\zeta_1+ \cdots
+\zeta_k +k} \frac{d\sigma}{\sigma} \int_{S^{k-1}_+(w'')}
\omega^\zeta \Omega_0(\omega) \int_{\gamma_s}X^{\bI} dX \int_{\Pi}
s^{\bz-\b1} e^{\sigma(\omega_1(f_1(X) +s_1)+\cdots \omega_k(f_k(X)
+s_k))} ds,$$ ici le cycle $\Pi$ est choisi en sorte qu'il \'evite
les p\^oles de l'int\'egrale $I_{X^{\bI}, \gamma}^{\zeta} (s )$ et
qu'il soit homologue \`a ${\bf R}^k$ (cf. \cite{Tsi2}, (2.1.4)).
On se sert de la notation $\gamma^\Pi:= \cup_{s\in
\Pi}(s,\gamma_s).$ Il ne faut pas le confondre avec l'onglet de
Lefschetz, car $\gamma^\Pi$ est plus t\^ot un tuyau sans onglet.

Et puis,
nous re\'ecrivons, \`a constante pr\`es, la derni\`ere expression
par
$$\int_{(\bR_+)^{k} \times \gamma^\Pi}
e^{\Psi(T)} X^{\bI+\b1} y^{\zeta+\b1}s^{\bz} \frac{dX}{X^\b1}
\wedge \frac{dy}{y^\b1} \wedge \frac{ds}{s^\b1}    \leqno(1.3.8)$$
o\`u
$$\Psi(T) = T_1(X,s,y) + \cdots + T_L(X,s,y)= F(X,s,y),  \leqno(1.3.9)$$
dans lequel chaque terme $T_i(X,s,y)$ repr\'esente un
mon\^ome des variables
$X,s,y$ du polyn\^ome $(1.2.1).$
Nous allons appeler le polyn\^ome (1.1.6) (not\'e par $\Psi(T)$)
{\it fonction de phase}.

{\bf 1.4. Lemmes combinatoires}
Ici on pr\'epare quelques lemmes de caract\`ere combinatoire
qui sont utiles pour le calcul de la transform\'ee de Mellin
$(1.3.7).$

D\'esormais on note par $L$  le nombre des mon\^omes en
$\;(X,s,y)\;$ qui prennent part \`a la fonction de phase. Dans le
cas de la singularit\'e d'A'Campo (0.2), on a $L=(N+1)k,$ avec
$$ T_1 = y_1 x_1^{\alpha_{1,1}},T_2 = y_1 x_2^{\alpha_{2,1}}, \cdots,
T_{L} = y_k s_k. \leqno(1.4.1)$$ D'ici bas pour (0.2),
on suppose
toujours que la matrice $$ \left [\begin {array}{ccc}
\alpha_{1,1}&\cdots&\alpha_{1,N}
\\\noalign{\medskip}
\vdots & \cdots &\vdots\\\noalign{\medskip}
\alpha_{k,1} &\cdots&\alpha_{k,N}
\end {array}\right ],\leqno(1.4.2)$$
est de rang $k.$

Si le nombre $L$ est strictement plus petit que le nombre des variables
$(X,s,y)$ $= M+2k,$ alors nous appelons le cas ``abondant''.
Par contre, si $M+2k <L$ le cas ``d\'eficitaire''.

Dans le cas d\'eficitaire, on peut
d\'efinir le nombre $m$ comme le nombre minimal des nouvelles variables
$$x'' =(x'_1, \cdots, x'_{m}) \leqno(1.4.3)$$ pour
rendre le nombre des variables
\'egal \`a celui des termes selon la
 fa\c{c}on introduite au \S {\bf 1.1}.
Par exemple, la relation (1.4.1) devient comme suit:
$$ T_1 = y_1x_1' x_1^{\alpha_{1,1}},
T_2 = y_1x_2' x_2^{\alpha_{2,1}}, \cdots,
T_{(N+1)k-1} = y_k x'_m x_N^{\alpha_{N,k}},
T_{(N+1)k} = y_k s_k. \leqno(1.4.4)$$
Pour simplifier nous introduisons l'\'ecriture des polyn\^omes
$$ f_i(X)= X^{\vec{\alpha}_{1,i}} +\cdots + X^{\vec{\alpha}_{\tau_i,i}},
1 \leq i \leq k,$$ o\`u $X^{\vec{\alpha}_{j,i}}=
X_1^{\alpha_{j,i,1}}\cdots X_M^{ \alpha_{j,i,M}}=$ $
x_1^{\alpha_{j,i,1}}\cdots x_N^{\alpha_{j,i,N}}
x_1'^{\alpha_{j,i,N+1}}\cdots x_{M-N}'^{\alpha_{j,i,M}}.$ Le
nombre $\tau_i$ d\'enote le nombre des termes qui sont pr\'esents
dans l'expression $f_i(X).$ On a ainsi $ \tau_1 +\cdots +\tau_k+k
=L.$

\begin{prop}
Pour la vari\'et\'e torique $\bP_{\Sigma}$ (resp. $\bP_{\tilde 
\Sigma}$) associ\'ee \`a
${\mathcal X}_s$ d\'efinie par (0.2), (1.4.1), on a $ dim {\bf
D}(\Sigma)$ = $ dim {\bf
D}(\tilde \Sigma)+1 =(N-1)(k-1),$ et $$ m =dim {\bf D}(\Sigma) -1 =
dim {\bf D}(\tilde \Sigma)=
(N-1)(k-1)-1 . \leqno(1.4.5) $$ \label{prop14}
\end{prop}
C'est \`a dire, l'addition convenable (voir Lemme 3.5. ci-dessous)
des nouvelles variables $x' = (x'',x'_{m+1})=
(x'_1, \cdots, x'_{m}, x'_{m+1})$
\`a $f_1(x), \cdots, f_k(x)$  rend notre polyn\^ome $G(X',0,y)$
de $(1.1.4)$
quasihomog\`ene.

{\bf D\'emonstration} D'abord on remarque que $\vec \alpha_{j,i} =
(0, \cdots, 0, \cdots, 0, {\rlap{\ ${}^{j
\atop{\hbox{${}^{\vee}$}}}$}\cdots,}\alpha_{j,i,j},0, \cdots,0),$
$1 \leq  i \leq k, $ $1 \leq  j \leq N, $  selon la notation
ci-dessus. Pour d\'emontrer la proposition, il suffit de trouver
l'ensemble minimal de vecteurs suppl\'ementaires ${\vec e}_{j,i},
(j,i) \in \tilde A,$ $|\tilde A| = m+1$ tel que $ \vec
\alpha_{j,i} -\vec \alpha_{j+1,i} +\vec e_{j,i}  $ poss\`ede un
unique vecteur vertical dans ${\bf C}^{N+m+1}.$ Si la condition
$(1.4.2)$ est satisfaite, on voit que l'espace vectoriel
engendr\'e par
$$ \vec{\alpha}_{j,1} -\vec{\alpha}_{j+1,1}, 1 \leq j \leq N-1$$
$$\vec{\alpha}_{j,i} -\vec{\alpha}_{j+1,i} + \vec{e}_{j,i},  1 \leq j \leq
N-1, 2 \leq i\leq k,$$ avec $\bigwedge_{1 \leq j\leq N-1 , 2 \leq
i\leq k}\vec e_{j,i} \not = 0$ poss\`ede unique vecteur vertical
$$ \bigwedge_{1 \leq j \leq N-1} \big(\vec{\alpha}_{j,1} -\vec{\alpha}_{j+1,1} \big)
\bigwedge_{1 \leq j\leq N-1 , 2 \leq i\leq k}
\big(\vec{\alpha}_{j,i} -\vec{\alpha}_{j+1,i} + \vec{e}_{j,i}
\big)$$ dans ${\bf C}^{N+(N-1)(k-1)}.$ On en d\'eduit que la
dimension de l'espace desir\'e doit \^etre $m+1 = (N-1)(k-1).$
{\bf C.Q.F.D.}

{\bf 1.5}
En suite, on va consid\'erer une param\'etrisation simple de la
vari\'et\'e d\'efinie par (0.1). Pour fixer la situation, d\'esormais on tient uniquement compte
du cas d\'eficitaire des variables. Les cas abondants des variables seront
trait\'es ailleurs (par exemple $(4.2.1), (4.2.2)$ ci-dessous).
On note
$$ \Xi := ^t(x_1, \cdots, x_N,  x_1', \cdots,  x'_m,  s_1, \cdots,
 s_k,y_1, \cdots,   y_k) , \leqno(1.5.1)$$
$$Log\; T := ^t(log\; T_1, \cdots, log\; T_L)
\leqno(1.5.2)$$
$$Log\; \Xi := ^t(\log\; x_1, \cdots, \log\; x_N,   log\; x_1', \cdots, log\; x'_m,
\log\; s_1, \cdots,
log\; s_k,\log\; y_1, \cdots, log\; y_k). \leqno(1.5.3)$$ Ici on
oublie la variable $x'_{m+1}$ en la fixant $x'_{m+1}=1.$ Alors,
une \'equation lin\'eaire \'equivalente \`a $(1.4.4)$ s'\'ecrit
comme,
$$ log\; T_1 = log\; y_1 + log\; x_1' +{\alpha_{1,1}} log\; x_1,
log\; T_2= log\; y_1+ log\; x_2' + {\alpha_{2,1}} log\; x_2, \cdots,
log\; T_{(N+1)k} = log\; y_k + log\; s_k. \leqno(1.5.4)$$

Re\'ecrivons  cette relation lin\'eaire (1.5.4) au moyen d'une matrice
${\sf L} \in End({\bf Z}^L),$
$$ Log\; T= {\sf L}\cdot Log\; X. \leqno(1.5.5)$$
Dor\'enavant, les colonnes des matrices ${\sf L},{\sf L}^{-1} $
etc. seront toujours rang\'ees en correspondance avec l'ordre des variables
donn\'e ci-dessus  dans (1.5.1), (1.5.2), (1.5.3).

Bien entendu, on peut modifier, avec une certaine libert\'e, de
mon\^omes du polyn\^ome $f$ en y ajoutant des param\`etres $(x'_1,
\cdots, x'_{m}),$ comme on verra dans l'exemple  suivant. C'est
\`a dire, on peut ajouter $x'$ comme chez (1.6.3) ou bien (1.6.7)
ci-dessous. On va revenir \`a cette question  plus tard au \S 3.2.
Voir Remarque ~\ref{remark41} \`a propos de la compactification de
Terasoma.

N\'eanmoins, dans le cas de la singularit\'e (0.2), le choix des
mon\^omes auxquels
on multipli\'e les param\`etres suppl\'ementaires est un peu d\'elicat.
Notamm
ent, il faut suivre les r\`egles suivantes
pour \'eviter la d\'eg\'en\'erescence de la matrice $\sf L$
de la relation
(1.5.5).
\begin{lem}
Pour (0.2) et $(1.4.3)$, on a une matrice non-d\'eg\'en\'er\'ee $\sf L$
si on suit les r\`egles suivantes:
\par {\bf a.} Pour l'indice $i \in\{1,\cdots, N\}$
fix\'e, il faut choisir au moins un des mon\^omes
$x_j^{\alpha_{j,i,j}}, 1 \leq i \leq k$ que on ne modifie pas.
\par {\bf b.} Pour l'indice $i \in\{1,\cdots, k\}$
fix\'e, il faut choisir au moins un des termes
$x_j^{\alpha_{j,i,j}}, 1 \leq j \leq N$ que on ne modifie pas.
\label{lem121}
\end{lem}
Par exemple, on peut modifier les termes de $y_\ell(f_\ell(x)+s_\ell)$
 pour $\ell$ assez grand de la  fa\c{c}on suivante:
$$y_\ell(f_\ell(x,x')+s_\ell)= y_\ell
(x_1^{\alpha_{1\ell}}+\sum_{i=2}^N
x'_{i,\ell}x_i^{\alpha_{i\ell}}+s_\ell). \leqno(1.5.6)$$  Ainsi on
cr\'ee un polyn\^ome qui d\'epend des $(N-1)$ nouveaux
param\`etres $(x'_{2,\ell}, \cdots, x'_{N,\ell})$ qui
n'appara\^issent pas aux $f_1(x,x')+s_1, \cdots
f_{\ell-1}(x,x')+s_{\ell-1}.$ On arrive \`a une notion naturelle
qui contient le cas (0.2).
\begin{dfn}
Nous appelons une IC (0.1) simpliciable, s'il existe un polyn\^ome
$(1.1.4)$  tel que $dim {\bf D}(\Sigma)= L-2k -N+1 $ et tel que la
matrice $\sf L$ d\'efinie par $(1.5.5)$ soit
non-d\'eg\'en\'er\'ee. Nous appelons tel polyn\^ome $(1.1.4)$
simplicial.
\end{dfn}
Dans le cas o\`u l'IC est simpliciable, on peut \'enoncer une nouvelle
proposition, analogue \`a la Proposition~\ref{prop14}. C'est \`a dire
une relation entre le nombre des variables n\'ecessaires pour la
quasihomog\'en\'eit\'e et le nombre $m$ de $(1.4.3)$.
\begin{prop}
Pour l'IC simpliciable plong\'ee  dans $\bP_{\Sigma}$ (resp. 
$\bP_{\tilde \Sigma}$),
 on a la
relation suivante;
$$m= dim \; {\bf D}(\Sigma)-1= dim {\bf
D}(\tilde \Sigma)=M-N.\leqno(1.5.7)$$
\end{prop}
Dor\'enavant, on va donc identifier $M= N+m.$

On se souvient ici de la notion de non-d\'eg\'en\'erescence
d'un polyn\^ome qui s'entra\^ine de sa simpliciabilit\'e.

\begin{dfn}
L'hypersurface d\'efinie par un polyn\^ome $g(x)= \sum_{\alpha \in supp(g)}
g_\alpha x^\alpha$ $\in \bC[x_1, \cdots, 
x_n]$ se dit non-d\'eg\'en\'er\'e
si et seulement si pour n'importe quel $\xi \in \bR^n_{\geq 0}$
l'inclusion suivante a lieu,
$$ \{x \in\bC^n;
x_1 \frac{\partial
g^\xi}{\partial x_1}= \cdots=  x_n \frac{\partial g^\xi}{\partial
x_n}=0 \}\subset \{x \in\bC^n;x_1 \cdots x_n =0\}$$
o\`u  $g^\xi(x)= \sum_{\{\beta; <\beta, \xi> \leq <\alpha, \xi>,
\;{\rm pour \;\;tout}\;\;\alpha \in supp(g) \}} g_\alpha x^\alpha.$
\label{dfn15}
\end{dfn}
\begin{prop}
Si la matrice $\sf L$ est non-d\'eg\'en\'er\'ee, l'hypersurface 
$Z_{F(X,0,y)}$ est non-d\'eg\'en\'er\'ee au sens de la D\'efinition
~\ref{dfn15}.
\label{prop18}
\end{prop}

\par {\bf D\'emonstration} D'abord on remarque
que pour $\xi = (\overbrace{0,\cdots,0}^{M},\overbrace{1,\cdots,1}^{k})$
on a $F^\xi(X,0,y)=F(x,0,y).$ Il est n\'ecessaire que 
l'\'equation suivante ait aucune solution dans $\bT^{M+k}$ pour la 
non-d\'eg\'en\'erescence de l'hypersurface $Z_{F(X,0,y)},$
$$ X_1 \frac{\partial F(X,0,y)}{\partial X_1}= \cdots=
X_M \frac{\partial F(X,0,y)}{\partial X_M}=
y_1 f_1(X)= \cdots= y_k f_k(X)=0.$$ 
Cett'\'equation entra\^ine,
$$ T_1(X,y)=T_2(X,y)= \cdots = T_{\tau_1}(X,y)= T_{d^1+1}(X,y),\cdots,
T_{d^2+1}(X,y), \cdots, T_{d^k-1}(X,y)=0, \leqno(1.5.8)$$
avec $d^q:= \sum^q_{i=1}\tau_i+ q.$ Pour resoudre l'\'equation $(1.5.8)$
dans $\bT^{M+k},$ il suffit de chercher la solution non-triviale du syst\`eme,
$$ {\sf L} Log\; \Xi_0= ^t(-\infty, \cdots, -\infty ),$$
$$Log\; \Xi_0 := ^t(\log\; x_1, \cdots, \log\; x_N,   
log\; x_1', \cdots, log\; x'_m, 0, \cdots, 0,\log\; y_1, \cdots, log\; y_k). $$
Si $det {\sf L} \not =0$ une solution quelconque 
$Log\; \Xi_0$ contient un \'el\'ement $-\infty.$
Cela veut dire que l'\'equation $(1.5.8)$ n'a pas de la solution dans
$\bT^{M+k}.$ {\bf C.Q.F.D.}

\par
{\bf 1.6.  Exemple} On consid\`ere un exemple de l'IC
non-d\'eg\'en\'er\'ee d\'efini par 3 polyn\^omes qui d\'ependent
de deux variables,
$$f_i(x) = x_1^{\alpha_{i}}+x_2^{\beta_{i}}\;, 1 \leq i \leq 3. \leqno(1.6.1)$$
Alors la fonction de phase a 9 termes, pourtant elle d\'epend de 8
variables-param\`etres:
$y_1(x_1^{\alpha_{1}}+x_2^{\beta_{1}}+s_1)+
y_2(x_1^{\alpha_{2}}+x_2^{\beta_{2}}+s_2)+
y_3(x_1^{\alpha_{3}}+x_2^{\beta_{3}}+s_3).$

 On se trouve au cas d\'eficitaire.
On ajoute un autre param\`etre $x'_1$
qui met notre polyn\^ome sous forme convenable. Notamment, on regarde
$$y_1(x_1^{\alpha_{1}}+x_2^{\beta_{1}}+s_1)+
y_2(x_1^{\alpha_{2}}+x_2^{\beta_{2}}+s_2)+
y_3(x_1^{\alpha_{3}}+x'_1x_2^{\beta_{3}}+s_3),  \leqno(1.6.2)$$ au
lieu de la fonction de phase initiale.Alors on a
$$ T_1=y_1x_1^{\alpha_{1}}, T_2=y_1x_2^{\beta_{1}}, T_3=y_1s_1,
T_4=y_2x_1^{\alpha_{2}}, T_5= y_2x_2^{\beta_{2}}, T_6=y_2s_2,
T_7=y_3x_1^{\alpha_{3}}, T_8=y_3x'_1x_2^{\beta_{3}}, T_9=s_3.
\leqno(1.6.3)$$ Dans ce cas l\`a, on est ramen\'e \`a resoudre
l'\'equation lin\'eaire suivant:
$$Log\; T
={\sf L_3}\cdot
Log \;X
. \leqno(1.6.4)$$
o\`u
$${\sf L_3}= \left [\begin {array}{ccccccccc} {\it \alpha_1}&0&0&0&0&0&1&0&0
\\\noalign{\medskip}0&{\it \beta_1}&0&0&0&0&1&0&0\\\noalign{\medskip}0&0&0&1
&0&0&1&0&0\\\noalign{\medskip}{\it \alpha_2}&0&0&0&0&0&0&1&0
\\\noalign{\medskip}0&{\it \beta_2}&0&0&0&0&0&1&0\\
\noalign{\medskip}0&0&0&0&1&0&0&1&0\\\noalign{\medskip}{\it
\alpha_3}&0&0&0&0&0&0&0&1
\\\noalign{\medskip}0&{\it \beta_3}&1&0&0&0&0&0&1\\\noalign{\medskip}0&0&0&0
&0&1&0&0&1
\end {array}\right ], \leqno(1.6.5)$$
$$Log \;T= ^t(log\;T_1, log\;T_2,log\;T_3,log\;T_4,log\;T_5,log\;T_6,log\;T_7,log\;T_8,log\;T_9)$$
$$Log\;X=
^t(log\; x_1, log\; x_2, log\; x'_1, log\; s_1,log\; s_2, log\;
s_3, log\; y_1, log\; y_2,log\; y_3).$$ Autrement dit, on a
relation:
$$ Log\; X= {\sf L_3}^{-1} \cdot Log\; T,$$
avec
$$
{\sf L_3}^{-1}=$$
$$\frac{1}{\delta(1,2)}\left [\begin {array}{ccccccccc} {\it \beta_2}&-{\it \beta_2}&0&-{\it \beta_1}&{\it
\beta_1}&0&0&0&0\\
\noalign{\medskip}{\it \alpha_2}&-{\it \alpha_2}&0&-{\it \alpha_1}&
{\it \alpha_1
}&0&0&0&0\\
\noalign{\medskip} -\delta (2,3) &\delta(2,3)&0&\delta(1,3)&
-\delta(1,3)&0&-\delta(1,2)&\delta(1,2)&0\\
\noalign{\medskip}{\it
 \beta_1}\,{\it \alpha_2}&-{\it \alpha_1}\,{\it \beta_2}
&\delta(1,2)&-{\it \alpha_1}\,{\it \beta_1}&{\it \alpha_1}\,{
\it \beta_1}&0&0&0&0\\
\noalign{\medskip}{\it \beta_2}\,{\it \alpha_2}&-{\it \beta_2}\,{
\it \alpha_2}&0&-{\it \alpha_1}\,{\it \beta_2}&{\it \beta_1}\,{\it \alpha_2}&\delta(1,2)&0&0&0\\
\noalign{\medskip}{\it \alpha_3}\,{\it \beta_2}&-{\it
\alpha_3}\,{\it \beta_2}&0&-{\it \beta_1}\,{\it \alpha_3}&{\it \beta_1}\,{\it \alpha_3}&0&-\delta(1,2)&0&\delta(1,2)\\
\noalign{\medskip}-{\it \beta_1}\,{\it \alpha_2}&{\it \alpha_1}\,{\it \beta_2}&0&{\it \alpha_1}
\,{\it \beta_1}&-{\it \alpha_1}\,{\it \beta_1}&0&0&0&0\\
\noalign{\medskip}-{\it \beta_2}\,
{\it \alpha_2}&{\it \beta_2}\,{\it \alpha_2}&0&{\it \alpha_1}\,{\it \beta_2}&-{\it \beta_1}\,{\it \alpha_2}&0
&0&0&0\\
\noalign{\medskip}-{\it \alpha_3}\,{\it \beta_2}&{\it \alpha_3}\,{\it \beta_2}&0&{
\it \beta_1}\,{\it \alpha_3}&-{\it \beta_1}\,{\it \alpha_3}&0&\delta(1,2)&0&0\\
\end {array}\right ]. \leqno(1.6.6)$$
Ici on s'est servi de la notation $\delta(i,j) = \alpha_i\beta_j-
\alpha_j\beta_i, 1 \leq i,j \leq 3.$ Pourtant on note que $det({\sf L_3})=
\delta(1,2).$
Si on regarde une autre param\'etrisation:
$$y_1(x_1^{\alpha_{1}}+x_2^{\beta_{1}}+s_1)+
y_2(x_1^{\alpha_{2}}+{\tilde x'_1}x_2^{\beta_{2}}+s_2)+
y_3(x_1^{\alpha_{3}}+x_2^{\beta_{3}}+s_3). \leqno(1.6.7)$$ On aura
$Log\; {\tilde T}={\sf L_2}\cdot(log\;{\tilde X'})$ avec la
matrice:
$$
{\sf L_2}= \left [\begin {array}{ccccccccc} {\it
\alpha_1}&0&0&0&0&0&1&0&0
\\\noalign{\medskip}0&{\it \beta_1}&0&0&0&0&1&0&0\\\noalign{\medskip}0&0&0&1
&0&0&1&0&0\\\noalign{\medskip}{\it \alpha_2}&0&0&0&0&0&0&1&0
\\\noalign{\medskip}0&{\it \beta_2}&1&0&0&0&0&1&0\\
\noalign{\medskip}0&0&0&0&1&0&0&1&0\\\noalign{\medskip}{\it
\alpha_3}&0&0&0&0&0&0&0&1
\\\noalign{\medskip}0&{\it \beta_3}&0&0&0&0&0&0&1\\\noalign{\medskip}0&0&0&0
&0&1&0&0&1
\end {array}\right ],
 \leqno(1.6.8)$$

$${\sf L_2}^{-1}= \leqno(1.6.9)$$
$$  \frac{1
}{\delta(1,3)}\left [\begin {array}{ccccccccc}
{\it \beta_3}&-{\it \beta_3}&0&0&0&0&-{\it \beta_1}&{\it \beta_1}&0\\
\noalign{\medskip}{\it \alpha_3}&-{\it \alpha_3}&0&0&0&0&-{\it \alpha_1}&{\it \alpha_1}&0\\
\noalign{\medskip}\delta(2,3)&
-\delta(2,3)&0&-\delta(1,3)&\delta(1,3) &0&
\delta(1,2)&-\delta(1,2)&0\\
\noalign{\medskip}{ \it
\beta_1}\,{\it \alpha_3}&-{\it \alpha_1}\,{\it \beta_3}&
\delta(1,3)&0&0&0&-{\it
\beta_1}\,{\it \alpha_1}&{\it \beta_1}\,{\it \alpha_1}&0\\
\noalign{\medskip}{\it \beta_3}\,{\it \alpha_2}&-{\it
\beta_3}\,{\it \alpha_2}&0&-\delta(1,3) &0&\delta(1,3)&-{ \it
\beta_1}\,{\it \alpha_2}&{\it \beta_1}\,{\it
\alpha_2}&0\\\noalign{\medskip}{\it \alpha_3}\,{\it \beta_3}&-{\it
\alpha_3}\,{
\it \beta_3}&0&0&0&0&-{\it \alpha_1}\,{\it \beta_3}&\alpha_3 \beta_1&\delta(1,3)\\
\noalign{\medskip}-{\it \beta_1}\,{\it \alpha_3}&{\it \alpha_1}\,{\it \beta_3}&0&0&0&0&{\it \beta_1}\,{\it \alpha_1}&-
{\it \beta_1}\,{\it \alpha_1}&0\\
\noalign{\medskip}-{\it \beta_3}\,{\it \alpha_2}&{\it \beta_3}\,{\it \alpha_2}&0&\delta(1,3)&0&0&{\it \beta_1}\,{\it \alpha_2}&-{\it \beta_1}\,{\it \alpha_2}&0
\\\noalign{\medskip}-{\it \alpha_3}\,{\it \beta_3}&{\it \alpha_3}\,{\it \beta_3}&0&0&0&0&{\it \alpha_1}\,{\it \beta_3}&-{
\it \beta_1}\,{\it \alpha_3}&0\\
\end {array}\right ]$$

$${\sf L_3}\cdot{\sf C_{32}}={\sf L_2} \leqno(1.6.10)$$
avec $\sf C_{32}= -\frac{\delta(1,3) }{\delta(1,2)}.$ On a
$det({\sf L_2}) = -\delta(1,3).$

{
\center{\section{ Repr\'esentation Mellin-Barnes-Pincherle de
l'int\'egrale fibre}}
}

{\bf 2.1  L'int\'egrale- fibre en tant qu'une fonction hyperg\'eom\'etrique
g\'en\'eralis\'ee }

Nous notons par $X^{\bI}$ le mon\^ome $x_1^{i_1}\cdots x_N^{i_N}
{x'}_{1}^{i'_1}
\cdots {x'}_{m}^{i_m}.$
Dans le cas de l'IC simpliciable, on a  $m =M-N.$
On reprend la forme

$$
M_{{\bI},\gamma_s}^\zeta ({\bz} ):=\int_\Pi s^{\bz } I_{X^{\bI},
\gamma_s}^{\zeta} (s ) \frac{ds}{s^{\b1}},
\leqno(2.1.1)$$
pour certain cycle $\Pi$ qui \'evite les p\^oles de
$I_{X^{\bI}, \gamma}^{\bzeta} (s).$
D'ailleurs il serait utile de voir le calcul du chapitre pr\'ec\'edent,
en liaison avec la notion  des fonctions
hyperg\'eom\'etriques g\'en\'eralis\'ees (FHG) au sens de
Mellin-Barnes-Pincherle
\cite{AK}, \cite{Nor}.
Par cette formulation la FHG de Gauss s'exprime par l'int\'egrale,
$$ _2 F_1(\alpha,\beta,\gamma|s)=
\frac {1}{2\pi i } \int_{z_0 - i\infty}^{  z_0 + i\infty }(-s)^z
\frac {\Gamma ( z+ \alpha)\Gamma ( z+ \beta )\Gamma ( -z)}{\Gamma
( z+ \gamma)} dz , \;\; - \Re  \alpha, - \Re  \beta < z_0.$$
L'ensemble des indices des colonnes et des rayons de la matrice
$\sf L$ sera not\'e par $I,$
$$ I:= \{1, \cdots, L\}.$$ Ici on se souvient de la relation
$L=N+m+2k= M+2k.$
La notion suivante facilitera notre formulation des r\'esultats.
\begin{dfn}
Nous appelons une fonction $g(z)= g(z_1, \cdots, z_k)$
$\Delta-$ p\'eriodique pour
$\Delta \in \bZ_{>0},$ si
$$g(z)= h(e^{\pi \sqrt -1
\frac{z_1}{\Delta}}, \cdots,
e^{\pi \sqrt -1
\frac{z_k}{\Delta}}),$$
pour une fonction rationelle $h(\zeta_1, \cdots, \zeta_k)$
quelconque.
\end{dfn}

Pour (0.1) une IC simpliciable, on a l'\'enonc\'e comme suit.
\begin{prop}
1)Il existe un cycle $\Pi \in H_{k}({\bf T}^{k} \setminus S.S.
I_{X^{\bI}, \gamma}^{\bzeta} (s ))$ tel que
la transform\'ee de Mellin $(2.1.1)$
se repr\'esente comme produit des $\Gamma-$ fonctions
\`a un facteur $g(z)$ de fonction $\Delta-$ p\'eriodique pr\`es:
$$  M_{{\bI}, \gamma}^\zeta (\bz )= g(z) \prod_{a \in I}
\Gamma\bigl( {\mathcal L}_a({\bI, \bz, \bzeta})\bigr),$$ avec
$${\mathcal L}_a({\bI, \bz, \bzeta} ) =\frac{\sum_{j=1}^N A_j^a (i_j+1)
+\sum_{j=1}^{m} C_j^a(i'_j+1)
+\sum_{\ell=1}^k \left( B_\ell^a z_\ell + D_\ell^a(\zeta_\ell+1)\right)
}{\Delta}, a \in I. \leqno(2.1.2)$$
Ici on a la matrice aux \'el\'ements entiers,
$$^t{\sf T}=(A_1^a,\cdots,
A_{N}^a, C_1^a , \cdots, C_{m}^a, B_1^a
, \cdots, B_k^a ,D_1^a
, \cdots, D_k^a )_{1 \leq
a \leq L}, \leqno(2.1.3)$$
avec $p.g.d.c.(A_1^a,\cdots,
A_{N}^a, C_1^a , \cdots, C_{m}^a, B_1^a
, \cdots, B_k^a ,D_1^a
, \cdots, D_k^a )=1,$
pour tout $ a \in [1, L],$
o\`u la matrice $\Delta^{-1}{\sf T}$
est l'inverse de $ \sf L.$
Ainsi $\Delta >0$ est determin\'e d'une mani\`ere unique.
Pourtant les coefficients (2.1.2) satisfont les
propri\'et\'es suivantes
pour chaque indice $a \in I$ :

$\bf a$
Soit
${\mathcal L}_a({\bI, \bz, \bzeta} ) = \frac{\Delta}{\Delta}z_\ell,$
c.\`a.d. $A_1^a=\cdots=
A_{N}^a=0,$ $B_1^a = \cdots\lvup = B_k^a=0,$ $B_\ell^a=1.$

$\bf b$
Soit
$${\mathcal L}_a({\bI, \bz, \bzeta} )=
\frac{\sum_{j=1}^N A_j^a (i_j+1)+\sum_{j=1}^{m} C_j^a(i'_j+1)
+\sum_{\ell=1}^k B_\ell^a (z_\ell -\zeta_\ell-1) }{\Delta}$$

2) Pour chaque l'indice $1 \leq \ell \leq N, 1 \leq q \leq k,$
$1 \leq j \leq m$
fix\'e les \'egalit\'es suivantes ont lieu:
$$\sum_{a \in I} A_\ell^a =0,\; \sum_{a \in I}
B_q^a =0,\sum_{a \in I}
C_j^a =0.\leqno(2.1.4)$$

3)Parmi les fonction lin\'eaires ${\mathcal L}_a,$ $a \in I$
la relation suivante a lieu:
$$ \sum_{a \in I }{\mathcal L}_a(\bI, \bz, \zeta)= \zeta_1 +\cdots +
\zeta_k +k.$$
\label{prop21}
\end{prop}

{\bf D\'emonstration}

1) D'abord
on se souvient de la d\'efinition de la $\Gamma-$fonction,
$$ \int_{C_a} e^{-T_a}T_a^{\sigma_a} \frac{dT_a}{T_a} = (1-
e^{2 \pi i \sigma_a})
\Gamma(\sigma_a),$$
pour l'unique cycle $C_a$ nontrivial
qui tourne autour de $T_a=0$ avec les asymptotes
$\Re T_a \rightarrow + \infty.$
On transforme l'int\'egrale
$(1.3.8)$ \`a l'aide de l'application $(1.5.5),$
$$
\int_{({\bR_+})^k \times  \gamma^\Pi } e^{\Psi(T(X,s,y))}
X^{\bI+\b1 } s^\bz y^{{\bzeta}+\b1} \frac{dX}{X^\b1} \wedge
\frac{dy}{y^\b1} \wedge \frac{ds}{s^\b1} \leqno(2.1.5)$$
$$
= (det {\sf L})^{-1}\int_{{\sl L}_\ast ({ \bR_+}^k
\times \gamma^\Pi  )}
e^{\sum_{a \in I}T_a} \prod_{a\in I}
T_a^{{\mathcal L}_a({\bI, \bz, \bzeta})}
\bigwedge_{a \in I} \frac{dT_a}{T_a}.$$
Evidemment l'int\'egrale $(2.1.5)$
est \'egale au $(1.3.8)$ \`a des fonctions $\Delta-$
p\'eriodiques pr\`es.
Ici ${\sl L}_\ast ({\bR_+}^k
\times \gamma^\Pi)$ d\'enote l'application d'un cha\^ine dans
$\bC^{M}_X\times \bC^{k}_s\times \bC^{k}_y$
au celui dans $\bC^{M+2k}_T$
induite par la transformation $(1.5.5).$
On peut consid\'erer une action naturelle
$\lambda:C_a \rightarrow \lambda(C_a)$
d\'efinie par la relation,
$$ \int_{\lambda(C_a)} e^{-T_a}T_a^{\sigma_a} \frac{dT_a}{T_a}
= \int_{(C_a)} e^{-T_a}(e^{2\pi  \sqrt -1 }T_a)^{\sigma_a} \frac{dT_a}{T_a}.$$
A l'aide de cett'action le cha\^ine ${\sf L}_\ast (
{\bR_+}^k \times
\gamma^\Pi ) $ s'av\`ere homologue au cha\^ine,
$$\sum_{(j_1^{(\rho)},\cdots,j_L^{(\rho)} )\in [1,\Delta]^L}
m_{j_1^{(\rho)},\cdots,j_L^{(\rho)}} \prod_{a=1}^k
\lambda^{j_a^{(\rho)}}(\bR_+)\prod_{a'=k+1}^L
\lambda^{j_{a'}^{(\rho)}}(C_{a'}),$$ avec $m_{j_1^{(\rho)},
\cdots, j_L^{(\rho)}}  \in \bZ. $ Celui-ci explique l'apparition
du facteur $g(z)=$ $\sum_{(j_1^{(\rho)},\cdots,j_L^{(\rho)} )\in
[1,\Delta]^L} m_{j_1^{(\rho)},\cdots,j_L^{(\rho)}}$ $
\prod_{a=1}^k$ $ e^{2 \pi \sqrt -1 j_a^{(\rho)}{\mathcal
L}_a({\bI, \bz, \zeta})}$ $ \prod_{a'=k+1}^L$ $ e^{2 \pi \sqrt -1
j_{a'}^{(\rho)}{\mathcal L}_{a'}({\bI, \bz, \zeta})} (1- e^{2 \pi
\sqrt -1 {\mathcal L}_{a'}({\bI, \bz, \zeta} ) })$ \`a part de la
fonction $\Gamma (\bullet).$

Dans la suite, on s'occupe de l'analyse des facteurs de la fonction $\Gamma$
chez l'int\'egrale (2.1.5).
Dans ce but, on repr\'esente la matrice ${\sf L}$ (resp.${\sf L}^{-1}$)
comme un ensemble des $L$ colonnes proprement rang\'ees:
$${\sf L} =({\vec v}_1,{\vec v}_2,, \cdots,{\vec v}_L),
{\sf L}^{-1} =({\vec w}_1,{\vec w}_2,, \cdots,{\vec w}_L), {\vec
w}_a =^t({w}_{a,1}, \cdots, {w}_{a,L}). \leqno(2.1.6)$$ Les
vecteurs colonnes de ${\sf L}^{-1} $ sont divis\'ees en 3 groupes:
\par {\bf 1 } les colonnes avec tous les \'el\'ements formellement non-z\'eros.
\par {\bf 2} avec unique \'el\'ement non-z\'ero $(=1)$
qui produit $z_i, 1 \leq i \leq k$.
\par {\bf 3} avec les \'el\'ements non-z\'eros
qui produisent une fonction lin\'eaire en $\zeta +\b1, \bI +\b1$
d'apr\`es (2.1.2).
\par
Dans l'argument plus loin, les premi\`eres deux groupes des colonnes
 nous
importent.

Le produit int\'erieur des vecteurs $(\bI, \bz, \bzeta)$ et ${\vec
w}_a$ d\'efinit la fonction lin\'eaire sous question:
$${\mathcal L}_a({\bI, \bz, \bzeta})= (\bI +\b1, \bz, \bzeta +\b1)
\cdot {\vec w}_a. \leqno(2.1.7)$$ La colonne qui correspond \`a la
variable $log\;s_i$ de $\sf L$ consiste en un seul \'el\'ement
non-z\'ero $(=1)$ \`a la position $\tau_1+ \cdots +\tau_i+i.$
Pourtant la colonne de $\sf L$ qui correspond \`a la variable $log
\;x'_\ell$ aussi consiste d'un \'el\'ement non-z\'ero $(=1)$ qui
n'est pas en position $\tau_1+ \cdots +\tau_i+i, (1 \leq i \leq
k).$
 Notons simplement,
cette correspondance par
$$ {\vec v}_{\rho(i)}= ^t(0,\cdots, 0,
{\rlap{\ ${}^{\sigma(i)\atop{\hbox{${}^{\vee}$}}}$}\cdots} 1, 0,
\cdots,0),$$ qui entra\^ine chez $\sf L^{-1},$
$$ {\vec w}_{\sigma(i)}= ^t(0,\cdots, 0,
{\rlap{\ ${}^{\rho(i)\atop{\hbox{${}^{\vee}$}}}$}\cdots} 1, 0,
\cdots,0).$$ Ici $\rho,\sigma: \{N+1, \cdots , M+k\} \rightarrow
I$ est l'injection qui envoye les num\'eros des colonnes
correspondant aux variables $s,x'$ \`a l'ensemble total des
indices $I$. En suite, on divise les colonnes de $\sf L^{-1}$ en
$k$ groupes $\Lambda_1, \cdots,\Lambda_k \subset I$ dont chacun
repr\'esente $\Lambda_b = \{\tau_1+\cdots +\tau_{b-1}+b ,
\cdots,\tau_1+\cdots + \tau_{b}+b\} \subset I.$ Pour ce groupe, on
peut constater les propri\'et\'es suivantes. $a)$ La colonne
${\vec v}_{M+k+b}= ^t(0,\cdots, 0,0, {\rlap{\ ${}^{\tau_1+\cdots
+\tau_{b-1}+b \atop{\hbox{${}^{\vee}$}}}$}\cdots}
0,1,1,\cdots,\cdots, {\rlap{\ ${}^{\tau_1+\cdots + \tau_{b}+b
\atop{\hbox{${}^{\vee}$}}}$}1}, \cdots,1, 0, \cdots,0),$ avec
$\tau_{b}+1,$ $(1 \leq b \leq k)$ \'el\'ements non-z\'eros $(=1).$
$b)$ Pour les vecteurs $\vec w_a$ du cas $\bf 1$ ci-dessus,
$$\sum_{a \in \Lambda_b}w_{a,j}= 0\;\;{\rm si}\; j \not= M+k+b, 1 \leq b \leq k, \leqno(2.1.8)$$
et puis il existe un autre vecteur du m\^eme groupe $\Lambda_b$
qui satisfait:
$$w_{\sigma(i),j}=\delta_{\rho(i),j},\;\;\leqno(2.1.9)$$
o\`u $\delta_{\cdot, \ast}$ est le symbole delta de Kronecker. Le
vecteur $(2.1.9)$ corresopnd au groupe {\bf 2}.

2) Les  1\`ere $, \cdots,$ ${M+k}-$i\`eme vecteurs rayons de la
matrice $\sf L^{-1}$ sont orthogonaux aux vecteurs $\vec
v_{M+k+1}, \cdots, \vec v_{M+2k} $ ci-dessus. Cela signifie les
relations (2.1.4).

3) L'\'enonc\'e se d\'eduit du point 2).
{\bf C.Q.F.D.}

\par
En vue de la Proposition ~\ref{prop21},  on peut introduire les sous-ensembles
des indices $a \in \{1,2, \cdots, M\}$ comme suit.
\begin{dfn}
Le sous-ensemble $I^+_q \subset \{1,2, \cdots, k\}$ (resp. $I^-_q, I^0_q$)
consiste en indices
$a$ tels que le coefficient $B^a_q$ de
${\mathcal L}_a(\bI, \bz, \bzeta)$ (2.1.2)
soit positif (resp. negatif, z\'ero). D'une fa\c{c}on analogue, on d\'efinit
le sous-ensemble
$J^+_r \subset \{1,2, \cdots, m\}$ (resp. $J^-_r, J^0_r$)
qui consiste en les indices
$a$ tels que le coefficient $C^a_r$  de
${\mathcal L}_a(\bI, \bz, \bzeta )$
soit positif (resp. negatif, z\'ero)
\label{dfn2}
\end{dfn}

{\bf 2.2. Syst\`eme de Horn}
Pour assurer la convergence de la transform\'ee de Mellin inverse de
$M_{{\bI}, \gamma}^\zeta(z)$ de (2.1.1),
on v\'erifie que la transofrm\'ee de Mellin $M_{{\bI},\gamma}^\zeta(z)$
admet l'\'estimation suivante quitte \`a multiplier une fonction $\Delta-$
p\'eriodique
$g(z).$
$$\mid M_{{\bI},\gamma}^\zeta(z)\mid < C_{\bI}
exp(-\epsilon \mid Im\; z \mid )\;\;
{\rm
lorsque}\;  Im\; z \rightarrow \infty,
\mbox{ dans un secteur d'ouverture }\;<2
\pi. $$ pour un certain $\epsilon >0,$

Avant d'\'enoncer un lemme  directement applicable
\`a notre situation, on se rappelle un lemme
\'el\'ementaire pour l'int\'egrale:
$$ \int_{z_0 - i\infty}^{  z_0 + i\infty } s^z g(z)
\prod_{j=1}^{\nu}\frac {\Gamma ( z+ \alpha_j)}{\Gamma ( z+ \rho_j)}
dz . \leqno(2.2.1)$$

\begin{lem}
Si on choit une des fonctions suivantes $g^{+}(z)$ (resp. $g^{-}(z)$) en tant
que $g(z),$
alors l'int\'egrand de (2.2.1) est de d\'ecroissance exponentielle lorsque
$ Im \;z $ tend vers $\infty $ dans le secteur $0 \leq arg \;z < 2\pi,$
(resp. $-\pi \leq arg \;z < \pi .$)

$$ g^{\pm}(z)= 1+ e ^{\pm 2 \pi i \beta_{\nu}}\prod_{j=1}^{\nu}
\frac {sin 2 \pi ( z+ \alpha_j)}{
sin 2\pi ( z+ \rho_j)}, $$
avec $ \beta_{\nu} =-1 + \sum_{j=1}^{\nu}(\rho_j  - \alpha_j)  $
\label{lem221}
\end{lem}

{\bf D\'emonstration}

Il suffit de se rappeler
$$\prod_{j=1}^{\nu}\frac {\Gamma ( x+ iy + \alpha_j)}{\Gamma ( x+iy + \rho_j)}
\rightarrow const. \mid y \mid^{-(\beta_{\nu}+1)}$$
lorsque  $y \rightarrow \pm \infty.$ Ici, on se sert de la formule
de Binet:
$$ log\; \Gamma(z+a) = \log\; \Gamma(z) + a \log\; z - \frac{a- a^2}{2z} + {\cal O}
( \mid z \mid^{-2}   )
$$
si $\mid z \mid >> 1,$ .
Le facteur $\mid s^{-(x+iy)} \mid = r^{-x}e^{\theta y},$
pour $s =  re^{i \theta} $ donne la contribution exponentiellement
d\'ecroissante dans chaque cas.
{\bf C.Q.F.D.}

Nous introduisons la notation
$$  {\mathcal L}_j(z)=
A_{j1}z_1 + A_{j2}z_2 +\cdots + A_{jk}z_k + A_{j0}, \; 1 \leq j \leq p$$
$${\mathcal M}_j (z)=B_{j1}z_1 + B_{j2}z_2 +\cdots + B_{jk}z_k + B_{j0},
\; 1 \leq j \leq r.$$
\begin{lem}
Pour que l'int\'egrale
$$ \int_{\check \Pi} s^{\bz} g(z)
\frac {\prod_{j=1}^{p}\Gamma ({\mathcal L}_j(z))}{\prod_{j=1}^{r}
\Gamma ( {\mathcal M}_j(z))} dz_1\wedge \cdots \wedge dz_k \leqno(2.2.2)$$
avec $g(z)$ une fonction $\Delta-$ p\'eriodique convenable
d\'efinisse une fonction \`a croissance polynomiale
\`a ses lieux singuliers
(y compris $\infty$) il suffit que les conditions suivantes soient safisfait.

i) Pour chaque $i>0$
$$\sum_{j=1}^p A_{j,i}=\sum_{j=1}^r B_{j,i} \leqno{(2.2.3)} $$

ii)
Le nombre
$$\alpha = min_{z \in S^{k-1}}\bigl(\sum_{j=1}^p |{\mathcal L}_j(z)- A_{j0}|-
\sum_{j=1}^r|{\mathcal M}_j(z)- B_{j0}| \bigr)\leqno(2.2.4)$$
soit non-negatif.

\label{lem}
\end{lem}

{\bf D\'emonstration} Pour voir l'existence d'une direction
d'int\'egration lelong laquelle l'int\'egrand de $(2.2.2)$ est de
d\'ecroissance exponentielle, on recourt \`a une astuce de
N\"orlund \cite{Nor} comme au lemme~\ref{lem221}. Si on utilise la
notation comme suit:
$$\beta_i = \sum_{j} A_{j,i}-\sum_{j} B_{j,i}$$
$$\eta = \sum_{j=1}^p A_{j,0}-\sum_{
j=1}^r B_{j,0} -\frac{1}{2}(p-r)$$
l'int\'egrand de $(2.2.2)$ se comporte comme
$ e^{\sum_{j=1}^k(-\frac{1}{2} \alpha\pi |u_j|+ \theta_j u_j)}
\prod_{j=1}^k |v_j|^{\beta_j v_j+\eta}(R_j/\rho_j)^\gamma $
pour certain $\gamma, \rho_1, \cdots, \rho_k \in \bf R$
et $s=(R_1e^{i \theta_1}, \cdots,
R_ke^{i \theta_k})$ lorsque le norme de la variable
$z_j = u_j+iv_j$ tend vers l'infini.
Il est facile de voir cela de la formule de Stirling
(Whittaker-Watson, Chapter XII, Exemple 44).
L'estimation ci-dessus entra\^ine imm\'ediatement le lemme. {\bf C.Q.F.D.}

\begin{remark}
Dans ~\cite{Tsi2} les auteurs ont \'etabli un th\'eor\`eme semblable \`a notre dans le cas o\`u $\alpha$ de (2.2.4) est strictement positif.
\end{remark}

Si on applique ce lemme \`a notre int\'egrale,
on voit qu'il existe un
cycle $\check \Pi$
$$I_{X^{\bI}, \gamma}^{\zeta} (s):
= \int_{\check \Pi} g(z)\frac{ \prod_{a\in I^+_q \cup I^0_q}\Gamma
\bigl({\mathcal L}_a (\bI , \bz , \zeta)\bigr)} {\prod_{\bar{a}\in
I^-_q }\Gamma \bigl (1-{\mathcal L}_{\bar{a}} (\bI , \bz ,
\zeta)\bigr)} s^{-\bz} dz, \leqno(2.2.5)$$ avec une fonction
$g(z)$ rationelle en $e^{\pi i {\mathcal L}_a (\bI , \bz ,
\zeta)}, a \in I$. Ici on se souvient de la relation
$\Gamma(z)\Gamma(1-z)= \frac{\pi} {sin\;\pi z}.$ Une autre
consequence importante de la proposition ci-dessus est le
th\'eor\`eme suivant:

\begin{thm}
L'int\'egrale $I_{X^{\bI}, \gamma}^{\bzeta}(s)$ satisfait
un syst\`eme des
\'equations du type de Horn comme suit:
$$ L_{q ,\bI}(\vartheta_s, s, \zeta)I_{X^{\bI}, \gamma}^{\zeta}(s)=
\Bigl[P_{q, \bI}( \vartheta_s, \zeta) -s_q^\Delta Q_{q, \bI}(
\vartheta_s, \zeta)\Bigr] I_{X^{\bI}, \gamma}^{\zeta}(s)=0,1 \leq
q \leq k \leqno(2.2.6)$$ avec
$$P_{q, \bI}(\vartheta_s, \zeta)=
 \prod_{a\in I^+_q} \prod_{j=0}^{B_q^a-1}
\bigl({\mathcal L}_a(\bI,-\vartheta_{s}, \zeta)+j\bigr)
\leqno(2.2.7)$$
$$Q_{q, \bI}( \vartheta_s, \zeta)
= \prod_{\bar{a}\in I^-_q}\prod_{j=0}^{-B_q^{\bar a}-1}
\bigl(-{\mathcal L}_{\bar{a}}(\bI,-\vartheta_{s}, \zeta)-j \bigr)
,\leqno(2.2.8)$$ o\`u $I^+_q, I^-_q,  1 \leq q \leq k. $ sont les
ensembles des indices d\'efinis dans la D\'efinition ~\ref{dfn2}.
Pourtant les degr\'es des deux op\'eratuers  $P_{q,
\bI}(\vartheta_s, \zeta),$ $Q_{q, \bI}( \vartheta_s, \zeta)$ sont
\'egaux. Notamment ils s'expriment par la relation suivante,
$$\sum_{a \in I^+_q} B_q^a = |\chi (X_q)|= |\chi(\tilde X_q)|, 
\leqno(2.2.9)$$
ici $X_q$ signifie la vari\'et\'e affine de dimension $M-1,$
$$ X_q :=
\{(X,s)|_{s_q=1} \in \bT^{M+k-1}; f_1(X)+s_1=\cdots = f_{q-1}(X)+s_{q-1}=$$
$$
= f_q(X)+1 = f_{q+1}(X)+s_{q+1}= \cdots =f_k(X)+s_k =0\}.$$
Nous notons par $\tilde X_q$ l'IC d\'efinie dans $\bT^{M+1},$
$$ \tilde X_q := \{(x, s_q) \in \bT^{M+1}; f_1(x) = \cdots= f_{q-1}(x)
=f_{q}(x)+s_q = f_{q+1}(x)= \cdots = f_{k}(x)=0 \}.$$

\label{thm23}
\end{thm}
{\bf D\'emonstration} Pour un cycle $\gamma$ qui r\'ealise
l'expression (2.1.1) de la transform\'ee de Mellin, l'\'enonc\'e
est une cons\'equence \'evidente du (2.2.5). Pour l'autre cycle
$\gamma'$, on applique l'argument de \cite{Baty1}, \cite{Baty2},
\cite{Sab1} 3.4.1, qui dit que l'int\'egrand de l'expression
$$\int_{\gamma'}\int_\Gamma
\Bigl[P_{q, \bI}( \vartheta_s, \zeta) -s_q^\Delta Q_{q, \bI}(
\vartheta_s, \zeta)\Bigr] (f_1(X) +s_1)^{-\zeta_1-1} \cdots
(f_k(X) +s_k)^{-\zeta_k-1} X^{\bI} dX ,$$ donne deriv\'ee d'une
forme rationelle sur ${\bf T}^{M} \setminus {X}_1 \cup \cdots \cup
{X}_k $ (notation de la Proposition ~\ref{prop124}).

Nous prouvons ici donc (2.2.9). Soit $E_q$ sous-espace de
$\bR^L_{(z,\zeta, \bI)}$ d\'efini par $z_q=0, \zeta=0.$ Alors on a
l'\'egalit\'e pour $a \in I^+_q$
$$ B_q^a= (M+k-1)! vol_{M+k-1}(\bigl<\vec{v_1}|_{E_q},{\rlap{\ 
${}^{a\atop{\hbox{${}^{\vee}$}}}$}\cdots},
\vec{v_L}|_{E_q}\bigr>),$$ o\`u $\vec{v_i}$ $i-$\`eme vecteur
rayon de la matrice $\sf L.$ En fait
$\bigl<\vec{v_1}|_{E_q},{\rlap{\
${}^{a\atop{\hbox{${}^{\vee}$}}}$}\cdots}  ,
\vec{v_L}|_{E_q}\bigr>$ est formellement un poly\`edre avec $L-1$
vertexes (y compris \{0\}) dans $E_q$ dont le volume est \'egal au
volume d'un simplexe avec $(M+k-1)$ vertexes  dans $E_q$ qui
appartiennent \`a l'un des poly\`edres de Newton
$\Delta(f_1(X)+s_1),
$$\Delta(f_{q-1}(X)+s_{q-1}), $ $\Delta(f_q(X)+1),$
$\Delta(f_{q+1}(X)+s_{q+1}),$ $\cdots,$ $\Delta(f_k(X)+s_k)$.

Dans la matrice $\sf L$ on multiplie les rayons qui correspondent aux termes de
 $f_\ell(X)$ par $\lambda_\ell$ $(\ell \not =q).$ Notons par $\sf L(\lambda) $
la matrice ainsi modi
fi\'ee. Alors $$det ({\sf L(\lambda)})
= \lambda_1^{\tau_1} \cdots \lambda_k^{\tau_k} det ({\sf L}),$$
o\`u $\tau_i =$ nombre des termes dans l'expression $f_{i}(X).$
Evidemment $\tau_1 + \cdots \tau_k+k = L.$
Si on calcule le $(M+q, \ell)-$mineur ( d\'eterminant de la matrice que l'on
obtient de $\sf L$ en enlevant $(M+q)-$i\`eme  colonne et $\ell-$i\`eme rayon
de $\sf L(\lambda)).$
On remarque ici que $M+1, \cdots M+k$ colonnes sont ceux qui correspondent \`a
$z_1,\cdots, z_k.$  Le $(M+q, \ell)-$mineur sous question nous donne
le $(M+q, \ell)-$mineur de $\sf L$ multipli\'e par
$\lambda_1^{\tau_1} \cdots
\lambda_\ell^{\tau_\ell-1}\lambda_k^{\tau_k}$
dont l'homog\'en\'eit\'e en $(\lambda_1, \cdots,\lambda_k)$
est \'egal \`a $L-k-1= M+k-1.$ D'apr\`es la d\'efinition
du volume mixte de Minkowski, le somme de tous
les simplexes possibles engendr\'es par
$M+k-1$ vertexes qui appartiennent \`a l'un des poly\`edres de Newton
  $\Delta(f_1(X)+s_1), $$\Delta(f_{q-1}(X)+s_{q-1}), $
$\Delta(f_q(X)+1),$  $\Delta(f_{q+1}(X)+s_{q+1}),$ $\cdots,$
$\Delta(f_k(X)+s_k)$ doit \^etre le somme des coefficients positifs
de $\lambda_1^{a_1} \cdots \lambda_k^{a_k} ,$ $$a_1 +\cdots + a_k = M+k-1,
a_1 \geq 1, \cdots, a_k \geq 1$$
chez les $(M+q, \ell)-$mineurs avec $ 1 \leq \ell \leq L.$
Ce somme est \'egal \`a $\sum_{a \in I^+_q} B_q^a.$
D'autre part d'apr\`es le th\'eor\`eme de Khovanski-Oka \cite{Kh1},
\cite {Oka}  le caract\'eristique d'Euler de l'intersection compl\`ete $X_q$
s'exprime par des donn\'ees combinatoires,
$$ \frac{1}{(M+k-1)!}|\chi(X_q)|=$$
$$vol_{M+k-1}\bigl(\Delta(f_1(X)+s_1),
\Delta(f_{q-1}(X)+s_{q-1}),
\Delta(f_q(X)+1),\Delta(f_{q+1}(X)+s_{q+1}),\cdots,
\Delta(f_k(X)+s_k)\bigr).$$
L'\'egalit\'e $|\chi(\tilde X_q)|=$ $|\chi(X_q)|$
est une cons\'equence facile d'alg\`ebre lin\'eaire.
{\bf C.Q.F.D.}

On peut formuler une condition de non-r\'esonance pour le
syst\`eme $(2.2.6).$ Notamment, l'ensemble des fonctions
lin\'eaires
$$ \bigcup_{a \in I_q^+}
\{{\cal L}_a(\bI,\bz,\zeta), {\cal L}_a(\bI,\bz,\zeta)+1, \cdots,
{\cal L}_a(\bI,\bz,\zeta)+ B_q^a-1 \},$$
ait aucune intersection avec l'ensemble
$$ \bigcup_{{\bar a} \in I_q^-}\{{\cal L}_{\bar a}(\bI,\bz,\zeta), {\cal
L}_{\bar a}(\bI,\bz,\zeta)+1, \cdots, {\cal L}_{\bar a}(\bI,\bz,\zeta)-
B_q^{\bar a}-1 \}.$$
Dans le cas de non-r\'esonance ci-dessus, le d\'egr\'e
de l'op\'erateur $L_{q, \bI}( \vartheta, \zeta)$ est \'egal
\`a $ |\chi(X_q)|.$

\begin{cor}
La vari\'et\'e caract\'eristique du syst\`eme (2.2.6) est \'egale
\`a la vari\'et\'e
$$\Lambda=\{(s,\xi) \in T^\ast\bC^k_s; \sigma(L_q)(s,\xi)
= \Bigl[P_{q, 0}( s\xi, 0) -s_q^\Delta Q_{q, 0}( s\xi, 0)\Bigr]
=0,1 \leq q \leq k\} $$ o\`u $s\xi=(s_1\xi_1, \cdots, s_k\xi_k).$
Le lieu singulier de (2.2.6) est contenu dans l'ensemble
${\mathcal R} =\{s \in \bC^k_s; R[\sigma(L_1)(s,\xi),
\cdots,\sigma(L_k)(s,\xi)]=0\}$ o\`u $R[\sigma(L_1)(s,\xi),
\cdots,\sigma(L_k)(s,\xi)] \in \bC [s]$ est le r\'esultant des
polyn\^omes $\sigma(L_1)(s,\xi), \cdots,\sigma(L_k)(s,\xi).$
\end{cor}
Malheureusement dans plusieurs cas le r\'esultant
$R[\sigma(L_1)(s,\xi), \cdots,\sigma(L_k)(s,\xi)]$ est
identiquement z\'ero. Pour d\'eduire des informations essentielles
sur le lieu singulier du sys\`eme (2.2.6), on introduit la notion
de ``r\'esultant r\'esiduel''. Soit $L_0(s,\xi)$ le p.g.d.c. des
polyn\^omes $\sigma(L_1)(s,\xi), \cdots,\sigma(L_k)(s,\xi)$ alors
le r\'esultant r\'esiduel est d\'efini comme le r\'esultant des
polyn\^omes $\frac{\sigma(L_1)(s,\xi)}{L_0(s,\xi)},
\cdots,\frac{\sigma(L_k)(s,\xi)} {L_0(s,\xi)}.$ Divers experiments
sur les cas concr\`ets nous permettent de formuler la conjecture
suivante.
\par
{\bf Conjecture}
 Le poly\`edre de Newton du r\'esultant r\'esiduel
$R[\frac{\sigma(L_1)(s,\xi)}{L_0(s,\xi)},
\cdots,\frac{\sigma(L_k)(s,\xi)}{L_0(s,\xi)}]$
d\'epourvu du facteur de forme $s_1^{d_1} \cdots s_k^{d_k}$
coincide avec le somme mixte de Minkowski $\sum_{q=1}^k\sum_{\bar a\in I_q^-}
\Delta(L_{\bar a}) $
pour $ \Delta(L_{\bar a}) :=\{z \in \bR^k_+; \sigma(L_q)(s,\xi)
\leq  \ell_{\bar a}
\}$ pour le choix propre des $\l_{\bar a}.$
\par Ici on d\'efinit le somme mixte de Minkowski pour deux poly\`edres
par $\Delta_1 + \Delta _2 =\{z_1+z_2; z_1 \in \Delta_1, z_2 \in \Delta _2\}.$
On peut trouver des exemples qui supportent cette conjecture
dans \cite{Tsikh} pour des syst\`emes de Horn qui n'ont pas forc\'ement
d'origine g\'eom\'etrique. Il est probable que une \'etude d\'etaill\'ee 
du polyt\^ope secondaire introduit par les auteurs du \cite{GZK}
peut donner une r\'eponse combinatoire \`a cette question.
\begin{cor}
Chaque  solution du syst\`eme (2.2.6) peut \^etre represent\'ee
comme  combinaison lin\'eaire des fonctions:
$$I_{X^{\bI}, \gamma_{\delta}}^{\bzeta} (s):
= \int_{\check \Pi_\delta} \frac{
\prod_{a\in I^+_q \cup I^0_q}\Gamma({\mathcal L}_a (\bI , \bz, \zeta))}
{\prod_{a\in I^-_q }\bigl(\Gamma (1-{\mathcal L}_a (\bI ,
\bz,\zeta)\bigr)} s^{-\bz} dz,$$ associ\'ee \`a un cycle $\check \Pi_\delta
  \subset \{z \in {\bf C}^k; {\mathcal L}_a (\bI , \bz, \zeta) \in
\bZ_{\leq 0}\}.$ \label{cor26}\end{cor} Il s'agit d'une  m\'ethode
g\'en\'eralis\'ee de Frobenius. Voir \cite{Sad1}, Chap. 2, 3.

Ici nous introduisons la notion du support de la transform\'ee de
Mellin $M^\zeta_{{\bI}, \gamma}(z)$ qui contribue essentiellement
\`a l'inversion d\'ecrite au  corollaire ~\ref{cor26}.
\begin{dfn} Pour un indice $q \in [1,k]$ fix\'e
on note par $supp^{(q)}(M^\zeta_{{\bI}, \gamma})(z)$ l'ensemble
des p\^oles d'ordre $\geq k$ de l'expression $$\frac{ \prod_{a\in
I^+_q}\Gamma({\mathcal L}_a (\bI , \bz, \zeta))} {\prod_{a\in
I^-_q }\bigl(\Gamma (1-{\mathcal L}_a (\bI , \bz,\zeta)\bigr)}.$$
\label{dfn26}
\end{dfn}
\begin{remark}
{\em On remarque ici que les r\'eseaux de l'ensemble
$supp^{(q)}(M^\zeta_{{\bI}, \gamma})(z)$ donnent naissance \`a la
bonne $\kappa$-filtration ($\kappa$ =k) introduite par C.Sabbah
\cite{Sab1}, \cite{Sab2}. } \label{remark22}
\end{remark}

\par Nous introduisons une fonction m\'eromorphe qui depend de
$\bz,$ $(\eta^{(1)}, \cdots ,\eta^{(k)} ),$ $(\eta'^{(1)}, \cdots
,\eta'^{(k)})$ $ \in$ $ \bZ^k,$ $(\alpha_1^{(\ell)},$ $ \cdots,
\alpha_k^{(\ell)}),$ $(\beta_1^{(\ell)}, \cdots, \beta_k^{(\ell)})
\in {\bf Q}^k $ et $\Delta$ un entier positif.
$$a(\bz)= \frac{\prod_{\ell=1}^k\Gamma(
\frac{<\alpha^{(\ell)},\bz>}{\Delta}+\eta^{(\ell)})}{
\prod_{\ell=1}^k\Gamma(\frac{ <\beta^{(\ell)},\bz>}{\Delta}+\eta'^{(\ell)})},
$$ o\`u $<\alpha^{(\ell)},\bz>= \alpha_1^{(\ell)}z_1 + \cdots +
\alpha_k^{(\ell)}z_k.$
L'expression
$<\beta^{(\ell)},\bz>$ se d\'efinit d'une fa\c{c}on analogue.
Nous introduisons la notation,
$z+\Delta e_r = (z_1, \cdots, z_{r-1}, z_r+\Delta, z_{r+1}, \cdots, z_k).$

\begin{lem}
Pour la fonction $a(\bz)$ on a la relation suivante:
$$\frac{a(\bz+ \Delta e_q)}{a (\bz)}=
\frac{\prod_{\ell:\alpha_q^{(\ell)} >0 }^k \prod_{j=0}^{\alpha_q^{(\ell)}-1}
\bigl(\frac{<\alpha^{(\ell)}, \bz>}{\Delta} +j+ \eta^{(\ell)}\bigr)
\prod_{m':\beta_q^{(m')} <0 }^k \prod_{j=1}^{-\beta_q^{(m')}}
\bigl(\frac{<\beta^{(m')}, \bz>}{\Delta} - j+ {\eta '}^{(m')} \bigr)}{
\prod_{m:\alpha_q^{(m)} <0 }^k \prod_{j=1}^{-\alpha_q^{(m)}}
\bigl(\frac{<\alpha^{(m)}, \bz>}{\Delta} -j +\eta^{(m)}\bigr)
\prod_{\ell':\beta_q^{(\ell')} >0 }^k \prod_{j=0}^{\beta_q^{(\ell')}-1}
\bigl(\frac{<\beta^{(\ell')}, \bz>}{\Delta} - j + {\eta '}^{(\ell')} \bigr)}.$$

\label{lem23}
\end{lem}
{\bf D\'emonstration}
L'\'enonc\'e s'entra\^ine de la relation recurrente bien connue:
$$\Gamma(\frac{\alpha(n+\Delta)}{\Delta}+\zeta)=
\Gamma(\frac{\alpha n}{\Delta}+\zeta)(\frac{\alpha n}{\Delta}+\zeta)
(\frac{\alpha n}{\Delta}+1+ \zeta)\cdots (\frac{\alpha n}{\Delta}+ \alpha-1+\zeta),
$$ si $\alpha>0.$
$$\Gamma(\frac{\alpha(n+\Delta)}{\Delta}+\zeta)=
\Gamma(\frac{\alpha n}{\Delta}+\zeta)(\frac{\alpha
n}{\Delta}+\zeta-1)^{-1} (\frac{\alpha n}{\Delta}+
\zeta-2)^{-1}\cdots (\frac{\alpha n}{\Delta}+ \zeta+ \alpha)^{-1},
$$ si $\alpha<0.${\bf C.Q.F.D.}
\par
La compatibilit\'e du syst\`eme (2.2.6) au sens d'Ore-Sato
(\cite{Sad1}) se d\'eduit de la proposition ci-dessus.
\begin{prop}
Les coefficients rationels
$$R_q(z)= \frac{P_{q, \bI}(\bz, \zeta)}{Q_{q, \bI}(\bz+\Delta e_q,
\zeta )}, \leqno(2.2.10)$$ d\'efinis pour les op\'erateurs
(2.2.7), (2.2.8) satisfont la condition de compabitilit\'e:
$$R_q(z+\Delta e_r)R_r(z)= R_r(z+\Delta e_q)R_q(z), \;\; q,r =1, \cdots, k.
\leqno(2.2.11) $$
\end{prop}

\vspace{2pc}
{
\center{\section{
Structure de Hodge de l'int\'egrale fibre }}
}
{\bf 3.1.  Structure de Hodge du groupe de cohomologie d'une
hypersurface dans un tore}

Nous nous souvenons des notions fondamentales sur la   structure de
Hodge mixte du groupe de cohomologie d'une hypersurface dans un tore
selon \cite{Baty1}, \cite{DX1}.
\par Soit  $\Pi$ un poly\`edre convexe de dimension $n$
dans ${\bR}^n$ avec tous les vertexes dans $\bZ^n.$ Soit $S_\Pi
\subset$ $\bC[x_1^{\pm}, \cdots, x_n^{\pm}]$ un sous-anneau de
l'anneau des polyn\^omes de Laurent d\'efini comme suivant:
$$S_\Pi:= \bC \oplus \bigoplus_{\frac{\vec \alpha}{k} \in \Pi,
\exists k \geq 1} \bC \cdot x^{\vec \alpha}.\leqno(3.1.1) $$

\par
 Nous notons par  $\Delta(f)$ l'envelope convexe de l'ensemble ${\vec
\alpha} \in supp(f)$ qui s'appelle poly\`edre de Newton d'un
polyn\^ome de Laurent $f(x).$ Nous introduisons l'ideal de Jacobi
comme suit:  $$ J_{f,\Delta (f)}= \bigl<x_1 \frac{\partial
f}{\partial x_1}, \cdots,  x_n \frac{\partial f}{\partial
x_n}\bigr>\cdot S_{\Delta(f)}.  \leqno(3.1.2)$$ Soit $\tau$ une
facette $\ell-$dimensionelle du $\Delta(f)$ et on d\'efinit $$
f^\tau(x)= \sum_{\vec \alpha \in \tau \cap supp(f)}a_{\vec \alpha}
x^{\vec \alpha}, \leqno(3.1.3)$$ o\`u $f(x) =\sum_{\vec \alpha \in
supp(f)}a_{\vec \alpha} x^{\vec \alpha}.$ Le polyn\^ome de Laurent
 $f(x)$ s'appelle $\Pi-$ r\'egulier, si $\Delta(f)=\Pi$ et pour chaque
 facette  $\ell-$dimensionelle $\tau \subset \Delta(f)$ ($\ell >0$)
les \'equations polynomiales:  $$ f^\tau(x)= x_1 \frac{\partial
f^\tau}{\partial x_1}= \cdots=  x_n \frac{\partial f^\tau}{\partial
x_n}=0,$$ ne poss\`edent pas de solutions communes dans $\bT^n =
(\bC^\times)^n.$

\begin{prop}
Soit $f$ un polyn\^ome de Laurent tel que $\Delta(f)=\Pi.$
Alors les conditions suivantes sont \'equivalentes.
\par
(i) Les \'el\'ements $x_1 \frac{\partial f}{\partial x_1}, \cdots,
x_n \frac{\partial f}{\partial x_n}$ donne naissance \`a une suite
r\'eguli\`ere dans $S_{\Delta(f)}.$
\par
(ii)
$$ dim \bigl(\frac{ S_\Delta (f)}{J_{f,\Delta (f)}}\bigr)= n! vol(\Pi).$$
\par
(iii) $f$ est $\Pi-$r\'egulier. \label{prop311}
\end{prop}
Pour un polyn\^ome $f$ $\Pi-$r\'egulier, il est  possible
d'introduire une filtration sur
  $\Pi= S_\Delta (f)$ qui est d\'efinie comme suit. L'inclusion $\vec \alpha \in
S_k $ soit vraie si et seulement si $ \frac{\vec \alpha}{k} \in
\Delta (f).$ Par cons\'equence nous avons une filtration
croissante;
$$\bC \cong \{0\}=S_0 \subset S_1 \subset \cdots \subset S_n \subset
\cdots,$$ dont une filtration d\'ecroissante est d\'eduit sur
$\frac{ S_\Delta (f)}{J_{f,\Delta (f)}}:$ $$ F^n\bigl(\frac{
S_\Delta (f)}{J_{f,\Delta (f)}}\bigr) \subset F^{n-1}\bigl(\frac{
S_\Delta (f)}{J_{f,\Delta (f)}}\bigr) \subset \cdots \subset F^0
\bigl(\frac{ S_\Delta (f)}{J_{f,\Delta (f)}}\bigr).$$ Elle
s'appelle filtration de Hodge de $\frac{ S_\Delta (f)}{J_{f,\Delta
(f)}}.$ Il faut remarquer ici que la filtration de Hodge se
termine par le $n-$i\`eme terme.  \par Nous nous rapelons la
notion du polyn\^ome d'Ehrhart:
\begin{dfn} Soit  $\Pi$ un polyt\^ope convexe de dimension
$n.$  Nous notons la s\'erie de Poincar\'e d'un alg\`ebre gradu\'e
$S_\Pi$ par $$P_\Pi(t)= \sum_{k \geq0} \ell(k\Pi)t^k,$$
$$Q_\Pi(t)= \sum_{k \geq0} \ell^\ast(k\Pi)t^k,$$ o\`u $\ell(k\Pi)$
(resp.$\ell^\ast(k\Pi)$ ) d\'enote le nombre des points entiers
dans $k\Pi.$ (resp. points entiers int\'erieurs dans $k\Pi$ ). Alors
nous appelons
$$\Psi_\Pi(t)= \sum_{k = 0}^n \psi_k(\Pi)t^k= (1-t)^{n+1}P_\Pi(t),$$
$$\Phi_\Pi(t)= \sum_{k = 0}^n \varphi_k(\Pi)t^k =
(1-t)^{n+1}Q_\Pi(t),$$ les polyn\^omes d'Ehrhart. \end{dfn}
Les polyn\^omes d'Erhart satisfont la
relation suivante,
$$t^{n+1}\Psi_\Pi(t^{-1})= \Phi_\Pi(t).$$ Dor\'enavant notre
objet central sera le groupe de cohomologie de l'hypersurface
$Z_f:=\{x \in \bT^n ; f(x)=0\}.$ On a un isomorphisme important
pour la filtration de Hodge filtration du $PH^{n-1}(Z_f).$
\begin{thm}(\cite{Baty1}) Pour la partie primitive $PH^{n-1}(Z_f)$
du $H^{n-1}(Z_f),$ on a isomorphisme suivant; $$
\frac{F^iPH^{n-1}(Z_f)}{F^{i+1}PH^{n-1}(Z_f)} \cong
Gr_F^{n-i}\bigl(\frac{ S_\Delta (f)}{J_{f,\Delta (f)}}\bigr)
=\frac{F^i\bigl(\frac{ S_\Delta (f)}{J_{f,\Delta (f)}}\bigr)}
{F^{i+1}\bigl(\frac{ S_\Delta (f)}{J_{f,\Delta (f)}}\bigr)}.
\leqno(3.1.4)$$ En plus $$ dim \;Gr_F^{n-i} \bigl(\frac{ S_\Delta
(f)}{J_{f,\Delta (f)}}\bigr)= \sum_{q \geq
0}h^{i,q}(PH^{n-1}(Z_f))=\psi_{n-i}(\Delta (f)),$$ pour $ i \leq
n-1.$ \label{thm312}
\end{thm}

Quant'\`a la filtration par le poids, nous avons la caract\'erisation
comme suit. Afin d'introduire la  notion du strate du support de
l'alg\`ebre $\frac{ S_\Delta (f)}{J_{f,\Delta (f)}}$ nous nous
servons de la convention d'identifier un polyn\^ome
$x^{\vec\alpha} \in S_\Delta (f)$ avec $\vec\alpha \in \bZ^n.$ Un
strate de dimension $(n-j)$ du $supp(S_\Delta (f))$ est d\'efini
comme l'ensemble des points $\vec i \in k \Delta (f),$
$k=1,2,\cdots$ tels que $\frac{\vec i}{k}$ se trouve sur la
facette $(n-j)-$dimensionelle du $\Delta (f)$ et au m\^eme temps
il n'appartient \`a aucune facette $(n-j-1)-$dimensionelle  $\Pi'
\subset \Delta (f).$

\begin{thm}
On peut identifier la filtration par le poids sur  $PH^{n-1}(Z_f)$
avec la filtration d\'ecroissante
$$0=W_{n-2}\subset W_{n-1} \subset
\cdots \subset W_{2n-2} = PH^{n-1}(Z_f),$$ telle que
$W_{n+i-1}\cong$$\{$ les points entiers situ\'es  sur le strate du
$supp\bigl(\frac{ S_\Delta (f)}{J_{f,\Delta (f)}}\bigr)$ de
dimension $\geq (n-i)$ mais sur aucun strate de dimension
$(n-i-1)$ $\}$ pour $0 \leq i \leq n-2.$ \label{thm313}\end{thm}
Le th\'eor\`eme ci-dessus est une cons\'equence facile du Theorem
8.2 \cite{Baty1}.  Tout d'abord nous notons que la suite exacte
suivante a lieu, $$ 0 \rightarrow H^n(\bT) \rightarrow H^n(\bT
\setminus Z_f) \stackrel{Res}{\rightarrow} H^{n-1}(Z_f)
\rightarrow 0.$$ L'application de r\'esidu de Poincar\'e $Res$
donne un morphisme de la structure de Hodge mixte du type
$(-1,-1),$ $$ Res(F^j\;H^n(\bT \setminus Z_f)) =
F^{j-1}\;H^{n-1}(Z_f), \;\;Res(W_j\;H^n(\bT \setminus Z_f)) =
W_{j-2}\;H^{n-1}(Z_f).$$ On a donc, $$ 0 \rightarrow
W_{n+i}\;H^n(\bT) \rightarrow W_{n+i}H^n(\bT \setminus Z_f)
\stackrel{Res}{\rightarrow} W_{n+i-2}H^{n-1}(Z_f) \rightarrow 0,$$
pour $i=2, \cdots , n-1$ o\`u $$W_{2n-1}\;H^n(\bT)= \cdots =
W_{n-1}\;H^n(\bT)=0, \leqno(3.1.5)$$ et $dim \;
W_{2n}\;H^n(\bT)=1.$ En tenant compte du $(3.1.5)$ on voit que
l'application de r\'esidu de Poincar\'e $Res$ entra\^ine
l'isomorphisme $$ Res:W_{n+i}H^n(\bT \setminus Z_f)
\stackrel{Res}{\rightarrow} W_{n+i-2}H^{n-1}(Z_f),$$ pour $i=1,
\cdots ,n-1.$ La structure combinatoire $W_{n+i}H^n(\bT \setminus
Z_f),$ $i=1, \cdots ,n-1$ est bien \'etablie par le Theorem 8.2
\cite{Baty1}.

{\bf 3.2. La transition}

Notons par $x'$ les param\`etres suppl\'ementaires qui ont \'et\'e
introduits dans \S 1 pour le cas d\'eficitaire des variables.
 On s'interesse ici \`a la fa\c{c}on de modifier les termes par
les variables suppl\'ementaires
$ x'= (x'_1, \cdots, x'_{m}).$
Nous supposons que le terme $T_{i_j}, 1 \leq j \leq m$
du $(1.3.9)$ contient le facteur $x'_j.$
Nous notons par l'int\'egrale- fibre
$I_{x^{\bI}, \gamma}^{\zeta} (s)$ pour ce choix de modification concret
$$ I_{x^{
\bI}, \gamma}^{\zeta}(s)=\int_\gamma
 (f_1(x,{x}') +s_1)^{-\zeta_1-1}
\cdots (f_k(x,{x}') +s_k)^{-\zeta_k-1} x^{\bI+\b1} \frac{dx}{x^{\b1}}
\frac{d{x}'}{{\tilde x}'^{\b1}} \leqno(3.2.1)$$

\begin{dfn} On consid\`ere la transform\'ee de Mellin (2.1.1)
pour l'int\'egrale (3.2.1) pour $q \in [1,k]$ fix\'e comme
(2.2.5). Nous notons l'envelope convexe de l'ensemble $
supp^{(q)}(M_{{\bI}, \gamma}^{\zeta} (z))$ (d\'efini dans la
D\'efinition ~\ref{dfn26}) qui contient infiniment beaucoup des
points $z$ avec $z_q < 0,$ par $\Sigma_{\bI}^q \subset {\bR}^k.$
Alors nous appelons le bord de cet ensemble $\partial
\Sigma_{\bI}^q$ les spectres de l'int\'egrale-fibre $I_{x^{\bI},
\gamma}^{\zeta} (s).$ \label{dfn321}
\end{dfn}

\begin{prop}
Les spectres de l'int\'egrale-fibre (3.2.1)
$I_{x^{\bI}, \gamma}^{\zeta} (s)$
sont donn\'es par
des hyperplanes affines:
$$ \partial \Sigma_{\bI}^q = \partial \Bigr\{\bigcup_{a \in I} \{z ;
{\mathcal L}_a({\bI, \bz, \zeta}) \leq 0,\; \rm {pour}\; \it{a \in
I^+_q, z_i \geq 0, i\not = q \;} \} \Bigl\}.$$ \label{prop321}
\end{prop}

{\bf D\'emonstration} Dans le domaine d\'efini comme  $\{z \in
\bR^M ; {\mathcal L}_a({\bI, \bz, \zeta})< 0,$ $ p^{(a)}_q >0,$
$z_q <0$ et $z_i >0, (i \not = q ) \}$ pour $a \in [1,L]=I $
fix\'e, on trouve infiniment beaucoup de points-p\^oles chez une
repr\'esentation de la transform\'ee de Mellin index\'ee par $q
\in [1,k],$
$$ M_{{\bI}, \gamma}^{\zeta} (z)=\frac{
\prod_{a\in I^+_q \cup I^0_q}\Gamma({\mathcal L}_a (\bI , \bz,
\zeta))} {\prod_{a\in I^-_q}
 \Gamma \bigl( 1-{\mathcal L}_a (\bI , \bz,
\zeta\bigr)}. $$ Pourtant ${\mathcal L}_a({\bI, \bz, \zeta})=0$
est une \'equation qui prend part \`a la d\'efinition du $\partial
\Sigma_{\bI}^q$ comme un ensemble semi-alg\'ebrique. {\bf
C.Q.F.D.}

Bien entendu les spectres $\partial \Sigma_{\bI}^q$ se diff\`ere
de ceux de $$\tilde I_{x^{ \bI}, \gamma}^{\zeta}(s)= \int_{\tilde
\gamma}
 (f_1(x,{\tilde x}') +s_1)^{-\zeta_1-1}
\cdots (f_k(x,{\tilde x}') +s_k)^{-\zeta_k-1} x^{\bI +\b1} \frac{dx}{x^{\b1}}
\frac{d{\tilde x}'}{{\tilde x}'^{\b1}}
$$ que nous  notons $\partial \tilde \Sigma_{\bI}^q$
qui est d\'efinie pour une autre fa\c{c}on d'ajouter des
param\`etres supplementaires
${\tilde x}'= ({\tilde x}'_1, \cdots,{\tilde x}'_{m})$
dans laquelle le terme $T_{{\tilde i}_j}, 1 \leq j \leq m$
du $(1.3.9)$ contient le facteur ${\tilde x}'_j.$

On a d\'ej\`a remarqu\'e telle sorte d'ambiguit\'e de choix des
param\`etres suppl\'ementaires dans l'exemple {\bf 1.6}, $(1.6.2)$
et $(1.6.7).$ Heureusement cette diff\'erence est controlable par
une proc\'edure assez simple.

\begin{prop}
Soient
$$ Log\; \Xi
= ^t(log\; x_1, \cdots, log\; x_N, log\; x_1', \cdots, log\;
x_{m}' , log\; s_1, \cdots, log\; s_k, log\; y_1, \cdots, log\;
y_k), $$ et
$$ Log\; \tilde \Xi = ^t(log\; x_1, \cdots, log\; x_N, log\;
{\tilde x_1}', \cdots, log\; {\tilde x_{m}'},  log\; s_1, \cdots,
log\; s_k,   log\; y_1, \cdots, log\; y_k), $$ deux syst\`emes des
variables-param\`etres diff\'erents. Soient les param\'etrisations
de Cayley correspondant:
$$ Log\; T = {\sf L}\cdot Log\;\; \Xi
=\tilde {\sf L}\cdot Log\;\;\tilde \Xi.$$
Alors les spectres de l'int\'egrale-fibre $I_{x^{\bI}, \gamma}^{\zeta} (s)$
sont transform\'ees dans ceux de $\tilde I_{x^{\bI}, \gamma}^{\zeta} (s)$
par la r\`egle suivante:
$$ \tilde {\sf L}\cdot {\sf L}^{-1}(\partial \Sigma_{\bI}^q )=
\partial \tilde \Sigma_{\bI}^q.$$
\label{prop33}
 \end{prop}

En fait cette transition peut \^etre mieux comprise dans le cadre
de la th\'eorie de Terasoma sur la compactification de l'espace
produit des fibres et celui   des param\`etres (voir Remarque
~\ref{remark41} ci-dessous).

{\bf 3.3   }
Il existe un syst\`eme de poids pour
les variables $(x,x'',s, y)$ pour que la relation suivante associ\'ee au
mon\^ome $T_{\nu}(y,X), 1 \leq \nu \leq L$ de $(1.3.9)$ soit vraie.
Il existe le champ d'Euler
$E^{(\nu)},$
$$ E^{(\nu)}(x,x'',1,y,s)= \leqno(3.3.1)$$
$$\sum_{i=1}^Nw_i^{(\nu)}x_i \frac{\partial}{\partial x_i}
+ \sum_{j=1}^{m}{w'}_j^{(\nu)} x'_j \frac{\partial}{\partial x'_j}-
x'_{m+1} \frac{\partial}{\partial x'_{m+1}}
+\sum_{\ell=1}^k (q_{\ell}^{(\nu)}
y_\ell \frac{\partial}{\partial y_\ell}+
{p}_{\ell}^{(\nu)} s_\ell \frac{\partial}{\partial s_\ell}),$$
qui agit sur le polyn\^ome
$F(x,x",1,y,s) +(x'_{m+1}-1) T_{\nu}(x,s,y)$
en sorte que
$$E^{(\nu)}(x,x",1,y,s)
(F(x,x",1,y,s) - T_{\nu}(x,s,y))=
0,$$

$$
E^{(\nu)}(x,x",1,y,s) T_{\nu}(x,s,y)= T_{\nu}(x,s,y).
\leqno(3.3.2)$$ Autement dit, le polyn\^ome $F(x,x",1,y,s)
+(x'_{m+1}-1) T_{\nu}(y,X)$ est quasihomog\`ene par rapport au
poids $(3.3.1)$ qui est d\'etermin\'e de telle fa\c{c}on que
${w'}_{m+1}^{(\nu)} =w^{(\nu)}(x'_{m+1})=-1.$
 Le remplacement de
$T_{\nu}(x,s,y)$ par $x'_{m+1} T_{\nu}(x,s,y)$ dans $F(x,x",1,y,s)$
entra\^ine la quasihomog\'en\'eit\'e.
Nous introduisons les vecteurs-poids pour chaque $\nu$:
$$ \vec{w}^{(\nu)}=({w}_1^{(\nu)}, \cdots , w_N^{(\nu)}),$$
$$\vec{w'}^{(\nu)}=({w'}_1^{(\nu)}, \cdots , {w'}_m^{(\nu)}),$$
$$ \vec{p}^{(\nu)}=({p}_1^{(\nu)}, \cdots , p_k^{(\nu)}). \leqno(3.3.3)$$
Nous appelons une quasihomog\'en\'eit\'e $w^{(\nu)}(\bullet)$
triviale s'il existe $i_\nu$ tel que   ${p}_i^{(\nu)}=0$
pour tous les indices
$i \not = i_\nu$ et ${p}_{i_\nu}^{(\nu)}>0$.
Nous d\'esignons par $I^T \subset I $ (resp. $I^{NT}$ )
l'ensemble des poids triviaux (resp. non-triviaux). Il est 
facile de voir que $|I^T| =k. $

Les autres quasihomog\'en\'eit\'es seront appel\'ees non-triviales.
\begin{lem}
1). La quasihomog\'en\'eit\'e d\'efinie par $(3.3.1)$ correspond au
vecteur-colonne
$\vec{v}_\nu$ de la matrice ${\sf L}^{-1},$ $(2.1.6).$Et elle donne
naissance \`a la fonction lin\'eaire ${\mathcal L}_\nu({\bI, \bz,
\bzeta}),$ $1 \leq \nu \leq L.$ \par 2) Il existe $M+1$ poids
non-triviaux.  C'est \`a dire  $|I^{NT}| =M+1. $ \label{lem31}
\end{lem}
{\bf D\'emonstration}
1). Pour simplifier, supposons que l'expression $y_1f_1(X)$
de (1.1.4) ne contient pas des variables $x'$. Alors le vecteur ${\vec v}_1$
de $(2.1.6)$ correspond \`a la quasihomog\'en\'eit\'e non-trivaile
qui provient de la multiplication $x_1^{\alpha_{1,1}} \rightarrow x'_{m+1}
x_1^{\alpha_{1,1}}$ car ${\vec v}_1$ est vertical \`a tout les
vecteurs horizontaux de $\sf L$ sauf le premier,
$$
\left [\begin {array}{cc} {\sf L}&1\\
      &0\\
      &\vdots\\
      &0
\end {array}\right ]
\left(\begin {array}{c} \vec{v}_1\\
           -1
\end {array}\right)=0. \leqno(3.3.4)$$
Ainsi en g\'en\'eral le vecteur ${\vec v}_i$
nous donne la quasihomog\'en\'eit\'e (soit triviale soit non-triviale)
qui provient
de la multiplication $T_i(x,s,y) \rightarrow x'_{m+1}T_i(x,s,y).$
2).
Pour chaque polyt\^ope simplicial, il existe $M+1$ fa\c{c}ons d'ajouter
un vertexe nouveau tel qu'il produise une quasihomog\'en\'eit\'e non-triviale.
En fait, les cas triviaux correspondent soit \`a la multiplication des types
$s_\ell \rightarrow x'_{m+1}s_\ell,$ ($k$ cas) $x_j'T_i(x,x'',s,y) \rightarrow
x'_{m+1}x_j'T_i(x,x'',s,y)$ ($m$ cas),
soit \`a la multiplication telle que les r\`egles
$\bf a$ et $\bf b$ du lemme 1.5 soient viol\'ees ($k-1$ cas).
Nous avons $L-(m+2k-1)= M+1.$
{\bf C.Q.F.D.}

Le syst\`eme de poids  $\vec{w}^{(\nu)}(\bullet),
\vec{w'}^{(\nu)}(\bullet),$ $\nu \in I$
est donn\'ee par les vecteurs-poids (3.3.3) ci-dessous.
\begin{thm}
1).
Pour un mon\^ome $ X^{\bI+\b1}y^{\zeta+\b1} \frac{ dX \wedge dy}{dF}
\in Gr_F^rPH^{M+k-1}
(Z_{F(X,0,y)}),$
l'inclusion suivante a lieu:
$$ \partial \Sigma^{\ell}_{\bI}
\supset \partial \Bigl\{ \bigcup_{\nu \in I^+_{\ell}} \bigl\{
\{ z;
{\mathcal L}_{\nu}({\bI, \bz,
 \bzeta})\geq 0, p^{(\nu)}_\ell >0 \} \cap_{q \not = \ell }  \{z_q \geq 0
\} \bigr\} \Bigr\}\leqno(3.3.5)$$
o\`u l'in\'egalit\'e comme suit est satisfaite,
$$ M+k-r < {\mathcal L}_{\nu}({\bI, 0,
 \bzeta})
\leq M+k-r+1. \leqno(3.3.6)$$
\par
2). Quant \`a la
filtration par le poids on a l'\'estimation suivante.
Pour l'\'el\'ement  $X^{\bI+\b1}y^{\zeta+\b1} \frac{dX \wedge dy}{dF} $
$\in$ $ Gr_F^r Gr^W_{M+k-1+q} PH^{M+k-1} (Z_{F(X,0,y)}),$
$1 \leq q \leq M+k-2,$
il existe $q-$ paire d'indices
$\{{\nu_1}, \cdots ,{\nu_q}\} \subset I^+_{\ell}  $
tel que  l'intersection des p\^oles de $M_\bI^\zeta(z)$
$$\partial \Sigma^{\ell}_{\bI}
\supset \{\bz \in \bC^k;
 {\mathcal L}_{{\nu_1}}({\bI, \bz,
 \bzeta})= \cdots
={\mathcal L}_{{\nu_q}}({\bI, \bz, \bzeta})= 0\} \not= \emptyset.$$
Par contre, si on prend $I^+_{\ell} \ni \nu_{q+1} \not \in $
$\{{\nu_1}, \cdots ,{\nu_q}\}$ quelconque, on a
$$  \bigcap_{j=1}^{q+1}
\{ {\mathcal L}_{{\nu_j}}({\bI, \bz, \bzeta})=0\} = \emptyset.$$
\label{thm31}
\end{thm}
{\bf D\'emonstration} 1) On applique le Th\'eor\`eme ~\ref{thm312}
au simplexe $\Pi $ d\'efini comme $\Pi = \Delta(F(X,0,y)+1).$ On
remarque la relation
$$
\Pi = \{\{(\bI+\b1, \zeta+\b1) \in \bR^{M+k};0 \leq
{\mathcal L}_{\nu}({\bI, 0, \bzeta}) \leq 1, \nu \in I^T\cup I^{NT}\}
$$ qui
d\'ecoule du fait que le vecteur $(\vec{w}^{(\nu)}, \vec{w'}^{(\nu)},
\vec{p}^{(\nu)},-\vec{p}^{(\nu)})$ non-trivial est vertical \`a
l'hyperplan engendr\'e par les vertexes du simplexe $\Pi$
except\'es le vertexe $\vec \ell_\nu$ qui correspond au $\nu-$i\`eme
rayon de la matrice $\sf L.$ 
On se souvient ici que $|I^T\cup I^{NT}|= M+k+1 $
d'apr\`es lemme 3.6, 2). Puisque pour tel $\nu,$ le mon\^ome 
$T_\nu(X,y)$ ne
d\'epend pas des variables $(s_1,\cdots,s_k),$
on peut consid\'erer la projection $pr(\vec \ell_\nu) \in \bR^{M+k}.$
On voit que
$$ \langle pr(\vec \ell_{\nu'}), (\vec{w}^{(\nu)}, \vec{w'}^{(\nu)},
-\vec{p}^{(\nu)})\rangle =0, \;\; si \;\;\nu' \not = \nu,
$$
$$ \langle pr(\vec \ell_{\nu}), (\vec{w}^{(\nu)}, \vec{w'}^{(\nu)},
-\vec{p}^{(\nu)})\rangle =1, \;\;
$$
 Ces relations nous donnent la d\'escription d\'esir\'ee du simplexe
$\Pi.$ On remarque ici,
$$ \Delta(F(X,0,y)) = \Pi \cap \{(\bI,\zeta) \in \bR^{M+k}; <\zeta, \b1>=1\},$$
et $ \Pi \cap \bZ^{M+k} = \{0\}\cup (\Delta(F(X,0,y))\cap \bZ^{M+k}).$
Il est donc possible d'appliquer les argments ci-dessus du 
$PH^{M+k-1} (Z_{F(X,0,y)+1})$
au groupe $PH^{M+k-1} (Z_{F(X,0,y)}).$
D'apr\`es la proposition ~\ref{prop321}, le bord
$\partial \Sigma_\bI^\ell$ est d\'efini par
${\mathcal L}_{\nu}({\bI, \bz, \bzeta})=0,$ $\nu \in
I_\ell^+.$
En revanche les hyperplans
$$\{(\bI +\b1, \zeta+\b1) \in \bR^{M+k};
{\mathcal L}_{\nu}({\bI, 0, \bzeta})= 0, \nu \in I_\ell^+ \cup
I_\ell^-\}$$ contiennent le bord du simplexe $(M+k-r)\cdot \Pi$ si
$X^{\bI+\b1}y^{\zeta +\b1} \frac{ dX \wedge dy}{dF} \in
Gr_F^rPH^{M+k-1}(Z_{F(X,0,y)}).$ On en d\'eduit que l'ensemble  
$\bigl\{ \{
z; {\mathcal L}_{\nu}({\bI, \bz, \bzeta}) \geq 0, p^{(\nu)}_\ell
>0 \} \cap_{q \not = \ell } \{z_q \geq 0 \} \bigr\}$ est contenu
dans l'ensemble compl\'ementaire de $\Sigma^\ell_\bI$ dans
$\bR^M.$

2)
Si le point $(i^0 +\b1, i'^0+\b1, \zeta^0+\b1)=(\bI^0+\b1, \bz^0+\b1) 
\in \bR^{M+k}$
satisfait un syst\`eme d'\'equations
$$
{\mathcal L}_{\nu_j}({\bI^0, 0, \bzeta^0})=
\langle \vec{p}^{(\nu_j)},-(\zeta^0+\b1)\rangle
+\langle \vec{w}^{(\nu_j)}, i^0+\b1 \rangle
+\langle \vec{w'}^{(\nu_j)}, i'^0+\b1\rangle=0,
\leqno(3.3.7)$$
pour $\nu_j,$ $j=1,2,\cdots,q \leq M+k-2,$ alors il existe une solution
$\bz^0$ de multiplicit\'e $q$ telle que
$$
{\mathcal L}_{\nu}({\bI^0, \bz^0, \bzeta^0})=
\langle \vec{p}^{(\nu_j)},\bz^0-(\zeta^0+\b1)\rangle
+\langle \vec{w}^{(\nu_j)}, i^0+\b1 \rangle
+\langle \vec{w'}^{(\nu_j)}, i'^0+\b1 \rangle=0. \leqno(3.3.7)'$$
Car ce syst\`eme l\`a est \'equivalent au $(3.3.7)$
apr\`es remplacement
de la variable $\zeta^0 +\b1$ par $\bz^0-(\zeta^0+\b1).$ \'Evidemment
tel $\bz^0$ s'associe \`a  un hyperplan de codimension $q$ dans
$\bR^{M+k}_{\bI,\zeta}\times \bR^{k}_{\bz}$ qui passe par
$(\bI^0 +\b1, \bz^0, \zeta^0+\b1).$
D'apr\`es le th\'eor\`eme ~\ref{thm313}, le fait que le point
$(i^0+\b1, i'^0+\b1, \zeta^0+\b1)$ satisfait le syst\`eme
$(3.3.7)$ est \'equivalent
 \`a l'appartenance du mon\^ome $x^{i^0+\b1}x'^{i'^0+\b1}
 y^{\zeta^0+\b1}\frac{dX\wedge dy}{dF} \in W_{M+k-1+q}$
du groupe $PH^{M+k-1}(Z_{F(X,0,y)}).$
Pour consid\'erer la graduation $Gr^W_{M+k-1+q},$ il suffit
d'exclure les strates qui correspondent \`a $W_{M+k-1+q'},$ $ 0 \leq q'
\leq q.$
{\bf C.Q.F.D.}
\begin{remark}
En se servant de l'isomorphisme $(1.3.6),$ on peut formuler
le th\'eor\`eme pour
$$ \frac {X^{\bI} dX}
{f_1(X)^{\zeta_1+1}
\cdots f_k(X)^{\zeta_k+1}}\in PH^M(\bT^M \setminus X_1
\cup \cdots \cup X_k),$$
par rapport \`a la structure de Hodge
induite sur $PH^M(\bT^M \setminus X_1
\cup \cdots \cup X_k)$.
\label{remark331}
\end{remark}

\begin {lem}
Supposons que
un syst\`eme des \'equations lin\'eaires
$$ {\mathcal L}_a({\bI, \bz, \zeta})=-n_a, a \in J_q \subset I, n_a \in
{\bf Z}_{\geq 0}, \leqno(3.3.8)$$
avec $q=p+k$
poss\`ede une solution $z= \dot z$ telle que $\dot z_i
\in \bZ_\leq 0$ pour au moins une indice $i \in \{1, \cdots, k\}.$
En plus pour $q=p+k+1$ le syst\`eme (3.3.8)
des \'equations n'a pas
des solutions comme ci-dessus. Alors la monodromie locale
de l'int\'egrale
$I_{x^{\bI}, \gamma}^{\zeta}(s)$
autour de l'origine $s=0$ a cellule de Jordan de taille $p+1.$
\label{lem321}
\end{lem}
En vue du lemme ci-dessus,
la monodromie locale
de l'int\'egrale
$I_{X^{\bI}, \gamma}^{\zeta}(s)$
autour de l'origine $s=0$ peut avoir une cellule de Jordan de taille
$q+1$
pour certains $\gamma \in H_{M-k}(X_s).$
Ce fait entra\^ine l'\'enonc\'e sur la filtration par le poids comme au
\cite{Var}.
\par
Nous formulons une interpr\'etation de notre 
Th\'eor\`eme ~\ref{thm31} en se servant de 
la formulation classique de P.Deligne \cite{Del}. 
Soient $X^{\bI^0+\b1}y^{\zeta^0+\b1}\frac{dX\wedge dy}{dF} \in Gr^W_{M+k-1+q}
PH^{M+k-1}(Z_{F(X,0,y)})$ et $(\bI^0+\b1, \bz^0, \zeta^0+\b1)$
un point qui satisfait l'\'equation $(3.3.7)'$ tel que les conditions du 
Th\'eor\`eme  3.7, 2) soient vraies. On consid\`ere 
 ${\mathcal L}_a({\bI, \bz, \zeta})$ comme fonction 
lin\'eaire d\'efinie dans l'espace affine $(\bI,\bz) \in \bC^{M+k}.$
Nous d\'esignons par $\tilde X$ un voisinage ouvert suffisamment petit
du point
$(\bI^0+\b1, \bz^0) \in \bC^{M+k}$ 
(pour chaque $\bz$ d'un point $(\bI,\bz) \in \tilde X$ les conditions 
du Th\'eor\`eme  3.7, 2) soient vraies)
et par $Y_{\nu_j}$ l'hyperplan,
$$Y_{\nu_j} := \{(\bI^0+\b1, \bz^0) \in \tilde X; 
{\mathcal L}_{\nu_j}({\bI, \bz, \zeta^0})=0    \}, \;\; 1 \leq j \leq q. $$
D\`es que $1 \leq q \leq M+k-2,$
la reunion des hyperplans,
$$ Y:= Y_{\nu_1} \cup Y_{\nu_2}\cup \cdots \cup Y_{\nu_q},$$
est un diviseur \`a croisements normaux. D'apr\`es 3.1, \cite{Del}, on peut 
d\'efinir $\Omega^1_{\tilde X} <\!Y\!> $ le sous ${\mathcal O}_{\tilde X}-$ 
module localement libre 
engendr\'e par  $\Omega^1_{\tilde X}$ et $\frac{d_{\bz,\bI} 
{\mathcal L}_{\nu_j}({\bI, \bz, \zeta^0})}
{{\mathcal L}_{\nu_j}({\bI, \bz, \zeta^0})}$
en telle sorte que 
 $\Omega^p_{\tilde X}<\!Y\!>=$  $\bigwedge^p 
\Omega^1_{\tilde X}<\!Y\!>$ soit un 
${\mathcal O}_{\tilde X}-$module des 
$p-$formes diff\'erentielles sur $\tilde X$
\`a p\^ole logarithmique le long de $Y.$
Nous consid\'erons une application $ M\!-\!Res$ (Mellin- R\'esidu),
$$ M\!-\!Res:\prod_{a \in I} {\mathcal L}_{a}({\bI, \bz, \zeta^0}) d\bz |\rightarrow \sum 
\bigwedge_{a \in \tilde I; |\tilde I|=M+k, \tilde I \subset I^T \cup I^{NT}} 
\frac{d_{\bz,\bI} {\mathcal L}_{a}({\bI, \bz, \zeta^0})}
{{\mathcal L}_{a}({\bI, \bz, \zeta^0})} \in  \Omega^{M+k}_{\tilde X}<\!Y\!>.$$
Si on se souvient de la filtration par le poids
$ W_q(\Omega^p_{\tilde X}<\!Y\!>)$ du \cite{Del}, 
un sous module de $\Omega^p_{\tilde X}<\!Y\!>$
form\'e des combinaisons lin\'eaires de produits,
$$ \alpha \wedge \frac{d_{\bz,\bI} {\mathcal L}_{\nu_1}({\bI, \bz, \zeta^0})}
{{\mathcal L}_{\nu_1}({\bI, \bz, \zeta^0})} \wedge \cdots 
\wedge \frac{d_{\bz,\bI} {\mathcal L}_{\nu_r}({\bI, \bz, \zeta^0})}
{{\mathcal L}_{\nu_r}({\bI, \bz, \zeta^0})},$$
avec $\alpha \in \Omega^{M+k-r}_{\tilde X},$ 
$0 \leq r \leq q,$ on a l'inclusion suivante,
$$ M\!-\!Res ({M}^{\zeta^0}_\bI(\bz) d\bz) \subset 
W_q(\Omega^{M+k}_{\tilde X}<\!Y\!>).
\leqno(3.3.9)$$
\vspace{2pc}
{
\center{\section{ Quelques applications
.}}
}

{\bf 4.1. $A-$Fonction Hyperg\'eom\'etrique de Gel'fand-Kapranov-Zelevinski}

Si nous rempla\c{c}ons les variables $x'$ par $s'=(s_1', \cdots,
s_{m}')$ consid\'er\'ees comme les variables de d\'eformation des
polyn\^omes $f_1(x), \cdots, f_k(x),$ nous pouvons regarder
l'int\'egrale $$  I_{x^{\bi}, \gamma}^{\zeta}(s,s')=
\int_{\gamma_{s,s'}}
 (f_1(x,{  s}') +s_1)^{-\zeta_1-1}
\cdots (f_k(x,{  s}') +s_k)^{-\zeta_k-1} x^{\bi+\b1} \frac{dx}{x^{\b1}},
\leqno(4.1.1)$$    et sa transform\'ee de Mellin
$$ M_{{\bi},\gamma}^\zeta ({\bz},{\bz'} ):=\int_\Pi s^{\bz }s'^{\bz'+1}
\tilde I_{x^{\bi}, \gamma}^{\zeta} (s,s')
\frac{ds}{s^{\b1}}\wedge \frac{ds'}{s'^{\b1}},
\leqno(4.1.2)$$ au lieu des $(1.3.5)$ et $(1.3.7).$
En fait, le calcul totalement parall\`ele \`a celui du \S 2 nous donne,
$$ M_{{\bi},\gamma}^\zeta ({\bz},{\bz'} )= g(\bz,\bz') \prod_{a \in I}
\Gamma({\cal L}_a (\bi,\bz,\bz',\zeta)),\leqno(4.1.3)$$ o\`u la
fonction lin\'eaire ${\mathcal L}_a(\bi,\bz, {\bz}', \bzeta)$ est
d\'efinie d 'une fa\c{c}on analogue \`a (2.1.2). Il s'en d\'eduit
que l'int\'egrale ${\tilde I}_{x^{\bi}, \gamma}^{\zeta}(s,s')$
satisfait le syst\`eme ci-dessous.
$${\mathcal L}_a({\bi,\bz, \bz', \bzeta} ):=
\frac{\sum_{j=1}^N A_j^a (i_j+1)
+\sum_{\ell=1}^k \left( B_\ell^a z_\ell + D_\ell^a(\zeta_\ell+1)\right)
+\sum_{j=1}^{m} C_j^a(z'_j+1)
}{\Delta}, a \in I \leqno(4.1.4)$$
$$ \tilde L_{q ,\bi}(\vartheta_s,\vartheta_{s'} s, s', \zeta)
{\tilde I}_{X^{\bi}, \gamma}^{\zeta}(s,s'):=
\Bigl[\tilde P_{q, \bi}( \vartheta_s, \vartheta_{s'}, \zeta)
-s_q^\Delta Q_{q, \bi}( \vartheta_s,\vartheta_{s'}, \zeta)\Bigr]
{\tilde I}_{x^{\bi}, \gamma}^{\zeta}(s,s')=0,1 \leq q \leq k \leqno(4.1.4)_1$$
avec
$$\tilde P_{q, \bi}(\vartheta_s, \vartheta_{s'}, \zeta)=
 \prod_{a\in I^+_q} \prod_{j=0}^{B_q^a-1}
\bigl({\mathcal L}_a(\bi,-\vartheta_{s}, -\vartheta_{s'},
\zeta)+j\bigr), \leqno(4.1.4)_2$$ $$\tilde Q_{q, \bi}( \vartheta_s,
\vartheta_{s'}, \zeta) = \prod_{\bar{a}\in
I^-_q}\prod_{j=0}^{-B_q^{\bar a}-1} \bigl(-{\mathcal
L}_{\bar{a}}(\bi,-\vartheta_{s},-\vartheta_{s'}, \zeta)-j \bigr)
,\leqno(4.1.4)_3$$
o\`u $I^+_q, I^-_q,  1 \leq q \leq k $
sont les ensembles des indices
d\'efinis dans la D\'efinition ~\ref{dfn2}.
$$ \tilde L'_{r ,\bi}(\vartheta_s,\vartheta_{s'} s, s', \zeta)
{\tilde I}_{X^{\bi}, \gamma}^{\zeta}(s,s'):=
\Bigl[\tilde P'_{r, \bi}( \vartheta_s, \vartheta_{s'}, \zeta)
-s_r'^\Delta Q'_{r, \bi}( \vartheta_s,\vartheta_{s'}, \zeta)\Bigr]
{\tilde I}_{x^{\bi}, \gamma}^{\zeta}(s)=0,1 \leq q \leq k \leqno(4.1.4)_4$$
$$\tilde P'_{q, \bi}(\vartheta_s, \vartheta_{s'}, \zeta)=
 \prod_{a\in J^+_r} \prod_{j=0}^{C_r^a-1}
\bigl({\mathcal L}_a(\bi,-\vartheta_{s}, -\vartheta_{s'}, \zeta)+j\bigr)
\leqno(4.1.4)_5$$
$$\tilde Q_{q, \bi}( \vartheta_s, \vartheta_{s'}, \zeta)
= \prod_{\bar{a}\in J^-_r}\prod_{j=0}^{-C_r^{\bar a}-1}
\bigl(-{\mathcal L}_{\bar{a}}(\bi,-\vartheta_{s},-\vartheta_{s'},
\zeta)-j \bigr). \leqno(4.1.4)_6$$

\begin{remark}

Comme $m = dim \;{\bf D}(\Sigma) -1$
 $= dim \;{\bf D}(\tilde \Sigma)$ (voir (1.4.5)), on peut
consid\'erer que le syst\`eme $(4.1.4)_{\ast}$ est d\'efini sur
$(\bC^\times)^k \times ({\bf D}(\Sigma)/\bC^\times) $
(resp. $(\bC^\times)^k \times {\bf D}(\tilde \Sigma)$)  
pour ${\bf D}(\Sigma)$ (resp. ${\bf D}(\tilde \Sigma)$) le 
groupe de N\'eron-Severi. \label{remark411}
\end{remark}
D'apr\`es le point 2.4 de \cite{Kap1} les expressions
$(4.1.3), (4.1.4)_{\ast}$
 donnent l'uniformisation de Horn du discriminant des
polyn\^omes $f_1(x,{s}') +s_1,\cdots, f_k(x,{s}') +s_k.$

$$ D_{s,s'} = \{(s,s') \in \bC^{k+m}; f_1(x,{s}') +s_1=\cdots=
f_k(x,{s}') +s_k=0,
rank
\left (
\begin {array}{c} grad_x \; f_1(x,{s}')\\
\cdots \\ grad_x \; f_k(x,{s}')\\
\end{array}
\right) < N, $$
\begin{flushleft}
$\rm{pour \; certain}\;$ $x\in
\bT^N\}.$
\end{flushleft}
Soit $\Delta_f(s,s')$ le polyn\^ome qui d\'efinit $D_{s,s'}$
sans multiplicit\'e.
Pour les fonctions suivantes qui sont $\Delta-$i\`eme racines des fonctions
rationelles,
$$ \psi_q(z,z') = \left(\frac{\prod _{a \in I_q^+}(\sum_{\ell=1}^k
B_\ell^a z_\ell
+\sum_{j=1}^{m} C_j^a z'_j)^{B_q^a}} {\prod _{\bar a \in I_q^-}
(\sum_{\ell=1}^k - B_\ell^{\bar a} z_\ell- \sum_{j=1}^{m} C_j^{\bar
a} z'_j)^{-B_q^{\bar a}}}\right)^{\frac{1}{\Delta}}, 1 \leq q \leq
k,$$
$$ \phi_r(z,z') = \left(\frac
{\prod _{a \in J_r^+}(
\sum_{\ell=1}^k B_\ell^a z_\ell
+\sum_{j=1}^{m} C_j^a z'_j)^{C_r^a}}
{\prod _{\bar a \in J_r^-}(\sum_{\ell=1}^k  -B_\ell^{\bar a} z_\ell
-\sum_{j=1}^{m} C_j^{\bar a} z'_j )^{-C_r^{\bar a}}
}\right)^{\frac{1}{\Delta}}, 1 \leq r \leq m.
$$
on d\'efinit l'application $h,$
$$
h: \bC^{k+m}\setminus \{0\} \rightarrow (\bC^{\times})^{k+m}, \leqno(4.1.5)$$
$$(z,z') \rightarrow (\psi_1(z,z'), \cdots ,
\psi_k(z,z'), \phi_1(z,z'), \cdots ,\phi_m(z,z')).$$
On se souvient ici que les ensembles
$I^+_q \subset \{1,2, \cdots, k\}$ (resp. $I^-_q$),
$J^+_r \subset \{1,2, \cdots, m\}$ (resp. $J^-_r$)
ont \'et\'e introduits dans la D\'efinition 3 du \S 2.1.

En vue de la propri\'et\'e (2.1.3), la fonction rationelle
$\psi_q(z,z')^\Delta$
(resp. $\phi_r(z,z')^\Delta$) est de poids homog\`ene z\'ero par
rapport aux variables $(z ,z')$ et donc on peut la regarder
comme l'application de $\bC P^{k+m-1}$ au lieu de celle de $\bC ^{k+m}.$
L'application $(4.1.5)$ n'est qu'une application inverse de l'application
de Gauss logarithmique;
$$ D_{s,s'}\cap (\bC^\times)^{k+m}
\rightarrow \bC P^{k+m-1},$$
$$ (s,s') \rightarrow \bigl(
s_1\frac{\partial}{\partial s_1}\Delta_f(s,s'): \cdots :
s_k\frac{\partial}{\partial s_k}\Delta_f(s,s'):
s_1'\frac{\partial}{\partial s_1'}\Delta_f(s,s'): \cdots:
s_m'\frac{\partial}{\partial s_m'}\Delta_f(s,s') \bigr).$$

Soit
$${\bar f}_\ell(x, \ba)=
a_{1,\ell}x^{\vec \alpha_{1,\ell}}+\cdots +a_{\tau_\ell,\ell}
x^{\vec \alpha_{\tau_\ell,\ell}}+a_{0,\ell}. \; 1 \leq \ell \leq k.$$

Pour la simplicit\'e on se sert de la notation
$\ba:=(a_{0,1}, \cdots, a_{\tau_k,k}) \in \bT^L.$
Nous regardons le cobord de Leray $\partial \gamma_\ba$ d'un cycle
$\gamma_\ba \in H_{N-k}(X_{\ba},\bZ)$ de l'IC $X_\ba=\{x \in \bT^N;
{\bar f}_1(x, \ba)= \cdots ={\bar f}_k(x, \ba)=0 \}.$

Alors on peut d\'efinir par
$$\Phi_{x^\bI, \gamma_{\ba}}^\zeta(a_{0,1}, \cdots, a_{\tau_k,k}):=$$
$$\int_{\partial \gamma_{\ba}} \prod_{\ell=1}^k
 {\bar f}_\ell(x, \ba)^{-\zeta_\ell-1}
x^{\bI +\b1} \frac{dx}{x^{\b1}},
$$
une $A-$ fonction hyperg\'eom\'etrique de Gel'fand-Zelevinski-Kapranov
\cite {GZK} associ\'ee aux polyn\^omes,
$$ f_\ell(x) = x^{\vec \alpha_{1,\ell}}+\cdots + x^{\vec \alpha_{\tau_\ell,\ell}},
\; 1 \leq \ell \leq k,$$
$$ x^{\bI}= x_1^{I_1}\cdots x_N^{I_N}, x^{\vec \alpha_{j,\ell}}=
x_1^{\alpha_{j,\ell,1}} \cdots x_N^{\alpha_{j,\ell,N}}.$$
On impose ici la condition de non-d\'eg\'en\'erescence
de la D\'efinition
~\ref{dfn1} pour l'IC ${\mathcal X}_s$ obtenue d'apr\`es la proc\'edure
d\'ecrite dans \S 1.1 pour n'importe quelle addition des variables
suppl\'ementaires $x'.$
Il est commode de consid\'erer la matrice suivante
analogue \`a la matrice $\sf L$ $(1.5.5),$
$${\sf M}(A):=
 \left [\begin {array}{cccccccccccccc}
1&1&\cdots&1&0&0&\cdots&0&0&\cdots&0&0&\cdots&0\\
0&0&\cdots&0&1&1&\cdots&1&0&\cdots&0&0&\cdots&0\\
0&0&\cdots&0&0&0&\cdots&0&1&\cdots&0&0&\cdots&0\\
\vdots&\vdots&\vdots&\vdots&\vdots&\vdots
&\vdots&\vdots&\vdots&\vdots&\vdots&\vdots&\vdots&\vdots\\
0&0&\cdots&0&0&0& \cdots&0&0&\cdots&1&1&\cdots&1\\
0&\alpha_{111}&\cdots&\alpha_{\tau_111}&0&\alpha_{121}& \cdots&
\alpha_{\tau_221}&0&\cdots&0&\alpha_{1k1}&\cdots&\alpha_{\tau_kk1}\\
\vdots&\vdots&\vdots&\vdots&\vdots&\vdots
&\vdots&\vdots&\vdots&\vdots&\vdots&\vdots&\vdots&\vdots\\
0&\alpha_{11N}&\cdots&\alpha_{\tau_11N}&0&\alpha_{12N}&\cdots&
\alpha_{\tau_22N} &0&
\cdots&0&\alpha_{1kN} &\cdots&\alpha_{\tau_kkN}\\ \end {array}\right ].
$$
Nous supposons d'ici bas toujours $rank({\sf M}(A))= k+N.$
En suite, on consid\`ere un r\'eseau $\Lambda \subset \bZ^L$
des $L-$vecteurs d\'efinis par le syst\`eme lin\'eaire
des \'equations ci-dessus:
$$ \sum^{\tau_q}_{i=0}b(j,q,\nu)=0,\;\; 1 \leq q \leq
k,$$ $$\sum_{q=1}^k\sum_{j=1}^{\tau_q}\alpha_{j q \ell}
b(j,q,\nu)=0,  1 \leq \ell \leq N.$$
Ici nous avons not\'e par $(b(0,1,\nu), \cdots,b(\tau_1,1,\nu),
b(0,2,\nu),\cdots, b(\tau_2,2,\nu),\cdots, b(\tau_k,k,\nu)),$
$1 \leq \nu \leq
m+k,$ une $\bZ$ base du $\Lambda.$

Pour un sous-ensemble $\bK \subset \{(0,1),\cdots, (k,\tau_k)\}$
tel que les colonnes ${\vec m_{j,q}}(A), (j,q) \in \bK$ de la matrice
${\sf M}(A)$
engendrent l'espace $\bR^{N+k}$ sur $\bR$ et $|\bK|= N+k$
on d\'efinit d'apr\`es \cite{GZK} l'ensemble des indices
(d\'et\'ermin\'es par la m\'ethode
 de Frobenius g\'en\'eralis\'ee),
$$ \Pi((\zeta+\b1, \bI +\b1), \bK)=\{((\lambda(0,1,\nu), \cdots,
\lambda(\tau_1,1,\nu),
\cdots, \lambda(\tau_k,k,\nu))\}_{
1 \leq \nu \leq |det ( {\vec m_{j,q}}(A))_{(j,q) \in \bK}| },$$
qui satisfont le syst\`eme des \'equations suivantes,
$$ \sum^{\tau_\nu}
_{j=0}\lambda(j,q,\nu) +\zeta_q +1
=0,\;\; 1 \leq q \leq k,$$
$$\sum_{q=1}^k\sum_{j=1}^{\tau_q}
\alpha_{j q \ell}\lambda(j, q, \nu)
-(I_\ell+1)=0,  1 \leq \ell \leq N.$$  Soit $T$ une triangulation du poly\`edre
de Newton $\Delta(F(x,\b1,y)+1)$ pour $F(x,\b1,y)$ d\'efini comme
au $(1.2.1)$ selon la d\'efinition \cite{GZK}, 1.2. Pourtant on impose
la condition que $ \lambda(j,q,\nu) \in \bZ $ pour $ (j,q) \not \in \bK.$
Soient $\bK_1, \bK_2 \in T$ deux simplexes diff\'erents de la
triangulation $T.$ On suppose que pour $(\lambda(0,1,\nu_p), \cdots,
\lambda(k,\tau_k, \nu_p)) \in \Pi((\zeta+\b1, \bI +\b1 ), \bK_p),$
$\lambda(j,q,\nu_p) \in
\bZ$ pour $(j,q) \not \in \bK_p,$ $(p=1,2)$ avec $ 1 \leq \nu_p \leq |det (
  {\vec m_\rho}(A))_{\rho \in \bK_p}|.$ On introduit la condition de
$T-$ non-r\'esonance sur $(\zeta+\b1, \bI+\b1)$ $$
(\lambda(0,1,\nu_1), \cdots,
\lambda(k,\tau_k, \nu_1))\not \equiv
(\lambda(0,1,\nu_2), \cdots,
\lambda(k,\tau_k, \nu_2))
\; mod \;\Lambda,  \leqno(4.1.6)$$ pour n'importe quel paire
$ (\lambda(0,1,\nu_p), \cdots,
\lambda(k,\tau_k, \nu_p))\in \Pi((\zeta+\b1, \bI +\b1 ),
\bK_p),$ $p=1,2.$ Une adaptation du theorem 3 \cite{GZK} \`a notre situation
peut \^etre formul\'ee comme suit.

\begin{thm}
La fonction $\Phi_{x^\bI, \gamma_a}^\zeta(\ba)$
satisfait le syst\`eme des \'equations comme suit.

1)$$\bigl(\sum_{j=0}^{\tau_q}a_{{ji}}\frac{\partial}{\partial a_{ji}}
+\zeta_q +1 \bigr)\Phi_{x^\bI, \gamma_a}^\zeta(\ba)=0, 1 \leq q \leq k,$$
$$
\bigl(\sum_{1 \leq q \leq k, 1 \leq j \leq \tau_q}
\alpha_{jq1}a_{{jq}}\frac{\partial}{\partial a_{jq}}
-(I_1+1) \bigr)\Phi_{x^\bI, \gamma_a}^\zeta(\ba)=
\cdots =\bigl(\sum_{1 \leq q \leq k, 1 \leq j \leq \tau_q}
\alpha_{jqN}a_{{jq}}\frac{\partial}{\partial a_{jq}}
-(I_N+1) \bigr)\Phi_{x^\bI, \gamma_a}^\zeta(\ba)= 0,$$
$$ \bigl(\prod_{\{(j,q);b(j,q,\nu)>0\}}
(\frac{\partial}{\partial a_{jq}})^{b(j,q,\nu)}
- \prod_{\{(j,q);b(j,q,\nu)<0\}}
(\frac{\partial}{\partial a_{jq}})^{-b(j,q,\nu)}\bigr)
\Phi_{x^\bI, \gamma_a}^\zeta(\ba)= 0,\; 1 \leq \nu \leq L-(k+N).$$

2) La dimension des solutions au point g\'en\'erique de $\ba \in
\bT^L$ du syst\`eme ci-dessus est \'egal \`a $(N+k)! vol_{N+k}
\Delta (F(x,\b1,y)+1)$ si la condition $(4.1.6)$ de
T-non-r\'esonance est satisfaite. \label{thm411}
\end{thm}

Dans la suite nous rangeons les param\`etres  $\ba
=(a_{0,1}, \cdots, a_{\tau_k,k})$ par leur ordre d'apparition
et d\'efinissons les param\`etres index\'es de nouveau
$a_{1}=a_{1,1}, \cdots,$ $a_{\tau_1}=a_{\tau_1,1},$ $a_{\tau_1+1}=
a_{0,1}, \cdots,$ $a_{L-1}=a_{\tau_k,k},$ $a_{L}=a_{0,k}.$
Nous introduisons la notation analogue du $(1.5.3),$
$$\Xi (A):= ^t(\log\; X_1, \cdots, \log\; X_N,  \log\; a_1, \cdots,
log\; a_L,\log\; U_1, \cdots,
log\; U_k). \leqno(4.1.7)$$
$$ \log\; T_1 = <\vec \alpha_{1,1}, \log\;X> + log\; a_1 + \log\; U_1,$$
$$\vdots$$
$$ \log\; T_{\tau_1} = <\vec \alpha_{1,\tau_1}, \log\;X> + log\; a_{\tau_1}
+ \log\; U_1,$$
$$\vdots$$
$$ \log\; T_L = log\; a_L  + \log\; U_k.$$

On consid\`ere l'\'equation
$$ {\sf L}(A)\cdot Log\; \Xi(A) = {\sf L}\cdot Log\; \Xi,$$
ici la matrice ${\sf L}(A)$ est construite comme suit.
Les colonnes $\vec \ell_i(A) = \vec w_i, 1 \leq i \leq N$
avec les vecteurs $\vec w_i$ d\'efinis comme la colonne de la matrice
$\sf L$ dans (1.5.5). Pour les colonnes du num\'ero $N+1$ \`a $N+L$
$(\vec \ell_{N+1}(A), \cdots, \vec \ell_{N+L}(A))= id_L.$
Les colonnes $\vec \ell_{N+L+j}(A)=
^t(\overbrace{0,\cdots, 0,0,}^{\tau_1+\cdots +\tau_{j-1}+j-1}
,\overbrace{1,1, \cdots, 1}^
{\tau_{j}+1}, 0, \cdots,0), 1 \leq j \leq k.$
Autrement dit, la matrice ${\sf L}(A)$ est obtenue apr\`es l'insertion
d'une matrice $id_L$ dans la matrice transpos\'ee $^t{\sf M}(A)$
entre la $k-$i\`eme et la $(k+1)-$i\`eme colonne quitte \`a des
permutations n\'ecessaires des colonnes apr\`es l'insertion.

 \begin{prop} Il existe un cycle $\gamma_a$ tel que
l'\'egalit\'e suivante soit vraie pour une int\'egrale d\'efinie dans
$(4.1.1),$ $$ \Phi_{x^\bI, \gamma_\ba}^\zeta(\ba) = B_\bI^\zeta(\ba)
I_{x^\bI, \gamma}^\zeta(s(\ba), s'(\ba)), \leqno(4.1.8)$$ ici $$
s_\ell(\ba)= \prod_{j=1}^L a_j^{\sigma_{j, N+\ell}}, \;\;, 1\leq \ell
\leq k,$$ $$ s_\rho'(\ba)= \prod_{j=1}^L a_j^{\sigma_{j, N+k+\rho}},
\;\;, 1\leq \rho \leq m,$$ $$B_\bI^\zeta(\ba)= \prod_{\ell=1}^N
(\prod_{j=1}^L a_j^{\sigma_{j,\ell}})^{I_\ell+1 }
\prod_{\nu=1}^k (\prod_{j=1}^L a_j^{\sigma_{j,N+k+m+\nu}})^{\zeta_\nu+1}.$$
Les exposants $\sigma_{j,\ell}$ sont d\'efinis par la relation suivante
$${\sf L}^{-1}\cdot {\sf L}(A)
= \left [\begin {array}{ccccccccc}
1&\cdots&0&\sigma_{1,1}&\cdots&\sigma_{L,1}&0&&\\
\vdots&\ddots&\vdots&\vdots&\vdots&\vdots&\vdots&\vdots&\vdots\\
0&\cdots&1&\sigma_{1,N}& \cdots& \sigma_{L,N}&0&&\\
0&\cdots&0&\sigma_{1,N+1}& \cdots& \sigma_{L,N+1}&0&&\\
\vdots&\ddots&\vdots&\vdots&\vdots&\vdots&\vdots&\vdots&\vdots\\
0&\cdots&0&\sigma_{1,N+k+m}& \cdots& \sigma_{L,N+k+m}&0&&\\
0&\cdots&0&\sigma_{1,N+k+m+1}& \cdots& \sigma_{L,N+k+m+1}&1&\cdots&0\\
\vdots&\ddots&\vdots&\vdots&\vdots&\vdots&\vdots&\vdots&\vdots\\
0&\cdots&0&\sigma_{1,L}& \cdots& \sigma_{L,L}&0&\cdots&1\\
\end {array}\right ].
\leqno(4.1.7)$$
La transition du cycle $\gamma(a)$ \`a
$\gamma$ est g\'er\'ee par les transformations,
$$ X_i = (\prod_{j=1}^L a_j^{\sigma_{j,i}})^{-1}
 \cdot x_i.$$
\label{prop42}
\end{prop}
{\bf D\'emonstration}
Il suffit de remarquer la relation
$$ x^{\bI+\b1} y^{\zeta +\b1} \frac{dx}{x^\b1}\wedge \frac{dy}{y^\b1}
=  B_\bI^\zeta(\ba) X^{\bI+\b1} U^{\zeta +\b1} \frac{dX}{X^\b1}\wedge
\frac{dU}{U^\b1}.$$
{\bf C.Q.F.D.}
\par On peut donc conclure que (au moins localement
sur une carte $a_j \not = 0$ pour $j \in I, |I|= k+m$)
$A-$fonction hyperg\'eom\'etrique
du type de GZK $(4.1.1)$ est exprim\'ee au moyen d'une int\'egrale-fibre
annul\'ee par un syst\`eme de Horn. On peut trouver un pareil \'enonc\'e
chez \cite{Kap1}.
\begin{cor}
Dans le cas o\`u  $dim\;\bD(\Sigma) =1,$ c.\`a.d. $m=0$
au (1.5.7), la condition de $T-$non r\'esonance  $(4.1.6)$
est \'equivalente \`a la condition, 
$$  x^{\bI +\b1} y^{\zeta +\b1} \frac{dx}{x^\b1}\wedge \frac{dy}{y^\b1} 
\in Gr^W_{N+k-1} H^{N+k-1}(Z_{F(x,0,y))}).$$
\label{cor43}
\end{cor}
\begin{cor}
La dimension de l'espace des solutions du syst\`eme du
Th\'eor\`eme ~\ref{thm411} au point g\'en\'erique est \'egal \`a
 $ |\chi ( Z_{F(x, \b1, y )} )|$ si la condition de $T-$
 non-r\'esonance $(4.1.6)$ est satisfaite. \end{cor} {\bf
D\'emonstration} Il faut consid\'erer l'envelope convexe des vecteurs
qui correspondent aux vertexes des polyn\^omes $ y_1(f_1(x)+1)+$
$\cdots$ $+y_k(f_k(x)+1).$ C'est \`a dire $(\vec \alpha_{1,1},1, 0,
\cdots,0),$$\cdots,$$(\vec \alpha_{\tau_1,1}, 1, 0,\cdots,0),$
$(\vec \alpha_{1,2},0,1, 0, \cdots,0),$ $\cdots,$ $(\vec \alpha_{\tau_k,k},0,
\cdots,0,1) \in \bZ^{N+k}.$ Ils se trouvent sur l'hyperplan $\zeta_1 +
\cdots +\zeta_k=1.$ Il est donc possible de regarder le volume
$(N+k-1)$ dimensionel $(N+k-1)! vol_{N+k-1}(\Delta(F(x,\b1,y))$ qui est \'egal
au $(N+k)! vol_{N+k}(\Delta(F(x,\b1,y)+1).$
Le caract\'eristique d'Euler est exprim\'e comme suit
$$ |\chi(Z_{F(x,\b1,y)})|= (N+k-1)!
vol_{N+k-1}\bigl(\Delta(F(x,\b1,y)
)\bigr),$$
d'apr\`es Khovanski-Oka \cite{Kh1}, \cite{Oka}. {\bf C.Q.F.D.}

Le corollaire ci-dessus avec la Proposition ~\ref{prop42}
nous donne
\begin{cor}
La dimension du syst\`eme $(4.1.4)_1,(4.1.4)_4 $
est \'egal \`a
$ |\chi(Z_{F(x,\b1,y)})|$ si la condition de $T-$non-r\'esonance
$(4.1.6)$ est satisfaite. \end{cor}

\begin{remark} (Compactification de Terasoma)
{\em Dans {\bf 3.2} nous avons regard\'e la relation entre deux
syst\`emes diff\'erents de multiplier $L-(N+k)$ termes par
nouvelles variables $(x_1', \cdots, x_m')$ et $({\tilde x}_1',
\cdots, {\tilde x}_m').$ Nous consid\'erons un ensemble des
triangulations $ {\mathcal J}=$ $\{$ la fa\c{c}on de choisir
$(N+k)-$termes de l'expression $y_1f_1(x)+ \cdots +y_kf_k(x)$
$\}.$ Alors \'evidemment un $(N+k)-$simplexe, $T_i$ $i \in
{\mathcal J}$ correspond au chaque choix des $L-(N+k)$ termes
auxquels les nouvelles variables sont multipli\'ees.

Par exemple dans $(1.6.2)$ les termes de num\'eros $1,2,4,5,7$
sont libres de ces variables, pourtant dans $(1.6.2)$ les termes
de num\'eros $1,2,4,7,8$ les sont. Les deux tores param\'etr\'es
par les variables $(s_1', \cdots, s_m') \in \bT^m= S_{i_1}'$ et
$({\tilde s}_1', \cdots, {\tilde s}_m') \in \bT^m = S_{i_2}'$ sont
ainsi index\'es par $ i_1, i_2 \in {\mathcal J}.$ Terasoma
\cite{Ter4} compactifie l'espace $\coprod_{\ba \in \bT^L} (X_\ba,
\ba)$ au moyen d'un espace produit $\prod_{i \in {\mathcal J}}
\bP_\Pi \times \prod_{i \in {\mathcal J}}S_{i}'. $ Ici $\Pi$
d\'enote le diagramme de Newton $\Delta(y_1f_1(x)+ \cdots
+y_kf_k(x))$ et $\bP_\Pi$ la vari\'et\'e torique y associ\'ee. Il
colle deux chartes $\bP_\Pi \times S_{i_1}'$ et $\bP_\Pi \times
S_{i_2}'$ dans leur compactification (not\'ee par $ \tilde
{\mathcal { H}_J}$) au moyen de la transition mentionn\'ee dans la
Proposition ~\ref{prop33} quitte \`a modifier $(x_1', \cdots,
x_m')$ et $({\tilde x}_1', \cdots, {\tilde x}_m')$ dans $(s_1',
\cdots, s_m')$ et $({\tilde s}_1', \cdots, {\tilde s}_m').$}
\label{remark41}
\end{remark}


{\bf 4.2 Application \`a la sym\'etrie de miroir
pour les espaces projectifs}

Comme chez Givental' \cite{Giv}
Nous regardons l'intersection compl\`ete d\'efinie par $r+1$ fonctions:
$$ X_s := \{(x_0, \cdots x_n) \in {\bT}^{n+1};
f_0(x) + s=0, f_1(x)+1=0, \cdots, f_r(x)+1=0\}.
\leqno(4.2.1)$$
o\`u
$$f_0(x)= x_0 x_1 \cdots x_n, f_i(x)=x_{\ell_1+\cdots+\ell_{i-1}}+ \cdots
+ x_{\ell_1+\cdots+
\ell_{i-1}+\ell_i-1}\; 1 \leq i \leq r. \leqno(4.2.2)$$
Ici on comprend que $\ell_0=0.$
On consid\`ere la fonction de phase comme au $ (1.1.5)',$
$$ y_0(f_0(x) +s)+ \sum_{i=1}^ry_i(f_i(x)+1)=T_0+ T_1 + \cdots + T_{n+r+3}.$$
La transform\'ee de Mellin de l'int\'egrale fibre associ\'ee
\`a $X_s$ ci-dessus peut \^etre calcul\'e par la matrice; $\sf L,$
$$ Log\; T= {\sf L}\cdot Log\; \Xi,$$
ici la notation est comme au \S 1.
$${\sf L}= \left [\begin {array}{ccccccccccccc}
1&1&\overbrace{\cdots}^{\ell_1-3}&1&1&\cdots&1&0&1&0&0&\cdots&0\\
0&0&\cdots                       &0&0&\cdots&0&1&1&0&0&\cdots&0\\
1&0&\cdots                       &0&0&\cdots&0&0&0&1&0&\cdots&0\\
0&1&\cdots                       &0&0&\cdots&0&0&0&1&0&\cdots&0\\
\vdots&\vdots&\ddots&\vdots&\vdots&\cdots&\vdots&\vdots&\vdots&\vdots&\vdots
&\cdots&\vdots\\
0&0&\cdots&1&0&\cdots&0&0&0&1&0&\cdots&0\\
0&0&\cdots&0&0&\cdots&0&0&0&1&0&\cdots&0\\
0&0&\cdots&0&1&\cdots&0&0&0&0&1&\cdots&0\\
\vdots&\vdots&\vdots&\vdots&\vdots&\ddots&\vdots&\vdots&\vdots&\vdots&\vdots
&\cdots&\vdots\\
0&0&\cdots&0&0&\cdots&1&0&0&0&0&\cdots&0\\
\end {array}\right ].\leqno(4.2.3)$$
Alors la transform\'ee de Mellin de l'int\'egrale fibre associ\'ee
\`a $X_s$ ci-dessus peut \^etre calcul\'e par la matrice; $(\sf L)^{-1}$,
qui se calcule facilement;
$${\sf L}^{-1}= \left [ \begin {array}{cccccccccccccc}
0&0&1&0&\overbrace{\cdots}^{\ell_1-3}&0&-1&0&0&\overbrace{\cdots}^{\ell_2-3}
&0&0&\overbrace{\cdots}^{n+r-\ell_1-\ell_2 -2}&0\\
0&0&0&1&  \cdots                     &0&-1&0&0&\cdots&0&0&\cdots&0\\
\vdots&\vdots&\vdots&\vdots&\cdots&\vdots&\vdots&\vdots&\vdots&\cdots&\vdots
&\vdots&\cdots&\vdots\\
0&0&0&0&  \cdots                     &1&-1&0&0&\cdots&0&0&\cdots&0\\
0&0&0&0&  \cdots                     &0& 0&1&0&\cdots&0&-1&\cdots&0\\
0&0&0&0&  \cdots                     &0& 0&0&1&\cdots&0&-1&\cdots&0\\
\vdots&\vdots&\vdots&\vdots&\cdots&\vdots&\vdots&\vdots&\vdots&\cdots&\vdots
&\vdots&\cdots&\vdots\\
0&0&0&0&  \cdots                     &0& 0&0&0&\cdots&1&-1&\cdots&0\\
\vdots&\vdots&\vdots&\vdots&\cdots&\vdots&\vdots&\vdots&\vdots&\cdots&\vdots
&\vdots&\cdots&\vdots\\
-1&1&1&1&  \cdots                     &1& -\ell_1&1&1&\cdots&1&-\ell_2
&\cdots&-\ell_r\\
1&0&-1&-1&  \cdots                     &-1& \ell_1&-1&-1&\cdots&-1&\ell_2
&\cdots& \ell_r\\
0&0&0&0&  \cdots                     &0& 1&0&0&\cdots&0&0&\cdots&0\\
0&0&0&0&  \cdots                     &0& 0&0&0&\cdots&0&1&\cdots&0\\
\vdots&\vdots&\vdots&\vdots&\cdots&\vdots&\vdots&\vdots&\vdots&\cdots&\vdots
&\vdots&\cdots&\vdots\\
0&0&0&0&  \cdots                     &0& 0&0&0&\cdots&0&0&\cdots&1\\
\end {array}\right ]. \leqno(4.2.4)$$

En r\'esultat nous avons la transofrm\'ee de Mellin
$M_{\bI,\gamma}^\zeta(z)$ de l'int\'egrale
p\'eriode $I_{x^\bI,\gamma}^\zeta(s)$ pour l'IC
$(4.2.1)$ admet une
expression comme suit \`a des fonctions p\'eriodiques pr\`es,
$$M_{\bI,\gamma}^\zeta(z)
=\prod_{i=0}^n \Gamma(z -\zeta_0 +I_i)
  \prod_{j=1}^r \Gamma( -\sum_{i=\ell_1+\cdots+\ell_{j-1}}^{
\ell_1+\cdots+\ell_j-1}I_i -\ell_j z+ \ell_j \zeta_0 +\zeta_j+1),$$
en particulier
$$ M_{0,\gamma}^0(z)= \frac{\Gamma(z)^{n+1}}
{\prod_{\nu=1}^r \Gamma(\ell_\nu z)}, \leqno(4.2.5)$$
toujours \`a des fonctions p\'eriodiques pr\`es.
Cette formule a \'et\'e reclam\'ee dans \cite{Giv}.
On se sert ici des notations des \S2 et \S3.

Si on regarde les vecteurs colonnes de la matrice ${\sf L}^{-1}$
apr\`es en avoir d\'epourvu les $r+2$ vecteurs rayons d'en bas, on obtient
une s\'erie de vecteurs
$ (1,0,\cdots,0),$$ (0,1,0,\cdots,0),$
$(0,\cdots, 0,\rlap{\ ${}^{\ell_1\atop{\hbox{${}^{\vee}$}}}$}1, \cdots,0),$
$\overbrace{(-1,\cdots ,-1}^{\ell_1}, 0, \cdots,0),$
$\cdots,$
$ \overbrace{(0,\cdots, 0}^{\ell_1+\cdots+\ell_{i-1}-1},
1,0,\cdots,0),$
$ \overbrace{(0,\cdots, 0}^{\ell_1+\cdots+\ell_{i-1}},
1,0,\cdots,0),$ $\cdots,$
$ \overbrace{(0,\cdots, 0}^{\ell_1+\cdots+\ell_{i}},
1,0,\cdots,0),$
$ \overbrace{(0,\cdots, 0}^{\ell_1+\cdots+\ell_{i-1}},
\overbrace{-1,\cdots, -1}^{\ell_{i}},0,\cdots,0),$ $1 \leq i \leq r.$
Ces $r$ groupes des vecteurs donnent naissance \`a un \'eventail
$\Sigma$ qui d\'efinit lavari\'et\'e torique
 $\bP_\Sigma= \bP^{\ell_1} \times \bP^{\ell_2}
\times \cdots \times
\bP^{\ell_r}.$
En particulier, on peut constater que l'int\'egrale-fibre associ\'ee \`a
l'IC (4.2.2) avec $r=1$ coincide avec la quantum cohomologie
de la vari\'et\'e $\bP^n.$ Voir \cite{Tan02}.


{\bf 4.3 L'hypoth\`ese de Berglund, Candelas et des autres.}

Dans l'article \cite{Ber} les auteurs ont propos\'e une
hypoth\`ese naturelle sur l'int\'egrale-fibre associ\'ee \`a la
vari\'et\'e de Calabi-Yau $X_s$ qui est une sym\'etrie de miroir
par rapport \`a la vari\'et\' e de Calabi-Yau g\'en\'erique
multi-quasihomog\`ene $Y$ de codimension $\ell$
d\'efinie dans le produit des espaces
projectifs quasihomog\`enes $\bP^{(\tau_1)}_{(g_1^{(1)},\ldots,
g_{\tau_1+1}^{(1)})} \times \ldots \times
\bP^{(\tau_k)}_{(g_1^{(k)},\ldots, g_{\tau_k+1}^{(k)})}$ dont la
codimension est $\ell$. Dans cette section, nous verifions leur
hypoth\`ese pour le cas $\ell =k.$

Tout d\'abord nous consid\'erons le syst\`eme suivant des
\'equations d\'efinies dans $({\bC^\times})^n,$

\begin{flushleft}
$$ f_1(x) = x^{\vec v^{(1)}_1} + \cdots + x^{\vec
v^{(1)}_{\tau_1}}, \leqno(4.3.1)$$
$$  f_2(x) =\prod_{j \in I^{(1)}} x_j  +1,$$
$$ \vdots,$$
$$ f_{2i-1}(x) = x^{\vec v^{(i)}_1} + \cdots + x^{\vec
v^{(i)}_{\tau_i}},$$
$$  f_{2i}(x) = \prod_{j \in I^{(i)}} x_j  +1,$$
$$ \vdots,$$
$$ f_{2k-1}(x) = x^{\vec v^{(k)}_1} + \cdots + x^{\vec
v^{(k)}_{\tau_k}},$$
$$  f_{2k}(x) = \prod_{j \in I^{(k)}} x_j+1.$$
\end{flushleft}
Ici, $I^{(j)}, 1 \leq j \leq k$ sont les ensembles d'indices
compl\'ementaires l'un \`a l'autre tels que $\cup_{q \in
[1,k]}I^{(q)}=\{1, \cdots, n \}$ et $I^{(q)}\cap I^{(q')}  =
\emptyset $ si $q \not = q'.$ Dor\'enavant nous nous serverons des
notations $\tilde \tau_\nu := |I^{(\nu)}|$ et $b^q :=
\sum_{\nu=1}^q \tau_\nu.$ On suppose que
$$ \sum_{\nu=1}^k \tau_\nu =b^k = n.$$
 L'\'equation
$f_{2j-1}(x)$ (resp.$f_{2j}(x)$) est donn\'ee par les mon\^omes
avec les puissances $\vec v^{(j)}_1,\cdots, \vec v^{(j)}_{\tau_j}
$ $\in \bZ^n$ (resp. $\vec v^{(j)}_{\tau_j+1} $ $\in \bZ^n$)
telles que pour le vecteur poids $\vec g^{(q)}=
(\overbrace{0,\cdots,0}^{b^{q-1}}, g_1^{(q)}, \cdots ,
g_{\tau_q}^{(q)}, \overbrace{0,\cdots,0}^{n-b^{q}})$ $\in
\bZ^{n}_{\geq 0}, $  $1 \leq q \leq k$ la relation ci-dessous soit
vraie,
$$ Q^{(q)}_j:=<\vec v^{(j)}_1, \vec g^{(q)}> = \cdots =
<\vec v^{(j)}_{\tau_j}, \vec g^{(q)}>=
<\vec v^{(j)}_{\tau_{j+1}},\vec g^{(q)}>
 \;\;\;, 1 \leq j \leq k, \leqno(4.3.2)$$
en sorte que $ g_{\tau_q+1}^{(q)}=Q^{(q)}_q.$
 Cela signifie que
le point $\vec v^{(j)}_{\tau_j+1}$ appartient \`a l'hyperplan de
dimension $\tau_j-1$  engendr\'e par $\vec v^{(j)}_1, \vec
v^{(j)}_1,\cdots, \vec v^{(j)}_{\tau_j}.$
 L'\'egalit\'e suivante  est demand\'ee
si on suppose que $X_s$ soit une vari\'et\'e de  Calabi-Yau:
$$ \sum_{j=1}^kQ^{(q)}_j = \sum_{i=1}^{\tau_q}g_i^{(q)}
= <\vec g^{(q)},(\overbrace{0,\cdots,0}^{b^{q-1}},
\overbrace{1,\cdots,1}^{\tau_{q}},\overbrace{0,\cdots,0}^{n-b^{q}})>,
1 \leq q \leq k. \leqno(4.3.3)$$ On ne suppose pas, pourtant,
cette condition dans l'argument suivant.

Du surcro\^it, on impose que  pour chaque \'el\'ement $ \lambda
\in Aut (X_j),$ du groupe de l'automorphisme de l'hypersurface
$X_j =\{x \in \bT^{n}; f_{2j-1}(x)=0 \}$ la relation
$(\lambda_{\ast}f_{2j-1})(x)= \lambda_{\ast}(x^{\vec
v^{(j)}_{\tau_j+1}})$ soit valable. Si on applique la technique de
Cayley au syst\`eme ci-dessus, on obtient un polyn\^ome
$$ F(x,s,y) = \sum_{j=1}^{k} y_{2j-1}(f_{2j-1}(x)+s_j) +
\sum_{j=1}^{k} y_{2j}f_{2j}(x), \leqno(4.3.4)$$ avec $L= n+3k$
termes. Pour faciliter la manipulation avec la matrice $\sf L$
ainsi construite pour $F(x,s,y)$ nous introduisons la notation
 $a^\nu:=\tau_1 + \cdots +\tau_{\nu}+3\nu = b^\nu+ 3\nu.$
En particulier, le $a^{i-1}+1-$\`eme terme de $F(x,s,y)$
correspond \`a $y_{2i-1}x^{\vec v^{(i)}_1}$ et le $a^{i}-3=
(a^{i-1}+\tau_i)-$\`eme terme - $y_{2i-1}x^{\vec
v^{(i)}_{\tau_i}}.$ Le $(a^{i}-1)-$\`eme terme - $y_{2i} \prod_{j
\in I^{(i)}} x_j.$Le $(a^{i}-2)-$\`eme terme - $y_{2i}s_i.$ Le
$a^{i}-$\`eme terme - $y_{2i}$.

On consid\`ere l'\'equation pour $\Xi = (\xi_1, \cdots, \xi_L),$
$$ ^t{\sf L}\cdot {\Xi} = ^t(1,\cdots,1, z_1, \cdots
,z_k),\leqno(4.3.5)$$ \'equivalente \`a la relation, $${\Xi}
=^t{\sf L}^{-1} \cdot ^t(1,\cdots,1, z_1, \cdots ,z_k).
\leqno(4.3.6)$$ Ainsi on obtient d'expressions lin\'eaires
$(\xi_1(z), \cdots, \xi_L(z))
=({\mathcal L}_1(0,z,0), \cdots, {\mathcal L}_L(0,z,0) $
pour $\bi  =0, \zeta =0$ en termes des fonctions lin\'eaires
d\'efinies dans la Proposition 2.1. Si on re\'ecrit la
matrice $\sf L$ en tenant compte des $(4.3.1), (4.3.4)$ on a le vecteur rayon
comme suit:
$$ \bar
v_{q}^{(\nu)}= (v_{q,1}^{(\nu)}, v_{q,2}^{(\nu)},\cdots,
v_{q,\tau_1}^{(\nu)}, \cdots, v_{q,n}^{(\nu)},\overbrace {0,
\cdots,0}^{2\nu-2},1, \overbrace{0, \cdots, 0}^{3k-2\nu
+1}),\leqno(4.3.7)$$ $1 \leq \nu \leq k, 1 \leq q \leq \tau_\nu.$
En se servant des notations du $(3.1.13), $ $\S 3$ pour les
vecteurs poids de la matrice ${\sf L}^{-1},$ on peut d\'eduire les
syst\`emes suivants pour chaque $\nu$ et $q$ fix\'es.
$$ \leqno(4.3.8)\begin{array}{ccc} v_{q,1}^{(\nu)} w_1^{(a^\nu)} +
v_{q,2}^{(\nu)}w_2^{(a^\nu)} + \cdots+ v_{q,n}^{(\nu)}
w_{n}^{(a^\nu)}&=&-1\;\; pour\;\; a^{(\nu-1)}+1 \leq q \leq
a^{(\nu-1)}+\tau_\nu,
\\ &=&0 \;\;\rm{ailleurs}.
\end{array} $$
$$
v_{1,j}^{(1)} p_q^{1} + v_{2,j}^{(1)}p_q^{2} + \cdots+
v_{\tau_1,j}^{(1)} p_q^{(\tau_1)} + v_{1,j}^{(2)}p_q^{(\tau_1+4)}
+ \cdots+ v_{\tau_2,j}^{(2)} p_{q}^{(\tau_1+\tau_2+3)}+ \cdots +
v_{\tau_k,j}^{(k)} p_q^{(a^k-3)} \leqno(4.3.9) $$
$$\begin{array}{ccc}
&=&-1\;\; pour
\;\;j \in I^{(q)}\\ &=&0.  \;\; pour j \not \in I^{(q)}
\end{array} $$
Il faut remarquer que le syst\`eme $(4.3.8)$ pour $\nu$ fix\'e
(resp. $(4.3.9)$ pour un $q$ fix\'e) consiste de $n-$\'equations
lin\'eaires ind\'ependantes par rapport aux inconnus $ \{
w_1^{(a^{\nu})}, \cdots, w_n^{(a^{\nu})} \},$ (resp. $ \{
p_q^{(j)};j \in I_\Lambda \}.$ Ici nous avons not\'e par
$I_\Lambda$ l'ensemble des indices $ \{1, \cdots, L\} \setminus
\cup_{\nu=1}^k\{a^{\nu}-2,a^{\nu}-1,a^{\nu} \}.$ On voit  encore
des autres  conditions n\'ec\'essaires sur $w_j^{(a^{\nu-1})}$, $$
\sum_{j \in I^{(\nu)}} w_j^{(a^{\nu}-1)}=1 , \leqno(4.3.10) $$
$$ \sum_{j \in
I^{(\nu)}}w_j^{(a^{\nu})}=-1 , \leqno(4.3.11) $$ puisque le produit
de la $(a^{\nu}-1)-$i\`eme colonne de ${\sf L}^{-1}$ avec  le
$(a^{\nu}-1)-$i\`eme rayon  de ${\sf L}$ est \'egal \`a $1.$ D'autre
part $n+2\nu-$i\`eme colonne de ${\sf L}$ avec  le  $j-$i\`eme
rayon ($1 \leq j \leq n+2$) de ${\sf L}^{-1}$ est \'egal \`a $0$
qui signifie $w_j^{(a^{\nu}-1)} +w_j^{(a^{\nu})} =0.$
En plus on voit que
$$
\sum_{j \in I^{(\nu)}} w_j^{(a^{\nu'}-1)}= \sum_{j \in
I^{(\nu)}}w_j^{(a^{\nu'})}=0,$$
pour $\nu' \not = \nu.$
On d\'efinit une $(n \times k)-$ matrice $V^\Lambda$ comme suit:
$$ V^\Lambda := ( ^tv_{\tau_1+1}^{(1)}, \cdots,
^tv_{\tau_k+1}^{(k)}), \leqno(4.3.12)$$ o\`u
$^tv_{\tau_q+1}^{(q)}$ est un $n-$ vecteur qui correspond au
$supp(f_{2\nu}) \setminus \{0\}.$ En vu de la
quasihomog\'en\'eit\'e $(4.3.2)$ on peut d\'efinir une $k\times
k$ matrice comme suivant:
$$  \hat{Q} :=
\left (\begin {array}{ccc}
Q^{(1)}_1&\cdots&Q^{(1)}_k\\
\vdots&\vdots&\vdots\\
Q^{(k)}_1&\cdots&Q^{(k)}_k\\
\end {array}\right )
= \left (\begin {array}{c}
\vec g^{(1)}\\
\vdots\\
\vec g^{(k)}\\
\end {array}\right )
\cdot V^\Lambda. \leqno(4.3.13)$$
Pour la simplicit\'e de la formulation, nous nous servons d'une matrice
diagonale,
$$ G=diag(g^{(1)}_1,\cdots, g^{(1)}_{\tau_1}, g^{(2)}_1,
\cdots, g^{(2)}_{\tau_2}, \cdots, g^{(k)}_1,\cdots,
g^{(k)}_{\tau_k} ). \leqno(4.3.14)$$
Dans la suite,
nous allons consid\'erer le cas
avec $rang (\hat{Q})=k.$ Introduisons une matrice $n\times n$ :
$$  {\sf L}_\Lambda :=\left (\begin {array}{c}
\vec v_1^{(1)}\\\vec v_2^{(1)}\\
\vdots\\
\vec v_{\tau_1}^{(1)}\\
\vdots\\
 \vec v_{\tau_k}^{(k)}\\
\end {array}\right ). \leqno(4.3.15)$$ Pour la commodit\'e nous nous servons d'une $n
\times k $ matrice,
$$\bQ := {\sf L}_\Lambda (\;^t\vec g^{(1)}, \cdots\;^t\vec g^{(k)}) = {\hat Q} \cdot \left (\begin {array}{c}
\vec g^{(1)}\\
\vdots\\
\vec g^{(k)}\\
\end {array}\right ) \cdot G^{-1}=$$
$$
= \left (\begin {array}{cccccccccc} Q^{(1)}_1&\overbrace{ \cdots
}^{\tau_1-2} & Q^{(1)}_1 & Q^{(1)}_2& \overbrace{\cdots
}^{\tau_2-2}& Q^{(1)}_2 & \cdots &  Q^{(1)}_k&\overbrace{ \cdots
}^{\tau_k-2}& Q^{(1)}_k\\
Q^{(2)}_1&\overbrace{ \cdots }^{\tau_1-2} & Q^{(2)}_1 & Q^{(2)}_2&
\overbrace{\cdots }^{\tau_2-2}& Q^{(2)}_2 & \cdots &
Q^{(2)}_k&\overbrace{ \cdots
}^{\tau_k-2}& Q^{(2)}_k\\
\vdots&\vdots&\vdots &\vdots&\vdots&\vdots
&\vdots&\vdots&\vdots &\vdots \\
Q^{(k)}_1&\overbrace{ \cdots }^{\tau_1-2} & Q^{(k)}_1 & Q^{(k)}_2&
\overbrace{\cdots }^{\tau_2-2}& Q^{(k)}_2 & \cdots &
Q^{(k)}_k&\overbrace{ \cdots
}^{\tau_k-2}& Q^{(k)}_k\\
\end {array}\right ).$$
Dans cette situation, on
d\'eduit du $(4.3.9)$ le syst\`eme ci-dessous:
$$
\left (\begin {array}{c}
\vec P_1\\
\vec P_2\\
\vdots\\
\vec P_k\\
\end {array}\right )
\cdot {\sf L}_\Lambda =- \;^tV^\Lambda, \leqno(4.3.16)$$
avec $\vec P_q =(p_q^{(1)},p_q^{(2)}, \cdots, p_q^{(\tau_1)}, \cdots,
p_q^{(a^k-3)} ).$
Si on note par $\check{\sf L}_j$ le $j-$\`eme vecteur colonne de la matrice
${\sf L}_\Lambda,$ on d\'eduit dir\'ectement du $(4.3.16)$
l'\'equation:
$$ \langle  z_1 \vec P_1 + z_2 \vec P_1 + \cdots + z_k
\vec P_k, \check{\sf L}_j
\rangle+z_q=0\;\; {\rm{si}} \;\; j \in I^{(q)}.$$

On construit une matrice $^T\sf L$ que l'on obtient de la matrice
transpos\'ee  $^t\sf L$ apr\`es des permutations propres des colonnes en
telle sorte que chaque rayon de $^T \sf L$ corresponde au vertexe d'un
polyn\^ome
$$ ^TF(x,s,y)=\sum_{j=1}^{k} y_{2j-1}(^Tf_{2j-1}(x)+s_j) +
\sum_{j=1}^{k} y_{2j}\;^Tf_{2j}(x), \leqno(4.3.4)^T$$
pour les polyn\^omes,
 \begin{flushleft}
$$ ^Tf_{2q-1}(x) = x^{^T\vec v^{(q)}_1} + \cdots + x^{^T\vec
   v^{(q)}_{\tilde \tau_q}}, \leqno (4.3.1)^T$$
$$ ^Tf_{2q}(x) =\prod_{\ell \in ^TI^{(q)}} x_\ell  +1,\;\; 1 \leq q
\leq k,\; 1 \leq q \leq k.  $$
\end{flushleft}

D'une facon analogue aux $(4.3.12) - (4.3.15)$ on peut d\'efinir le
syst\`eme des poids $$ (^T\vec g^{(1)}, \cdots , ^T\vec g^{(k)}), $$
$$^T{Q}_j^{(q)}= \langle ^T\vec v^{(j)}_r,  ^T\vec g^{(q)}\rangle,\; 1
\leq r \leq \tilde \tau_q, 1 \leq j,q \leq k.$$
Il est facile \`a voir que les \'equations du $\;^T(4.3.1)$
d\'efinissent une IC dans
$\bP^{(\tilde \tau_1)}_{(^Tg_1^{(1)},\ldots,
^Tg_{\tilde \tau_1+1}^{(1)})} \times
\ldots \times \bP^{(\tilde \tau_k)}_{(^Tg_1^{(k)},\ldots,
^Tg_{\tilde \tau_k+1}^{(k)})}$
avec $ ^Tg_{\tilde \tau_q+1}^{(q)}= ^T{Q}_q^{(q)}.$
Nous introduisons des
matrices analogues au cas $(4.3.1),$ $$^TV^\Lambda := (
^t(^Tv_{\tilde \tau_1+1}^{(1)}), \cdots, ^t(^Tv_{\tilde
\tau_k+1}^{(k)})), \leqno(4.3.12)^T$$ o\`u
$^Tv_{\tau_q+1}^{(q)}$ est un $n-$ vecteur qui correspond au
$supp(^Tf_{2q}) \setminus \{0\}.$
$$  \;^T{\hat Q} :=
\left (\begin {array}{ccc}
^TQ^{(1)}_1&\cdots&^TQ^{(1)}_k\\
\vdots&\vdots&\vdots\\
^TQ^{(k)}_1&\cdots&^TQ^{(k)}_k\\
\end {array}\right )
= \left (\begin {array}{c}
^T\vec g^{(1)}\\
\vdots\\
^T \vec g^{(k)}\\
\end {array}\right )
\cdot ^TV^\Lambda. \leqno (4.3.13)^T$$
$$ ^TG=diag
(^Tg^{(1)}_1,\cdots, ^Tg^{(1)}_{\tilde \tau_1}, ^Tg^{(2)}_1, \cdots,
^Tg^{(2)}_{\tilde \tau_2}, \cdots, ^Tg^{(k)}_1,\cdots,
\;^Tg^{(k)}_{\tilde \tau_k}
).   \leqno(4.3.14)^T$$
$$  ^T{\sf L}_\Lambda :=\left (\begin {array}{c}
^T\vec v_1^{(1)}\\
^T\vec v_2^{(1)}\\
\vdots\\
^T\vec v_{\tilde \tau_1}^{(1)}\\
\vdots\\
^T\vec v_{\tilde \tau_k}^{(k)}\\
\end {array}\right ).  \leqno(4.3.15)^T $$ On remarque ici que par la construction
$^T{\sf L}_\Lambda = ^t{\sf L}_\Lambda.$ Enfin on introduit les
fonctions lin\'eaires $^T\Xi := (^T\xi_1(z), \cdots,^T\xi_L(z)),$
d\'efinies par la relation,
$$^T{\Xi}
= \;^t({^T\sf L})^{-1} \cdot ^t(1,\cdots,1, z_1, \cdots ,z_k).\leqno(4.3.6)^T$$
$$^T\bQ :=
\;^T{\sf L}_\Lambda (\;^t(^T \vec g^{(1)}), \cdots\;
^t(^T \vec g^{(k)})).$$

Nous pouvons introduire l'ensemble des indices $^TI_\Lambda$
analogue au $I_\Lambda.$
Nous avons les indices du $I_\Lambda$ $1=j_1 <\cdots <j_n =L-4$
et celles du  $^TI_\Lambda,$ $1=i_1 <\cdots <i_n =L-4.$
Dans cette situation, nous consid\'erons deux conditions sur
les fonctions $(4.3.1),$
$$
\left (\begin {array}{c}
^T\xi_{i_1}(z)\\
\vdots\\
^T\xi_{i_n}(z)\\
\end {array}\right )
= G \cdot V^\Lambda \cdot \left (\begin {array}{c}
^T\xi^{(1)}(z)\\
\vdots\\
^T \xi^{(k)}(z)\\
\end {array}\right ), \leqno(4.3.17)$$
pour certaines fonctions lin\'eaires $(^T\xi^{(1)}(z), \cdots,
^T\xi^{(k)}(z))$ et
$$
\left (\begin {array}{c}
\xi_{j_1}(z)\\
\vdots\\
\xi_{j_n}(z)\\
\end {array}\right )
=^TG \cdot ^TV^\Lambda \cdot \left (\begin {array}{c}
\xi^{(1)}(z)\\
\vdots\\ \xi^{(k)}(z)\\
\end {array}\right ), \leqno(4.3.17)'$$
pour des fonctions lin\'eaires $(\xi^{(1)}(z), \cdots,
\xi^{(k)}(z)).$ En plus nous nous servons de la condition
$$   \;^TV^\Lambda = V^\Lambda, \leqno(4.3.17)''$$
apr\`es une permutation propre des variables lors de la
construction du polyn\^ome $\;^TF(x,s,y)$ par $(4.3.1)^T.$
On remarque que sous la condition $(4.3.17)''$
l'\'egalit\'e suivante a lieu entre les
nombres des termes pr\'esents dans $(4.3.1)$ et $(4.3.1)^T$
$$ \tau_\nu = \tilde \tau_\nu,\;\; 1 \leq \nu \leq k ,$$
car
$$ \overbrace{(1, \cdots, 1)}^n \cdot V^\Lambda =
(\tilde \tau_1, \cdots, \tilde \tau_k),$$
$$ \overbrace{(1, \cdots, 1)}^n \cdot \;^TV^\Lambda =
(\tau_1, \cdots, \tau_k).$$

En plus nous imposons une condition suppl\'ementaire,
$$   rank(\hat Q)=rank(^T \hat Q)=k.  \leqno(4.3.17)'''$$
\par

\begin{thm}
La transform\'ee de Mellin $M_{0,\gamma}^0(z)$ de l'int\'egrale
p\'eriode $I_{x^0,\gamma}^0(s)$ pour l'IC  $(4.3.1)$ admet une
expression comme suit \`a des fonctions $\Delta-$ p\'eriodiques pr\`es, si
elle satisfait les conditions $(4.3.17),$$(4.3.17)',$$(4.3.17)'',$
$(4.3.17)'''.$
$$ M_{0,\gamma}^0(z(\xi))= \frac{\prod_{\nu=1}^k \prod_{j=1}^{\tau_\nu}
\Gamma(\;^Tg_j^{(\nu)}\xi^{(\nu)})}{
\prod_{\nu=1}^k \Gamma(\sum_{q=1}^k
\;^TQ_q^{(\nu)} \xi^{(q)})}.
\leqno(4.3.18)$$ Ici
$z(\xi) =(z_1(\xi), \cdots , z_k(\xi))$ est une $k-$tuple
des fonctions
lin\'eaires en variables $\xi =(\xi^{(1)}, \cdots , \xi^{(k)})$
d\'efinies par la relation $(4.3.6).$

D'une fa\c{c}on sym\'etrique $ M_{0,\;^T\gamma}^0(z(\xi))$
pour l'IC, $(4.3.1)^T$ admet une expression comme suit
\`a des fonctions $\Delta$ p\'eriodiques pr\`es,
$$
M_{0,\;^T\gamma}^0(z(\xi))=
\frac{\prod_{\nu=1}^k \prod_{j=1}^{\tau_\nu}
\Gamma(g_j^{(\nu)}\;^T\xi^{(\nu)})}
{\prod_{\nu=1}^k \Gamma(\sum_{q=1}^k
Q_q^{(\nu)} \;^T\xi^{(q)})}. \leqno(4.3.18)^T$$
Les $^T\xi^{(\nu)}$ sont d\'efinies par la relation $(4.3.6)^T.$
\label{thm43} \end{thm}
 Afin de faciliter l'argument de l'\'epreuve, nous formulons un
\'enonc\'e auxiliaire.
\begin{lem}
Sous les conditions du Th\'eor\`eme ~\ref{thm43},
la transform\'ee de Mellin $M_{0,\gamma}^0(z)$ de l'int\'egrale
p\'eriode $I_{x^0,\gamma}^0(s)$ pour l'IC  $(4.3.1)$ admet une
expression comme suit \`a des fonctions $\Delta$
p\'eriodiques pr\`es,
$$M_{0,\gamma}^{0}(z)  = \prod_{i=1}^k \Gamma(z_i)
\prod_{j \in I_\Lambda} \Gamma(\xi_j(z)).$$
\label{lem43}
\end{lem}
{\bf D\'emonstration}
Il faut donc d\'emontrer
que
les solutions sp\'eciales du syst\`eme $(4.3.5)$
satisfont
$$ \xi_{a^{\nu}-2}(z)
=\xi_{a^\nu-1}(z)=z_{\nu},\;\;\xi_{a^\nu}(z) = 1-z_\nu. \leqno(4.3.19)$$

Pour le voir, on remarque que
le syst\`eme ci-dessous est une  cons\'equence directe de la relation
${\sf L}{\sf L}^{-1}= id_L,$
$$ \sum_{i=1}^n v_{q,i}^{(\nu)} w_i^{(a^\nu)} +
q_{2\mu-1}^{(a^\nu)} =0, \;\; a^{(\mu-1)}+1 \leq q \leq a^{\mu-1}+\tau_\mu,
\;\; \mu\in [1,k].\leqno(4.3.20)$$
$$- \sum_{i=1}^n w_i^{(a^\nu)}= \sum_{i=1}^n w_i^{(a^{\nu}-1)}= 1.$$
La derni\`ere se d\'eduit des $(4.3.10)$ et $(4.3.11).$
Soit $\nu(j) \in [1,k]$
tel que $j$ appartienne au $I^{(\nu(j))},$ alors on a
$$ \sum_{i \in I_\Lambda} v_{i,j}^{(\nu(j))} p_q^{(i)} +
p_q^{(a^{\nu(j)-1})} =0,\leqno(4.3.21)$$
pour $q \in [1,k].$

On en d\'eduit les relations suivantes.
$$q_{2\nu-1}^{(a^\nu)}=q_{2\nu}^{(a^\nu)}=1,
\leqno(4.3.22)$$
$$q_{2\nu-1}^{(a^{\nu'})}=q_{2\nu}^{(a^{\nu'})}=0,\;\;
\rm{pour}\; \nu' \not = \nu.
\leqno(4.3.23)$$
L'\'egalit\'e $q_{2\nu}^{(a^{\nu'})}=0$ se d\'eduit du fait que
le produit de la
$a^{\nu-1}-$i\`eme colonne de ${\sf L}^{-1}$ avec  le
$a^{\nu'}-$i\`eme rayon  de ${\sf L}$ est \'egal \`a $0.$
D'autre part, le produit de la
$(n+2\nu)-$i\`eme colonne de ${\sf L}$ avec  le
$a^{\nu'}-$i\`eme rayon  de ${\sf L}^{-1}$ est \'egal \`a
$q_{2\nu-1}^{(a^{\nu'})}+q_{2\nu}^{(a^{\nu'})}=0.$
Cela d\'emontre $(4.3.23).$
Le produit de la
$a^{\nu}-$i\`eme colonne de ${\sf L}^{-1}$ avec  le
$a^{\nu'}-$i\`eme rayon  de ${\sf L}$ est \'egal \`a $
q_{2\nu-1}^{(a^{\nu'})} +p_{\nu}^{(a^{\nu'})}=0.$
D'o\`u s'entra\^ine l'\'egalit\'e $p_{\nu}^{(a^{\nu'})}=0$
pour $\nu \not = \nu'.$
Le produit de la
$(a^{\nu}-1)-$i\`eme colonne de ${\sf L}^{-1}$ avec  le
$(a^{\nu'}-2)-$i\`eme rayon  de ${\sf L}$ est \'egal \`a $
q_{2\nu'-1}^{(a^{\nu}-1)} +p_{\nu'}^{(a^{\nu}-1)}=0.$
Par consequence, on a
$$ p_{\nu'}^{(a^{\nu}-1)}=0 \;\;\rm{pour}\;\; \nu \not = \nu', \leqno(4.3.24)$$
$$ p_{\nu}^{(a^{\nu}-1)}=-q_{2\nu-1}^{(a^{\nu}-1)}  =1,\;\;
p_{\nu}^{(a^{\nu})}=-1,$$
car $ p_{\nu}^{(a^{\nu})}+p_{\nu}^{(a^{\nu}-1)} =0. $
En plus on voit que
$$ q_{r}^{(a^{\nu} -1)}=0,\;\; r \not = 2\nu-1. \leqno(4.3.25)$$

En somme, gr\^ace aux  $(4.3.20), (4.3.22), (4.3.24),$
$$ \xi_{a^\nu}(z)= \sum_{i=1}^n w_{i}^{(a^{\nu})} + \sum_{r=1}^{2k}
q_{r}^{(a^{\nu})} + \sum_{q=1}^k  p_{q}^{(a^{\nu})} z_q
= -1 +1 + 1 -z_\nu =  1 -z_\nu.$$
D'autre part $(4.3.20), (4.3.25), (4.3.24),$
entra\^inent
$$  \xi_{a^\nu-1}(z)= 1+0-1 + z_\nu= z_\nu.$$

Quant'\`a la fonction $  \xi_{a^\nu-2}(z)$ il est facile de voir que tout
les \'el\'ements de la $ (a^\nu-2)-$i\`eme  colonne de ${\sf L}^{-1}$
consiste des z\'eros  sauf le $(n+2k+\nu)-$i\`eme \'el\'ement qui
est \'egal \`a $1.$

On a donc l'\'egalit\'e
$$M_{0,\gamma}^{0}(z)  = \prod_{i=1}^k\Gamma(z_i)^2 \Gamma(1-z_i)
\prod_{j \in I_\Lambda} \Gamma(\xi_j(z)) = \prod_{i=1}^k \frac{\pi
}{sin\;\pi z_i} \Gamma(z_i) \prod_{j \in I_\Lambda}
\Gamma(\xi_j(z)).$$ {\bf C.Q.F.D.}

{\bf D\'emonstration du Th\'eor\`eme ~\ref{thm43}}
On d\'emontre la relation suivante,
$$
\left (\begin {array}{c}
1-z_1\\
\vdots\\
1-z_k\\
\end {array}\right )
=^T{\hat Q}\left (\begin {array}{c}
\xi^{(1)}(z)\\
\vdots\\ \xi^{(k)}(z)\\
\end {array}\right ). $$
En fait $$^t(^T \bQ)\cdot ^t(^T{\sf L}_\Lambda)^{-1}\cdot \;^TV^\Lambda
\left (\begin {array}{c}
1-z_1\\
\vdots\\
1-z_k\\
\end {array}\right )
=  \left (\begin {array}{c}
\;^T\vec g^{(1)}\\
\vdots\\ \;^T\vec g^{(k)}\\
\end {array} \right )
\cdot ^TV^\Lambda
\left (\begin {array}{c}
1-z_1\\
\vdots\\
1-z_k\\
\end {array}\right ) =
\;^T{\hat Q}\left (\begin {array}{c}
1-z_1\\
\vdots\\
1-z_k\\
\end {array}\right )$$
qui se d\'eduit du $ (4.3.6)^T$ et du $ (4.3.13)^T.$ D'autre part
gr\^ace \`a la condition $(4.3.17)''$ on a
$$ ^t(^T \bQ)\cdot ^t(^T{\sf L}_\Lambda)^{-1}\cdot \;^TV^\Lambda
\left (\begin {array}{c}
1-z_1\\
\vdots\\
1-z_k\\
\end {array}\right )=
 ^t(^T \bQ)\cdot {\sf L}_\Lambda^{-1}\cdot V^\Lambda
\left (\begin {array}{c}
1-z_1\\
\vdots\\
1-z_k\\
\end {array}\right )=
^T \hat Q \left(\begin {array}{c}
\;^T\vec g^{(1)}\\
\vdots\\ \;^T\vec g^{(k)}\\
\end {array} \right ) \cdot {^TG}^{-1}
\left(\begin {array}{c}
\;^T\xi_{j_1}(z)\\
\vdots\\ \;^T\xi_{j_n}(z)\\
\end {array} \right ).
$$
La derni\`ere \'egalit\'e s'explique par la relation
$^t(^T \bQ) = ^T \hat Q \left(\begin {array}{c}
\;^T\vec g^{(1)}\\
\vdots\\ \;^T\vec g^{(k)}\\
\end {array} \right ) \cdot {^TG}^{-1} $
et par une consequence directe du $(4.3.16),$
$$ {\sf L}_\Lambda^{-1}\cdot V^\Lambda
\left (\begin {array}{c}
1-z_1\\
\vdots\\
1-z_k\\
\end {array}\right )= \left(\begin {array}{c}
\;^T\xi_{j_1}(z)\\
\vdots\\ \;^T\xi_{j_n}(z)\\
\end {array} \right ). $$
 Le cha\^ine des \'egalit\'es ci-dessus est vrai
en vue du $(4.3.16)$ et du $(4.3.17),$
$$
^T \hat Q \left(\begin {array}{c}
\;^T\vec g^{(1)}\\
\vdots\\ \;^T\vec g^{(k)}\\
\end {array} \right ) \cdot {^TG}^{-1}
 \left(\begin {array}{c}
\;^T\xi_{j_1}(z)\\
\vdots\\ \;^T\xi_{j_n}(z)\\
\end {array} \right )=
^T \hat Q \left(\begin {array}{c}
\vec g^{(1)}\\
\vdots\\ \vec g^{(k)}\\
\end {array} \right ) \cdot {G}^{-1}
 \left(\begin {array}{c}
\;^T\xi_{j_1}(z)\\
\vdots\\ \;^T\xi_{j_n}(z)\\
\end {array} \right )
$$

$$=
^T \hat Q \left(\begin {array}{c}
\vec g^{(1)}\\
\vdots\\
\vec g^{(k)}\\
\end {array} \right )
\cdot G^{-1} G V^{\Lambda} \left(\begin {array}{c}
\;^T\xi^{(1)}(z)\\
\vdots\\ \;^T\xi^{(k)}(z)\\
\end {array} \right )
= \;^T \hat Q \cdot \hat Q
 \left(\begin {array}{c}
\;^T\xi^{(1)}(z)\\
\vdots\\
\;^T\xi^{(k)}(z)\\
\end {array}\right ).$$ L'\'egalit\'e,
$$ \left(\begin {array}{c}
\vec g^{(1)}\\
\vdots\\
\vec g^{(k)}\\
\end {array} \right )
\cdot G^{-1}=\left (\begin {array}{cccccccccc} 1&\overbrace{
\cdots }^{\tau_1-2} & 1 & 0& \overbrace{\cdots }^{\tau_2-2}& 0 &
\cdots & 0&\overbrace{ \cdots}^{\tau_k-2}& 0\\
0&\overbrace{ \cdots }^{\tau_1-2} & 0 & 1& \overbrace{\cdots
}^{\tau_2-2}& 1& \cdots & 0&\overbrace{ \cdots
}^{\tau_k-2}& 0\\
\vdots&\vdots&\vdots &\vdots&\vdots&\vdots
&\vdots&\vdots&\vdots &\vdots \\
0 &\overbrace{ \cdots }^{\tau_1-2} &0 & 0&
\overbrace{\cdots}^{\tau_2-2}& 0& \cdots
& 1&\overbrace{ \cdots}^{\tau_k-2}& 1\\
\end {array}\right )=$$
$$ = \left(\begin {array}{c}
 \;^T \vec g^{(1)}\\
\vdots\\
\;^T \vec g^{(k)}\\
\end {array} \right ) \cdot \;^T G^{-1}
=\left (\begin {array}{cccccccccc} 1&\overbrace{ \cdots }^{\tilde
\tau_1-2} & 1 & 0& \overbrace{\cdots }^{\tilde \tau_2-2}& 0 &
\cdots & 0&\overbrace{ \cdots}^{\tilde \tau_k-2}& 0\\
0&\overbrace{ \cdots }^{\tilde \tau_1-2} & 0 & 1&
\overbrace{\cdots }^{\tilde \tau_2-2}& 1& \cdots & 0&\overbrace{
\cdots
}^{\tilde \tau_k-2}& 0\\
\vdots&\vdots&\vdots &\vdots&\vdots&\vdots
&\vdots&\vdots&\vdots &\vdots \\
0 &\overbrace{ \cdots }^{\tilde \tau_1-2} &0 & 0&
\overbrace{\cdots }^{\tilde \tau_2-2}& 0& \cdots & 1&\overbrace{
\cdots
}^{\tilde \tau_k-2}& 1\\
\end {array}\right ),$$
est une cons\'equence des relations $\tau_\nu =\tilde \tau_\nu, 1
\leq \nu \leq k $ mentionn\'ees comme corollaire du (4.3.17)".

Gr\^ace \`a  la condition $rank(^T \hat Q)=n$
celui-l\`a signifie
$$
\left (\begin {array}{c}
1-z_1\\
\vdots\\
1-z_k\\
\end {array}\right )
=
{\hat Q}
\left(\begin {array}{c}
\;^T\xi^{(1)}(z)\\
\vdots\\
\;^T\xi^{(k)}(z)\\
\end {array}\right ).
$$
D'une fa\c{c}on sym\'etrique, en supposant que  $rank(\hat Q)=n,$ on a
$$
\left (\begin {array}{c}
1-z_1\\
\vdots\\
1-z_k\\
\end {array}\right )
=
\;^T{\hat Q}
\left(\begin {array}{c}
\xi^{(1)}(z)\\
\vdots\\
\xi^{(k)}(z)\\
\end {array}\right ).
$$

On en d\'eduit directement la relation d\'esir\'ee,
$$
\prod_{\nu=1}^k \Gamma(z_\nu) =
\prod_{\nu=1}^k
\Gamma(1-\sum_{q=1}^k\;^TQ^{(\nu)}_q \xi^{(q)}(z))
=\frac{1}
{\prod_{\nu=1}^k \Gamma ( \sum_{q=1}^k\;^TQ^{(\nu)}_q \xi^{(q)}(z) )
sin \left (\;\pi \sum_{q=1}^k\;^TQ^{(\nu)}_q \xi^{(q)}(z) \right)}.$$

D'ailleurs le lemme ~\ref{lem43}
et la condition $(4.1.17)'$ nous donne
$$ \xi_{a^{\nu-1}+j}=
\;^Tg_j^{(\nu)}\xi^{(\nu)},\;\; j\in [1, \tau_\nu]$$
qui signifie
$$ \prod_{j \in I_\Lambda} \Gamma(\xi_j(z))
= \prod_{\nu=1}^k \prod_{j =1}^{\tau_\nu}
\Gamma(\;^Tg_j^{(\nu)}\xi^{(\nu)}).$$

La formule $(4.3.18)^T$ se d\'emontre d'une fa\c{c}on compl\`etement
parall\`ele.
{\bf C.Q.F.D.}

{\bf Conjecture}
Le pair $(G\cdot V^\Lambda, \;^TG\cdot \;^TV^\Lambda)$
correspond aux sym\'etries quantique ($Q_X$)et g\'eom\'etrique
(${\mathcal G}_X$) de l'IC $X, (4.3.1)$ pourtant le
pair  $(^TG\cdot\;^TV^\Lambda, \;G\cdot V^\Lambda)$
\`a celles-ci de l'IC  $Y,
(4.3.1)^T
$ au
sens de \cite{Ber}.

En fait, dans le cas $k=1, $ on peut facilement voir
la dualit\'e comme
suit existe
$$ Q_X = (\bZ_{\hat Q}; G),{\mathcal G}_X=(\bZ_{^T \hat Q};\;^T G), $$
$$ {\mathcal G}_Y = (\bZ_{\hat Q}; G),Q_Y=(\bZ_{^T \hat Q};\;^T G), $$
avec $|G|= \hat Q$ qui se d\'eduit de la condition de Calabi-Yau.

En liaison avec la sym\'etrie de miroir, on consid\`ere l'alg\`ebre 
structurel de l'IC de dimension $n-k$
not\'ee par $X$ et d\'efinie par $(4.3.1),$
$$ A_X:= \frac{\bC [x]}{(f_1+f_2-1, \cdots f_{2k-1}+f_{2k}-1) \bC 
[x]},$$
et la filtration naturelle sur lui,
$$ A^j_X := \{x^\alpha \in A_X; \sum_{q=1}^k \langle \alpha, \vec g^{(q)}
\rangle =j
 \},$$
avec le polyn\^ome de Poincar\'e,
$$ P_{A_X}(\lambda) = \sum_{j \in \bZ_{\geq 0}} dim (A_X^j) \lambda^j.$$
D'une mani\`ere analogue, on d\'efinit les notions
correspondantes de l'IC de la m\^eme dimension  
not\'ee par $Y$ et d\'efinie par $(4.3.1)^T,$
$$ A_Y:= \frac{\bC [x]}{(^Tf_1+ ^Tf_2-1, \cdots ^Tf_{2k-1}
+^Tf_{2k}-1) \bC 
[x]},$$
$$ A^j_Y := \{x^\alpha 
\in A_Y; \sum_{q=1}^k \langle \alpha, 
^T \vec g^{(q)} 
\rangle =j \},$$
$$ P_{A_Y}(\lambda) = \sum_{j \in \bZ_{\geq 0}} 
dim (A_Y^j) \lambda^j.$$
Alors le r\'esultat classique du \cite{Dolg} 
nous donne,
$$ P_{A_X}(\lambda) = \frac
{\prod_{q=1}^k \prod_{\nu=1}^k (1-\lambda^{Q^{(\nu)}_q})}
{\prod_{\nu=1}^k \prod_{j=1}^{\tau_\nu} (1-\lambda^{g^{(\nu)}_j})}, 
P_{A_Y}(\lambda) = \frac{\prod_{q=1}^k \prod_{\nu=1}^k 
(1-\lambda^{^TQ^{(\nu)}_q})}{\prod_{\nu=1}^k 
\prod_{j=1}^{\tau_\nu} (1-\lambda^{^Tg^{(\nu)}_j})}.
\leqno(4.3.26)$$

En suite on introduit les variables
$(t_1, \cdots, t_k) \in \bT^k$ telles que
$$ \prod^k_{\nu=1} s_\nu^{^TQ^{(\nu)}_q} = t_q, \;\; 1 \leq q \leq k.$$
Gr\^ace \`a la condition $(4.3.17)''$, cett'\'equation est toujours resoluble par rapport
aux variables $s=s(t).$
On consid\`ere la transform\'ee de Mellin inverse de $M_{0,\gamma}^0(z(\xi))$
associ\'ee \`a l'IC (4.3.1),
$$ U_\alpha(s)= 
\int_{\check \Pi_\alpha}
\frac{\prod_{\nu=1}^k \prod_{j=1}^{\tau_\nu}
\Gamma(\;^Tg_j^{(\nu)}\xi^{(\nu)}(z))}{
\prod_{\nu=1}^k \Gamma(\sum_{q=1}^k
\;^TQ_q^{(\nu)} \xi^{(q)}(z))} s^{-\bz +\b1} d\bz,$$
o\`u ${\check \Pi}_\alpha \subset \bT^k$ est un cycle qui \'evite
les lieux singuliers de l'int\'egrand.
Il est facile de voir que
$$ U_\alpha(s(t))= det (^T \hat Q)^{-1}\int_{^T \hat Q^{-1}(\check \Pi_\alpha)}
\frac{\prod_{\nu=1}^k \prod_{j=1}^{\tau_\nu}
\Gamma(\;^Tg_j^{(\nu)}\xi^{(\nu)})}{
\prod_{\nu=1}^k \Gamma(\sum_{q=1}^k
\;^TQ_q^{(\nu)} \xi^{(q)})} t_1^{-\xi^{(1)}}\cdots t_k^{-\xi^{(k)}}
d\xi^{(1)} \wedge \cdots \wedge d\xi^{(k)}. $$
La transform\'ee de Mellin inverse $ U_\alpha(s(t))$ est annul\'ee 
par le syst\`eme suivant
d'\'equations diff\'erentielles,
$$  L_\nu (t_\nu, \vartheta_t)U_\alpha(s(t))=0  
, \;\; 1 \leq \nu \leq k,$$
o\`u
$$ L_\nu (t_\nu, \vartheta_t)=  
\left(\prod_{j=1}^{\tau_\nu}\prod_{r=0}^{^Tg_j^{(\nu)}-1}
( - ^Tg_j^{(\nu)} \vartheta_{t_\nu} +r) - t_\nu 
\prod_{\mu=1}^k \prod_{r=0}^{\;^TQ_\nu^{(\mu)}-1}(\sum_{q=1}^k
\;^TQ_q^{(\mu)} \vartheta_{t_q} -r) \right), \;\; 1 \leq \nu \leq k.$$
Nous notons par  $|\chi(\tilde X_\nu)|$ le d\'egr\'e
de l'op\'erateur $L_\nu (t, \vartheta_t):$
$|\chi(\tilde X_\nu)|=\sum_{\mu=1}^k\;^TQ_\nu^{(\mu)}=     
\sum_{j=1}^{\tau_\nu}\;^Tg_j^{(\nu)}= \;^Tg_{\tau_\nu+1}^{(\nu)}. $ 
On  d\'efinit
 la restriction de l'op\'erateur
$ L_{\nu}(t, \vartheta_t)$ au tore 
$\bT = \{t \in \bC^k; t_i=0, i \not = \nu\} \setminus \{
t_\nu=0\}$ comme suit, 
$$ {\tilde L}_{\nu}(t_\nu, \vartheta_{t_{\nu}})
:= \left( \prod_{j=1}^{\tau_\nu}\prod_{r=0}^{^Tg_j^{(\nu)}-1}
( - ^Tg_j^{(\nu)} \vartheta_{t_\nu} +r) - t_\nu 
\prod_{\mu=1}^k \prod_{r=0}^{\;^TQ_\nu^{(\mu)}-1}(
\;^TQ_\nu^{(\mu)} \vartheta_{t_\nu} -r) \right).$$

Sur l'espace de dimension $|\chi(\tilde X_q)|$ 
des solutions annul\'ees par 
${\tilde L}_q(t_q, \vartheta_{t_q})$, on consid\`ere la 
monodromie $M_q^{(0)} \in GL(|\chi(\tilde X_q)|, \bC) $
(resp. 
$M_q^{(\infty)} \in GL(|\chi(\tilde X_q)|, \bC) $)
autour du point $t_q=0$ (resp. $t_q=\infty$).
Alors l'expression ci-dessus de l'op\'erateur  ${\tilde L}_{\nu}
(t_\nu, \vartheta_{t_{\nu}})$ entra\^ine,
$$ det (M_q^{(\infty)} - \lambda_q \cdot id_{|\chi(\tilde X_q)|})
= \prod_{\nu=1}^k 
(1-\lambda_q^{^TQ^{(\nu)}_q}),
det (M_q^{(0)} - \lambda_q \cdot id_{|\chi(\tilde X_q)|})
= \prod_{j=1}^{\tau_q} 
(1-\lambda_q^{^Tg^{(q)}_j}).$$
Par cons\'equence la fonction rationelle d\'efinie par
$$ M_X(\lambda_1, \cdots \lambda_k): =
\prod^k_{q=1} \frac
{det (M_q^{(\infty)} - \lambda_q \cdot id_{|\chi(\tilde X_q)|})}
{det (M_q^{(0)} - \lambda_q \cdot id_{|\chi(\tilde X_q)|})}$$
s'exprime
$$M_X(\lambda_1, \cdots \lambda_k)=
\frac{\prod_{q=1}^k \prod_{\nu=1}^k 
(1-\lambda_q^{^TQ^{(\nu)}_q})}{\prod_{q=1}^k 
\prod_{j=1}^{\tau_q} (1-\lambda_q^{^Tg^{(q)}_j})}. \leqno(4.3.27)$$
Pour la fonction rationelle 
$M_Y(\lambda_1, \cdots \lambda_k)$ d\'efinie d'une man\`ere analogue \`a 
la fonction $M_X(\lambda_1, \cdots \lambda_k),$ on a 
$$M_Y(\lambda_1, \cdots \lambda_k)=
\frac{\prod_{q=1}^k \prod_{\nu=1}^k 
(1-\lambda_q^{Q^{(\nu)}_q})}{\prod_{q=1}^k 
\prod_{j=1}^{\tau_q} (1-\lambda_q^{g^{(q)}_j})}.$$

Soit $\bar Y$ la compactification de l'IC $Y \subset \bT^{n}$
dans le produit des espaces projectifs quasihomog\`egenes
$\bP:=\bP^{(\tau_1)}_{(^Tg_1^{(1)},\ldots,
^Tg_{\tau_1+1}^{(1)})} \times
\ldots \times \bP^{(\tau_k)}_{(^Tg_1^{(k)},\ldots,
^Tg_{\tau_k+1}^{(k)})}.$ On peut d\'efinir le faisceau coh\'erent
${\mathcal O}_{\bP}(\zeta),$ 
$\zeta=(\zeta_1, \cdots, \zeta_k)$ dont la section sur un ouvert
$U_I=\{x \in \bC^n; x_i \not=0, i \in I   \}$
est donn\'ee par
$$ \Gamma(U_I, {\mathcal O}_{\bP}(\zeta)):= \{x^\alpha; \alpha=(\alpha_1,
\cdots,\alpha_n) \in \bZ^n, 
\alpha_i \geq 0 \;pour \; i \not \in I,\langle ^Tg^{(q)}, \alpha \rangle = 
\zeta_q, 1 \leq q \leq k \}.$$
On d\'efinit le faisceau coh\'erent
${\mathcal O}_{\bar Y}(\zeta) $ par la section sur l'ouvert
$U_I,$
$$ \Gamma(U_I, {\mathcal O}_{\bar Y}(\zeta)):= 
\{x^\alpha \in A_Y; \alpha= (\alpha_1,
\cdots,\alpha_n) \in \bZ^n, 
\alpha_i \geq 0 \;\;pour \;\;i \not \in I, \langle ^Tg^{(q)}, 
\alpha \rangle = 
\zeta_q, 1 \leq q \leq k \}.$$
Pour lui on introduit le caract\'eristique d'Euler,
$$ \chi({\mathcal O}_{\bar Y}(\zeta)):= 
\sum^{n-k}_{i=0}(-1)^i dim H^i({\mathcal O}_{\bar Y}(\zeta)).$$

D'apr\`es \cite{Dolg}, \cite{Gol} le polyn\^ome de Poincar\'e
du caract\'eristique d'Euler, 
$$ P{\mathcal O}_{\bar Y}(t_1, \cdots, t_k)
:= \sum_{\zeta \in (\bZ_{\geq 0})^k}\chi({\mathcal O}_{\bar Y}(\zeta))
t_1^{\zeta_1} \cdots t_k^{\zeta_k}$$ admet une expression comme suit,
$$ P{\mathcal O}_{\bar Y}(t_1, \cdots, t_k)= 
\frac{\prod_{\nu=1}^k \prod_{q=1}^k  
(1-t_\nu^{^TQ^{(\nu)}_q})}{\prod_{\nu=1}^k 
\prod_{j=1}^{\tau_\nu} (1-t_\nu^{^Tg^{(\nu)}_j})}. \leqno(4.3.28)$$
Pour le faisceau 
${\mathcal O}_{\bar X}(\zeta)$ d\'efini d'une fa\c{c}on analogue au 
${\mathcal O}_{\bar Y}(\zeta)$, on d\'efinit 
$$ P{\mathcal O}_{\bar X}(t_1, \cdots, t_k)
:= \sum_{\zeta \in (\bZ_{\geq 0})^k}\chi({\mathcal O}_{\bar X}(\zeta))
t_1^{\zeta_1} \cdots t_k^{\zeta_k}.$$ Pour celui-l\`a on a,
$$ P{\mathcal O}_{\bar X}(t_1, \cdots, t_k)= 
\frac{\prod_{\nu=1}^k \prod_{q=1}^k  
(1-t_\nu^{Q^{(\nu)}_q})}{\prod_{\nu=1}^k 
\prod_{j=1}^{\tau_\nu} (1-t_\nu^{g^{(\nu)}_j})}.$$
Si on compare $(4.3.26),$ $ (4.3.27)$
et $ (4.3.28)$ on obtient l'\'enonc\'e suivant.
\begin{thm}
Pour les IC  $X$ et $Y$ d\'efinies par $(4.3.1)$ et $(4.3.1)^T$
on a la relation suivante,
$$M_X(\lambda,\cdots, \lambda)=  P{\mathcal O}_{\bar Y}(\lambda,\cdots, 
\lambda)
= P_{A_Y}(\lambda),$$
$$M_Y(\lambda,\cdots, \lambda)=  P{\mathcal O}_{\bar X}(\lambda,\cdots, 
\lambda)
= P_{A_X}(\lambda).$$
\label{thm432}
\end{thm}

{\bf Exemple 4.3.1, La vari\'et\'e de Schimmrigk}
Comme exemple simple, mais peu trivial, on regarde l'exemple
suivant \'etudi\'e dans \cite{Ber},
$$ f_1(x)= \sum_{i=0}^3 x_i^3,
f_2(x)= x_1 x_2 x_3 +1,$$
$$ f_3(x)= \sum_{i=1}^3 x_i x_{i+3}^3,
 f_4(x)= x_0 x_4 x_5 x_6 +1.$$
Alors on a les matrices ci-dessous,

$${\sf L}=
\left (\begin {array}{ccccccccccccc} 3&0&0&0&0&0&0&1&0&0&0&0&0
\\\noalign{\medskip}0&3&0&0&0&0&0&1&0&0&0&0&0\\\noalign{\medskip}0&0&3
&0&0&0&0&1&0&0&0&0&0\\\noalign{\medskip}0&0&0&3&0&0&0&1&0&0&0&0&0
\\\noalign{\medskip}0&0&0&0&0&0&0&1&0&0&0&1&0\\\noalign{\medskip}0&1&1
&1&0&0&0&0&1&0&0&0&0\\\noalign{\medskip}0&0&0&0&0&0&0&0&1&0&0&0&0
\\\noalign{\medskip}0&1&0&0&3&0&0&0&0&1&0&0&0\\\noalign{\medskip}0&0&1
&0&0&3&0&0&0&1&0&0&0\\\noalign{\medskip}0&0&0&1&0&0&3&0&0&1&0&0&0
\\\noalign{\medskip}0&0&0&0&0&0&0&0&0&1&0&0&1\\\noalign{\medskip}1&0&0
&0&1&1&1&0&0&0&1&0&0\\\noalign{\medskip}0&0&0&0&0&0&0&0&0&0&1&0&0
\end {array}\right ).$$

$${\sf L}^{-1}=
\left (\begin {array}{ccccccccccccc} 1/3&-1/9&-1/9&-1/9&0&1/3&-1/3&0&0
&0&0&0&0\\\noalign{\medskip}0&2/9&-1/9&-1/9&0&1/3&-1/3&0&0&0&0&0&0
\\\noalign{\medskip}0&-1/9&2/9&-1/9&0&1/3&-1/3&0&0&0&0&0&0
\\\noalign{\medskip}0&-1/9&-1/9&2/9&0&1/3&-1/3&0&0&0&0&0&0
\\\noalign{\medskip}-1/9&-1/27&{\frac {2}{27}}&{\frac {2}{27}}&0&-1/9&
1/9&2/9&-1/9&-1/9&0&1/3&-1/3\\\noalign{\medskip}-1/9&{\frac {2}{27}}&-
1/27&{\frac {2}{27}}&0&-1/9&1/9&-1/9&2/9&-1/9&0&1/3&-1/3
\\\noalign{\medskip}-1/9&{\frac {2}{27}}&{\frac {2}{27}}&-1/27&0&-1/9&
1/9&-1/9&-1/9&2/9&0&1/3&-1/3\\\noalign{\medskip}0&1/3&1/3&1/3&0&-1&1&0
&0&0&0&0&0\\\noalign{\medskip}0&0&0&0&0&0&1&0&0&0&0&0&0
\\\noalign{\medskip}1/3&-1/9&-1/9&-1/9&0&0&0&1/3&1/3&1/3&0&-1&1
\\\noalign{\medskip}0&0&0&0&0&0&0&0&0&0&0&0&1\\\noalign{\medskip}0&-1/
3&-1/3&-1/3&1&1&-1&0&0&0&0&0&0\\\noalign{\medskip}-1/3&1/9&1/9&1/9&0&0
&0&-1/3&-1/3&-1/3&1&1&-1\end {array}\right ).$$

D'apr\`es la construction $(4.3.4)^T$
bas\'ee sur la matrice transpos\'ee
$\;^t(\sf L)$ il est facile de voir que
$ \;^Tf_\ell(x) = f_\ell(x), 1 \leq \ell \leq 4.$
La matrice $\;^T{\sf L}^{-1}= {\sf L}^{-1}$ nous donne
les fonctions lin\'eaires
$$\xi^{(1)} = -\frac{1}{3}(z_1-1) +\frac{1}{9}(z_2-1)
, \xi^{(2)} = -\frac{1}{3}(z_2-1).$$
On en d\'eduit facilement que pour $f(x)$ ainsi que $^Tf(x)$
on a la transform\'ee de Mellin de la p\'eriode,
$$ M_{0,\gamma}^{0}(z)=
\Gamma( -\frac{1}{3}(z_1-1) +\frac{1}{9}(z_2-1))^3
\Gamma(-\frac{1}{3}(z_2-1))^4 \Gamma(z_1)\Gamma(z_2)
= \frac{\Gamma(\xi^{(1)})^3\Gamma(\xi^{(2)})^4}{\Gamma(3\xi^{(1)}
+\xi^{(2)}) \Gamma(3\xi^{(2)}) },$$
\`a des fonctions $27-$p\'eriodiques pr\`es.

D'une mani\`ere analogue, on peut consid\'erer l'IC,
$$ f_1(x)= \sum_{i=0}^n x_i^n,
f_2(x)= x_1 x_2 \cdots x_n +1,$$
$$ f_3(x)= \sum_{i=1}^n x_i x_{i+n}^n,
 f_4(x)= x_0 x_{n+1} x_{n+2}\cdots x_{2n} +1,$$
dont l'int\'egrale fibre est exprim\'ee par
$$ M_{0,\gamma}^{0}(z)=
 \frac{\Gamma(\xi^{(1)})^n\Gamma(\xi^{(2)})^{n+1}}{\Gamma(n\xi^{(1)}
+\xi^{(2)}) \Gamma(n\xi^{(2)}) },$$
\`a des fonctions $n^n-$p\'eriodiques pr\`es.

Il faut remarquer que notre matrice $\sf L$ satisfait une condition
int\'eressante ci-dessous.

{\bf La condition de carr\'e magique}
Pour chaque $q   \in [1,k]$ fix\'e, on peut trouver une application
univalente $\sigma: b \in I_\Lambda  \rightarrow \{1,
\cdots,n \}$ telle que
$$ p_q^b= w_{\sigma(b)}^{a^q}, {\rm{pour\; tout}}\;b \in I_\Lambda.$$

Cette condition joue un r\^ole central
lors de l'interpretation de la dualit\'e bizarre d'Arnol'd sur
l'\'echange des nombres de Gabrielov et de Dolgachev
de point de vue de la sym\'etrie de miroir \cite{Kob}.

{\bf Exemple 4.3.2}
Nous  consid\'erons un exemple d'une hypersurface \'etudi\'e dans
\cite{Ber}.
Avec les notations ci-dessus, nous avons les donn\'ees suivantes.
$$ f_1(x)= x_1^7 + x_2^7x_4+x_3^7x_5+ x_4^3+x_5^3,\;\;
f_2=x_1x_2x_3x_4x_5.$$

$$
{\sf L}=\left [\begin {array}{cccccccc} 7&0&0&0&0&1&0&0\\\noalign{\medskip}0&7
&0&1&0&1&0&0\\\noalign{\medskip}0&0&7&0&1&1&0&0\\\noalign{\medskip}0&0
&0&3&0&1&0&0\\\noalign{\medskip}0&0&0&0&3&1&0&0\\\noalign{\medskip}0&0
&0&0&0&1&0&1\\\noalign{\medskip}1&1&1&1&1&0&1&0\\\noalign{\medskip}0&0
&0&0&0&0&1&0\end {array}\right ],$$

$$
{\sf L}^{-1}=
\left [\begin {array}{cccccccc} {\frac {6}{49}}&-1/49&-1/49&-{\frac {2
}{49}}&-{\frac {2}{49}}&0&1/7&-1/7\\\noalign{\medskip}-{\frac {2}{147}
}&{\frac {19}{147}}&-{\frac {2}{147}}&-{\frac {11}{147}}&-{\frac {4}{
147}}&0&2/21&-2/21\\\noalign{\medskip}-{\frac {2}{147}}&-{\frac {2}{
147}}&{\frac {19}{147}}&-{\frac {4}{147}}&-{\frac {11}{147}}&0&2/21&-2
/21\\\noalign{\medskip}-1/21&-1/21&-1/21&{\frac {5}{21}}&-2/21&0&1/3&-
1/3\\\noalign{\medskip}-1/21&-1/21&-1/21&-2/21&{\frac {5}{21}}&0&1/3&-
1/3\\\noalign{\medskip}1/7&1/7&1/7&2/7&2/7&0&-1&1\\\noalign{\medskip}0
&0&0&0&0&0&0&1\\\noalign{\medskip}-1/7&-1/7&-1/7&-2/7&-2/7&1&1&-1
\end {array}\right ]$$
On a donc la transform\'ee de Mellin de la
p\'eriode associ\'ee \`a  l'IC $\{x \in (\bC
^\times)^5; f_1(x)+s_1=0, f_2(x)+1=0\}$ exprim\'ee ci-dessous \`a des
$7-$fonctions p\'eriodiques pr\`es,
$$ M_{0,\gamma}^{0}(z)=
\Gamma( -\frac{1}{7}(z_1-1))^3
\Gamma(-\frac{2}{7}(z_1-1))^2 \Gamma(z_1) =
\frac{\Gamma(\xi^{(1)})^3\Gamma(2\xi^{(1)})^2}{\Gamma(7\xi^{(1)})}.$$

D'apr\`es la construction $(4.3.4)^T$ on voit que
$$ \;^Tf_1(x)= x_1^7 + x_2^7+x_3^7+ x_2 x_4^3+x_3x_5^3,\;\;
\;^Tf_2=x_1x_2x_3x_4x_5.$$

On a alors la transform\'ee de Mellin de la
p\'eriode associ\'ee \`a  l'IC $\{x \in (\bC
^\times)^5; \;^Tf_1(x)+s_1=0, \;^Tf_2(x)+1=0\}$
exprim\'ee ci-dessous \`a des
$147-$fonctions p\'eriodiques pr\`es,
$$ M_{0,\gamma}^{0}(z)=
\Gamma( -\frac{1}{7}(z_1-1))
\Gamma(-\frac{2}{21}(z_1-1))^2
\Gamma(-\frac{1}{3}(z_1-1))^2 \Gamma(z_1) =
\frac{\Gamma(3\xi^{(1)})\Gamma(2\xi^{(1)})^2
\Gamma(7\xi^{(1)})^2}{\Gamma(21\xi^{(1)})}.$$


\vspace{\fill}

%

\noindent

\begin{flushleft}
\begin{minipage}[t]{6.2cm}
  \begin{center}
{\footnotesize Indepent University of Moscow\\
Bol'shoj Vlasievskij pereulok 11,\\
 Moscow, 121002,\\
Russia\\
{\it E-mails}:  tanabe@mccme.ru, tanabe@mpim-bonn.mpg.de \\}
\end{center}
\end{minipage}
\end{flushleft}

\end{document}